\documentclass[11pt]{amsart}
\usepackage{fullpage}
\usepackage[utf8]{inputenc}
\usepackage[T1]{fontenc}
\usepackage{graphicx}
\usepackage{amsmath}
\usepackage{amsfonts}
\usepackage{graphicx}

\usepackage{color}

\usepackage{hyperref}

\newtheorem{theo}{Theorem}[section]
\newtheorem{cor}[theo]{Corollary}
\newtheorem{rem}[theo]{Remark}

\newtheorem{propo}[theo]{Proposition}
\newtheorem{lemme}[theo]{Lemma}
\newtheorem{defi}[theo]{Definition}
\newtheorem{ex}[theo]{Example}

\newtheorem{hyp}[theo]{Assumptions}
\newcommand{\E}{\mathbb{E}}
\newcommand{\R}{\mathbb{R}}
\newcommand{\PP}{\mathbb{P}}

\newcommand{\F}{\mathcal{F}}

\newcommand{\N}{\mathbb{N}}

\newcommand{\DM}{\mathcal{D}}
\newcommand{\PsiM}{\tilde{\Psi}^{(M)}}

\title[Full discretization and invariant laws of SPDEs]{Approximation of the invariant law of SPDE\tiny{s}\normalsize{: error
analysis using a Poisson equation for a full-discretization scheme}}
\author[C-E Br\'ehier]{Charles-Edouard Br\'ehier}
\address{Universit\'e Paris-Est, CERMICS (ENPC),  6-8-10 Avenue Blaise Pascal, Cit\'e Descartes ,  F-77455 Marne-la-Vall\'ee, France}
\address{INRIA Paris-Rocquencourt, Domaine de Voluceau - Rocquencourt, B.P. 105 - 78153 Le Chesnay, France}
\email{brehierc@cermics.enpc.fr,charles-edouard.brehier@inria.fr}
\author[M. Kopec]{Marie Kopec}
\address{ENS Rennes - IRMAR, Universit\'e Rennes 1\\Avenue Robert Schumann\\F-35170 Bruz\\France}
\email{marie.kopec@bretagne.ens-cachan.fr}
\date{}

\begin{document}

\begin{abstract}
We study the long-time behavior of fully discretized semilinear SPDEs with additive space-time white noise, which admit a unique invariant probability measure $\mu$. We show that the average of regular enough test functions with respect to the (possibly non unique) invariant laws of the approximations are close to the corresponding quantity for $\mu$.

More precisely, we analyze the rate of the convergence with respect to the different discretization parameters. Here we focus on the discretization in time thanks to a scheme of Euler type, and on a Finite Element discretization in space.

The results rely on the use of a Poisson equation, generalizing the approach of \cite{MattStuTre}; we obtain that the rates of convergence for the invariant laws are given by the weak order of the discretization on finite time intervals: order $1/2$ with respect to the time-step and order $1$ with respect to the mesh-size.
\end{abstract}

\maketitle

\section{Introduction}

In this article, we want to analyze in a quantitative way the effect of time and space discretization schemes on the knowledge of the unique invariant law of a semi-linear Stochastic PDE of parabolic type, written in the abstract form of \cite{DaP-Z1}
\begin{equation}\label{eqY}
\begin{gathered}
dX(t,x)=\bigl(AX(t,x)+F(X(t,x))\bigr)dt+dW(t), \quad 0<t\leq T,\\
X(0,x)=x.
\end{gathered}
\end{equation}
This process takes value in an infinite-dimensional Hilbert space $H$ - typically $H=L^{2}(0,1)$; $A$ is a negative, self-adjoint, unbounded linear operator on $H$, with a compact inverse - for instance, $A=\frac{\partial^2}{\partial \xi^2}$, given on the domain $H^2(0,1)\cap H_{0}^{1}(0,1)$ when homogeneous Dirichlet boundary conditions are applied. The coefficient $F$ is an operator on $H$, on which regularity conditions are assumed and given below.

Finally, $\bigl(W(t)\bigr)_{t\in[0,T]}$ is a cylindrical Wiener process on $H$, so that we have a space-time white noise in \eqref{eqY}.

More precise assumptions are given below. In a more general setting, one could for instance intend to remove boundedness of $F$, or include spatially correlated noise processes. However, we already consider the main technical and conceptual difficulties in this paper.

In our setting, it is known that the SPDE \eqref{eqY} admits a unique invariant probability measure $\mu$, and that convergence is exponentially fast. This result comes from the spatial non-degeneracy of the noise and from a dissipation relation satisfied by the drift part. Nevertheless, in general, no expression of $\mu$ is available for practical use; moreover, the support of this measure is an infinite-dimensional space.

The approximation of quantities like $\int_{H}\phi d\mu$ for bounded test functions $\phi$ is therefore complicated. The exponential convergence ensures that $\E\phi(X(t))$ tends to $\int_{H}\phi d\mu$ when $t$ tends to infinity, exponentially fast. However, using a Monte-Carlo method to compute $\E\phi(X(t))$ would require the ability to simulate the $H$-valued random variable $X(t)$, for a large value of $t$. Since it is not possible to have an exact simulation, we introduce two approximations:
\begin{itemize}
\item a discretization in time, in order to get an approximation of the law of the random variables $X(t)$ for different, fixed values of $t$, using a finite number of calculations; it done here with a semi-implicit Euler scheme;
\item a discretization in space, in order to replace the $H$-valued random variables with finite-dimensional ones. Here, it is performed with a finite element method.
\end{itemize}
The second approximation is specific to the case of SPDEs, while the first one has already been studied a lot in the case of SDEs.

Different techniques to control the error are available in the litterature. A first one is found in \cite{TalInv}, where an estimate of the weak error introduced by the numerical scheme is proved, holding for any value of the finite time $T$. The idea there is to expand the error thanks to the solution of the Kolmogorov equation associated with the diffusion, and to prove bounds on the spatial derivatives of this solution, with an exponential decrease with respect to the time variable.

This strategy has been generalized to SPDEs like \eqref{eqY} in \cite{B2}, when a semi-implicit Euler scheme is used. The main additional difficulty, when compared with the SDE case, is the need for tools introduced in \cite{deb}, for the estimate of the weak error at a fixed time $T$.

Using theses tools aims at proving that on a finite time interval the weak order of convergence is twice the strong one: in other words, laws at fixed times are approximated more accurately than the trajectories. These tools have also been used in \cite{WangGan} to treat the time-discretization in a slightly more-general setting, and in \cite{AndersLars} where discretization in space with a finite element method is studied. Basically, the two ingredients are the following:
\begin{itemize}
\item improved estimates on the derivatives of the solution of the Kolmogorov equations, with spatial regularization;
\item an integration by parts formula issued from Malliavin calculus, in order to transform some stochastic expressions with insufficient spatial regularity into more suitable ones.
\end{itemize}
These tools are fundamental to treat equations with nonlinear terms; they are used again in the present work. Notice that for linear equations a specific idea simplifies the proof - so that the second tool is not used - but can not be adapted for nonlinear parabolic equations like \eqref{eqY}: see \cite{deb-prin}, and \cite{deB-deb} where a stochastic Schr\" odinger equation is discretized.

Here, we are interested in another method for the approximation of the invariant measure: we want to follow the approach of \cite{MattStuTre}. There, the authors study the distance between time-averages along the realization of the numerical scheme of a test function $\phi$, ant its expected value with respect to the invariant law $\mu$. They introduce the solution $\Psi$ of the Poisson equation $\mathcal{L}\Psi=\phi-\int\phi d\mu$, where $\mathcal{L}$ is the infinitesimal generator of the SDE - the solvability of this elliptic or hypoelliptic PDE is ensured by ergodic properties. Then they show how to expand the error for various numerical methods, in a stochastic Taylor expansions fashion.

The use of a Poisson equation to prove convergence results of Law of Large Numbers type is classical, as explained in \cite{MattStuTre}. In the context of SPDEs, it has been used in \cite{B1} and \cite{Ce-F} for the study of the averaging principle for systems evolving with two separate time-scales.

Such a technique gives an approximation result for $\mu$, even if the numerical method is not ergodic, having possibly several invariant laws. In the SDE case, the study of ergodicity for time-discretizated processes has been the subject of \cite{HighMattStu}; there the author use general results on Markov chains, like the Harris Theorem. Up to our knowledge, no such study has been completed for SPDEs so far.

Our main result is the adaptation of the approach of \cite{MattStuTre} for SPDEs, with time and space approximation procedures: we essentially obtain the following result - a more precise statement is Theorem \ref{th}:
There exists a constant $C>0$ such that for any $\mathcal{C}^2_b(H)$ function $\phi$, any parameters $\tau\in(0,1)$ and $h\in(0,1)$, any time $N\geq 1$ and any initial condition $x\in H$
\begin{equation*}
\Big|\frac{1}{N}\sum_{m=0}^{N-1}\Bigl(\E\phi(X_{m}^{h})-\int_{H}\phi(z)\mu(dz)\Bigr)\Big|\leq C\Big(\tau^{1/2}+h+\frac{1}{N\tau}\Big).
\end{equation*}

The error is divided into three parts: the first (resp. the second) is due to the time (resp. space) discretization, while the last one goes to $0$ when time increases.

For the proof of this result, we need to study the Poisson equation associated with the SPDE \eqref{eqY}. More precisely, we work with Galerkin approximations, and show bounds that are independent of dimension. Moreover, the strategy described above for the study of weak approximation for infinite dimensional processes require some additional regularization properties for the solution of the Kolmogorov equation. Here, we need similar properties for the solutions of Poisson equations. When compared with \cite{AndersLars}, \cite{B2} or \cite{deb}, many error terms look the same: the method with Poisson equation does not really simplifies the proof, it just moves some technical problems to other places along the proof.
The main reason for studying only the additive noise case is the following. As soon as the equation is discretized either in time - \cite{deb} - or in space - \cite{AndersLars} - the possible diffusion coefficients must satisfy very strict conditions: they should be decomposed as the sum of a continuous affine function, and another function such that the second order derivative is controlled with respect to a very weak norm - namely, the norm associated with a negative power of the linear operator. Moreover, the treatment of such noise requires lengthier computations. We could do so here by adding our argument with the ones in \cite{AndersLars} and \cite{deb} but this would result only in hiding the main ideas of our work.

The paper is organized as follows: in Section \ref{sectass}, we precise the assumptions made on the coefficients of the equation, and we define the discretization method in Section \ref{sectdefnum}. In Section \ref{sectasymptu}, we study the asymptotic behavior of the solutions of the continuous and discrete time processes. In Section \ref{sectres}, we give the convergence results that we obtained. In Section \ref{sectdescproof} we show how the error is decomposed, and we present two essential tools: the Poisson equation, and an integration by parts from Malliavin calculus. Finally detailed proofs of the estimates are developed in Section \ref{sectdetproof}.


\section{Notations and assumptions}\label{sectass}

Let $\mathcal{D}\subset\mathbb{R}$ be a bounded, open interval; without restriction in the sequel we assume $\mathcal{D}=(0,1)$. Let $H=L^2(\mathcal{D})$, with norm and inner product denoted by $|.|_{H}$ and $\langle.,.\rangle_{H}$ or simply $|.|$ and $\langle.,.\rangle$.

We consider equations in the abstract form
\begin{equation}
\begin{gathered}
dX(t,x)=(AX(t,x)+F(X(t,x)))dt+dW(t)\\
X(0,x)=x.
\end{gathered}
\end{equation}

In the next paragraphs, we state the assumptions made on the coefficients $A$ and $F$ in \eqref{eqY}. We also recall basic facts on the cylindrical Wiener process $W$, and on the mild solution of the SPDE.

\subsection{Test functions}
To quantify the weak approximation, we use test functions - called \textit{admissible} - $\phi$ in the space $\mathcal{C}_{b}^{2}(H,\R)$ of functions from $H$ to $\R$ that are twice continuously differentiable, bounded, with first and second order bounded derivatives.

\begin{rem}\label{remident}
In the sequel, we often identify the first derivative $D\phi(x)\in \mathcal{L}(H,\R)$ with the gradient in the Hilbert space $H$, and the second derivative $D^2\phi(x)$ with a linear operator on $H$, via the formulae:
\begin{gather*}
\langle D\phi(x),h\rangle =D\phi(x).h \text{ for every }h\in H\\
\langle D^2\phi(x).h,k\rangle =D^2\phi(x).(h,k) \text{ for every }h,k\in H.
\end{gather*}
We then use the following notations, for an admissible test function $\phi$:
\begin{gather*}
\|\phi\|_{\infty}=\sup_{x\in H}|\phi(x)|_{H}\\
\|\phi\|_{1}=\sup_{x\in H}|D\phi(x)|_{H},\\
\|\phi\|_{2}=\sup_{x\in H}|D^2\phi(x)|_{\mathcal{L}(H)}.
\end{gather*}
\end{rem}

\subsection{Assumptions on the coefficients}

\subsubsection{The linear operator}

We denote by $\N=\left\{0,1,2,\ldots\right\}$ the set of nonnegative integers.

We suppose that the following properties are satisfied:
\begin{hyp}\label{hypB}
\begin{enumerate}
\item We assume that there exists a complete orthonormal system of elements of $H$ denoted by $(e_k)_{k\in \N}$, and a non-decreasing sequence of real positive numbers $(\lambda_k)_{k\in \N}$ such that:
$$Ae_{k}=-\lambda_{k}e_{k}\text{ for all } k\in\N.$$
\item The sequence $(\lambda_k)_{k\in\mathbb{N}}$ goes to $+\infty$ and
$$
\sum_{k=0}^{+\infty}\frac{1}{\lambda_{k}^{\alpha}}<+\infty\Leftrightarrow \alpha>1/2.$$
\end{enumerate}
\end{hyp}

The smallest eigenvalue of $-A$ is then $\lambda_0$.

\begin{ex}\label{exampleB}
We can choose $A=\frac{d^2}{dx^2}$, with the domain $H^2(0,1)\cap H_{0}^{1}(0,1)\subset L^2(0,1)$ - corresponding to homogeneous Dirichlet boundary conditions. In this case for any $k\in\N$ $\lambda_k=\pi^2 (k+1)^2$, and $e_k(\xi)=\sqrt{2}\sin((k+1)\pi\xi)$ - see \cite{Brezis}.
\end{ex}

In the following Definition, we introduce finite dimensional subspaces of $H$ and associated orthogonal projections; both are based on the spectral decomposition of $A$.
\begin{defi}\label{def_HPM}
For any $M\in\N$, we define $H_M$ the subspace of $H$ generated by $e_0,\ldots,e_M$,
$$H_M=\text{Span}\left\{e_k;0\leq k\leq M\right\}$$
and $P_M\in\mathcal{L}(H)$ the orthogonal projection onto $H_M$: for any $x=\sum_{k=0}^{+\infty}x_ke_{k}\in H$,
$$P_Mx=\sum_{k=0}^{M}x_ke_{k}.$$
\end{defi}

The domain $D(A)$ of $A$ is equal to $D(A)=\left\{x=\sum_{k=0}^{+\infty}x_ke_{k}\in H, \sum_{k=0}^{+\infty}(\lambda_k)^{2}|x_k|^2<+\infty\right\}$.
More generally, fractional powers of $-A$, are defined for $\alpha\in[0,1]$:
$$(-A)^\alpha x=\sum_{k=0}^{\infty}\lambda_{k}^{\alpha}x_ke_k\in H,$$
with the domains
$$D(-A)^\alpha=\left\{x=\sum_{k=0}^{+\infty}x_ke_{k}\in H,\quad |x|_{\alpha}^{2}:=\sum_{k=0}^{+\infty}(\lambda_k)^{2\alpha}|x_k|^2<+\infty\right\}.$$

\begin{ex}
In the case when $A$ is the Laplace operator with homogeneous Dirichlet boundary conditions on $H=L^2(\mathcal{D})$,
$$D((-A)^{1/2})=H_{0}^{1}(\mathcal{D})\quad D(A)=H_{0}^{1}(\mathcal{D})\cap H^{2}(\mathcal{D}).$$
\end{ex}

When $\alpha\in[0,1]$, it is also possible to define spaces $D(-A)^{-\alpha}$ and operators $(-A)^{-\alpha}$, with norm denoted by $|.|_{-\alpha}$; in particular when $x=\sum_{k=0}^{+\infty}x_ke_{k}\in H$, we have $(-A)^{-\alpha}x=\sum_{k=0}^{+\infty}\lambda_{k}^{-\alpha}x_kf_k$ and $|x|_{-\alpha}^{2}:=\sum_{k=0}^{+\infty}(\lambda_k)^{-2\alpha}|x_k|^2$.

The semi-group $(e^{tA})_{t\geq 0}$ can be defined by the Hille-Yosida Theorem - see \cite{Brezis}. We use the following spectral formula: if $x=\sum_{k=0}^{+\infty}x_ke_k\in H$, then for any $t\geq 0$
$$e^{tA}x=\sum_{k=0}^{+\infty}e^{-\lambda_k t}x_ke_k.$$

For any $t\geq 0$, $e^{tA}$ is a continuous linear operator in $H$, with operator norm $e^{-\lambda_0 t}$. The semi-group $(e^{tA})$ is used to define the solution $Z(t)=e^{tA}z$ of the linear Cauchy problem
$$\frac{dZ(t)}{dt}=AZ(t)\quad \text{with} \quad Z(0)=z.$$

To define solutions of more general PDEs of parabolic type, we use mild formulation, and Duhamel principle.

This semi-group enjoys some smoothing properties that we often use in this work. Here we recall a few important ones, which are easily obtained from the spectral formula given above:
\begin{propo}\label{proporegul}
Under Assumption \ref{hypB}, for any $\sigma\in[0,1]$, there exists $C_\sigma>0$ such that we have:
\begin{enumerate}
\item for any $t>0$ and $x\in H$,
$$|e^{tA}x|_{\sigma}\leq C_\sigma t^{-\sigma}e^{-\frac{\lambda_0}{2}t}|x|_{H}.$$
\item for any $0<s<t$ and $x\in H$,
$$|e^{tA}y-e^{sA}x|_{H}\leq C_\sigma\frac{(t-s)^\sigma}{s^\sigma}e^{-\frac{\lambda_0}{2}s}|x|_H.$$
\item for any $0<s<t$ and $y\in D(-A)^\sigma$,
$$|e^{tA}y-e^{sA}x|_{H}\leq C_\sigma(t-s)^\sigma e^{-\frac{\lambda_0}{2}s}|x|_{\sigma}.$$
\end{enumerate}
\end{propo}

\subsubsection{The nonlinear operator}

The nonlinear operator $F$ is assumed to satisfy some general assumptions, like in \cite{B2}. In Example \ref{exG}, we give the two main kind of operators that can be used in our framework.

{\color{red} Demander juste $F$ de classe $\mathcal{C}^2$ pour simplifier un peu, et ajouter une remarque?}

\begin{hyp}\label{hypG}
The function $F:H\rightarrow H$ is assumed to be bounded and Lipschitz continuous. We denote by $L_F$ the Lipschitz constant of $F$

We also define for each $M\geq 0$ a function $F_M:H_M\rightarrow H_M$, with $F_M(x)=P_MF(x)$ for any $x\in H_M$. We assume that each $F_M$ is twice differentiable, and that we have the following bounds on the derivatives, uniformly with respect to $M$:
\begin{itemize}
\item There exists a constant $C_1$ such that for any $M\geq 0$, $x\in H_M$ and $h\in H_M$
$$|DF_M(x).h|_{H}\leq C_1|h|_{H}.$$
\item There exists $\eta\in[0,1[$ and a constant $C_2$ such that for any $M\geq 0$, $x\in H_M$ and any $h,k\in H_M$ we have
$$|(-A)^{-\eta} D^2F_M(x).(h,k)|\leq C_2|h|_H|k|_H.$$
\item Moreover, there exists a constant $C_3$ such that for any $M\geq 0$, $x\in H_M$ and any $h,k\in H_M$
$$|D^2F_M(x).(h,k)|\leq C_3|h|_{(-B)^\eta}|k|_{H}.$$
\end{itemize}
\end{hyp}

\begin{rem}
Multiplicative noise with appropriate assumptions like in \cite{AndersLars} can be considered; however proofs of the required estimates become much more technical.
\end{rem}

Since $F$ is bounded, the following property is easily satisfied:
\begin{propo}[Dissipativity]\label{propodiss}
There exist $c>0$ and $C>0$ such that for any $x\in D(A)$
\begin{equation}\label{hypdiss}
\langle Ax+F(x),x\rangle \leq -c|x|^2+C.
\end{equation}
\end{propo}

We remark that we have uniform control with respect to the dimension $M$ of the bounds on $F_M:=P_M\circ F$ and on its derivatives, and that \eqref{hypdiss} is also satisfied for $F_M$, with constants $c$ and $C$ independent from $M$.

\begin{ex}\label{exG}
We give some fundamental examples of nonlinearities for which the previous assumptions are satisfied:
\begin{itemize}
\item A function $F:H\rightarrow H$ of class $\mathcal{C}^2$, bounded and with bounded derivatives, 
fits in the framework, with the choice $\eta=0$.
\item The function $F$ can be a \textbf{Nemytskii} operator: let $g:(0,1)\times \R\rightarrow \R$ be a measurable, bounded, function such that for almost every $\xi\in(0,1)$ $g(\xi,.)$ is twice continuously differentiable, with uniformly bounded derivatives.
Then $F(y)$ is defined for every $y\in H=L^2(0,1)$ by
$$F(x)(\xi)=g(\xi,x(\xi)).$$
In general, such functions are not Fr\'echet differentiable, but only G\^ateaux differentiable, with the following expressions:
\begin{gather*}
[DF(x).h](\xi)=\frac{\partial g}{\partial x}(\xi,x(\xi))h(\xi)\\
[D^2F(x).(h,k)](\xi)=\frac{\partial^2 g}{\partial x^2}(\xi,x(\xi))h(\xi)k(\xi).
\end{gather*}
If $h$ and $k$ are only $L^2$ functions, $D^2F(x).(h,k)$ may only be $L^1$; however if $h$ or $k$ is $L^\infty$, it is $L^2$.
The conditions in Assumption \ref{hypG} are then satisfied as soon as there exists $\eta<1$ such that $D(-A)^\eta$ is continuously embedded into $L^\infty(0,1)$ - it is the case for $A$ given in Example \ref{exampleB}, with $\eta>1/4$. Then the finite dimensional spaces $H_N$ are subspaces of $L^\infty$, and differentiability can be shown.
\end{itemize}
\end{ex}

\subsection{The cylindrical Wiener process and stochastic integration in $H$}\label{SectWiener}

In this section, we recall the definition of the cylindrical Wiener process and of stochastic integral on a separable Hilbert space $H$ with norm $|.|_H$. For more details, see \cite{DaP-Z1}.

We first fix a filtered probability space $(\Omega,\mathcal{F},(\mathcal{F}_t)_{t\geq 0},\PP)$. A cylindrical Wiener process on $H$ is defined with two elements:
\begin{itemize}
\item a complete orthonormal system of $H$, denoted by  $(q_i)_{i\in I}$, where $I$ is a subset of $\N$;
\item a family $(\beta_i)_{i\in I}$ of independent real Wiener processes with respect to the filtration $((\mathcal{F}_t)_{t\geq 0})$;
\end{itemize}

then $W$ is defined by
\begin{equation}\label{defWiener}
W(t)=\sum_{i\in I}\beta_i(t)q_i.
\end{equation}

When $I$ is a finite set, we recover the usual definition of Wiener processes in the finite dimensional space $\R^{|I|}$. However the subject here is the study of some Stochastic Partial Differential Equations, so that in the sequel the underlying Hilbert space $H$ is infinite dimensional; for instance when $H=L^{2}(0,1)$, an example of complete orthonormal system is $(q_k)=(\sqrt{2}\sin(k\pi.))_{k\geq 1}$ - see Example \ref{exampleB}.

A fundamental remark is that the series in \eqref{defWiener} does not converge in $H$; but if a linear operator $\Psi:H\rightarrow K$ is Hilbert-Schmidt, then $\Psi W(t)$ converges in $L^2(\Omega,H)$ for any $t\geq 0$.

Moreover, the resulting process does not depend on the choice of the complete orthonormal system $(q_i)_{i\in I}$.

We recall that a bounded linear operator $\Psi:H\rightarrow K$ is said to be Hilbert-Schmidt when
$$|\Psi|_{\mathcal{L}_{2}(H,K)}^{2}:=\sum_{k=0}^{+\infty}|\Psi(q_k)|_{K}^{2}<+\infty,$$
where the definition is independent of the choice of the orthonormal basis $(q_k)$ of $H$.
The space of Hilbert-Schmidt operators from $H$ to $K$ is denoted $\mathcal{L}_{2}(H,K)$; endowed with the norm $|.|_{\mathcal{L}_{2}(H,K)}$ it is an Hilbert space.

The stochastic integral $\int_{0}^{t}\Psi(s)dW(s)$ is defined in $K$ for predictible processes $\Psi$ with values in $\mathcal{L}_2(H,K)$ such that $\int_{0}^{t}|\Psi(s)|_{\mathcal{L}_2(H,K)}^{2}ds<+\infty$ a.s; moreover when $\Psi\in L^2(\Omega\times[0,t];\mathcal{L}_2(H,K))$, the following two properties hold:
\begin{gather*}
\E|\int_{0}^{t}\Psi(s)dW(s)|_{K}^{2}=\E\int_{0}^{t}|\Psi(s)|_{\mathcal{L}_2(H,K)}^{2}ds \text{ (It\^o isometry),}\\
\E\int_{0}^{t}\Psi(s)dW(s)=0.
\end{gather*}
A generalization of It\^o formula also holds - see \cite{DaP-Z1}.

For instance, if $v=\sum_{k\in\N}v_kq_k\in H$, we can define $$\langle W(t),v\rangle =\int_{0}^{t}\langle v,dW(s)\rangle =\sum_{k\in\N}\beta_k(t)v_k;$$
we then have the following space-time white noise property
$$\E\langle W(t),v_1\rangle \langle W(s),v_2\rangle =t\wedge s\langle v_1,v_2\rangle .$$

Therefore to be able to integrate a process with respect to $W$ requires some strong properties on the integrand; in our SPDE setting, the Hilbert-Schmidt properties follow from the assumptions made on the linear coefficients of the equations.

Thanks to Assumption \ref{hypB}, it is easy to show that the following stochastic integral is well-defined in $H$, for any $t\geq 0$:
\begin{equation}\label{stoconvWB}
W^A(t)=\int_{0}^{t}e^{(t-s)A}dW(s).
\end{equation}

It is called a stochastic convolution, and it is the unique mild solution of
$$dZ(t)=AZ(t)dt+dW(t)\quad \text{with} \quad Z(0)=0.$$

Under the second condition of Assumption \ref{hypB}, there exists $\delta>0$ such that for any $t>0$ we have $\int_{0}^{t}\frac{1}{s^\delta}|e^{sA}|_{\mathcal{L}_{2}(H)}^{2}ds<+\infty$; it can then be proved that $W^A$ has continuous trajectories - via the \textit{factorization method}, see \cite{DaP-Z1} - and that for any $1\leq p<+\infty$
\begin{equation}\label{momWBp}
\E\sup_{t\geq 0}|W^A(t)|_{H}^{p}<+\infty.
\end{equation}

We can now define solutions to equation \eqref{eqY}, thanks to the assumptions made on the coefficients: the following result is classical - see \cite{DaP-Z1}:
\begin{propo}
For every $T>0$, $x\in H$, the equation \eqref{eqY} admits a unique mild solution $X\in \text{L}^2(\Omega,\mathcal{C}([0,T],H))$:
\begin{equation}\label{eqmild}
X(t)=e^{tA}x+\int_{0}^{t}e^{(t-s)A}F(X(s))ds+\int_{0}^{t}e^{(t-s)A}dW(s).
\end{equation}
\end{propo}

\section{Definition of the discretization schemes}\label{sectdefnum} 

We will consider approximations in time and in space of the process $X$. In this Section, we introduce the corresponding schemes: a finite element approximation for discretization in space, and a semi-implicit Euler scheme for discretization in time. We also discuss the discretization in space using the spectral decomposition of the operator $A$.

\subsection{Discretization in space : finite element approximation}\label{sectfe}

We use the same framework as in \cite{AndersLars} and \cite{deb-prin}. For precise general references on Finite Element Methods, see for instance \cite{Ciarlet} and \cite{Ern}.

Let $(V_h)_{h\in(0,1)}$ be a family of spaces of continuous piecewise linear functions corresponding to a finite set of nodes (possibly not uniformly distributed) in $\mathcal{D}=(0,1)$ such that $V_h\subset H^1_0(\mathcal{D})=D((-A)^{1/2})$ - in other words, $0$ and $1$ should be included as nodes in the partition of $[0,1]$. The parameter $h$ denotes the mesh size, which is the length of the largest subinterval in the partition.

Let $P_h:H\rightarrow V_h$ denote the orthogonal projection onto the finite dimensional space $P_h$. According to the context, we also consider $P_h$ as a linear operator in $\mathcal{L}(H)$, since $V_h\in H$.

We finally define the approximation of the operator $A$: it is a linear operator $A_h\in\mathcal{L}(V_h)$.
\begin{defi}
The linear operator $A_h:V_h\rightarrow V_h$ is defined such that the following variational equality holds: for any $x_h\in V_h$ and $y_h\in V_h$
\begin{equation*}
\langle A_hx_h,y_h\rangle=\langle Ax_h,y_h\rangle
\end{equation*} 
\end{defi}

We recall a few important properties of the operator $V_h$:
\begin{propo}\label{propA_h}
For any $h\in(0,1)$, $A_h$ is symmetric, such that $-A_h$ is positive definite.

If $N_h$ is the dimension of $V_h$, we denote by $(e_i^h)_{i=0}^{N_h-1}\subset V_h$ an orthonormal eigenbasis corresponding to $-A_h$ with eigenvalues $0<\lambda_0^h\leq \lambda_1^h\leq... \leq\lambda_{N_h-1}^h$.

Then for any $h\in(0,1)$, we have $\lambda_{0}^{h}\geq \lambda_{0}$.
\end{propo}

Indeed, we have
$$ \lambda_0=\inf_{v,u\in H}\langle -Au,v\rangle\leq\inf_{u,v\in V_h} \langle -Au,v\rangle= \inf_{u,v\in V_h} \langle -A_hu,v\rangle=\lambda_0^h.$$

For any $h\in(0,1)$, $A_h$ generates a semi-group on $V_h$, which is denoted $(e^{tA_h})_{t\in\R^+}$. It is also not difficult to define fractional powers $(-A_{h})^{\alpha}$ of $-A_h$, for any $\alpha\in[-1,1]$: for any $x^h=\sum_{i=0}^{N_h-1}x_{i}^{h}e_{i}^{h}\in V_h$, we have
\begin{gather*}
e^{tA_h}x^h=\sum_{i=0}^{N_h-1}e^{-\lambda_{i}^{h}t}x_{i}^{h}e_{i}^{h};
(-A_h)^\alpha x^h=\sum_{i=0}^{N_h-1}(\lambda_{i}^{h})^{\alpha}x_{i}^{h}e_{i}^{h}.
\end{gather*}
The regularization estimates of Proposition \ref{proporegul} are then easily generalized to these semi-groups; moreover bounds are uniform with respect to the mesh size $h\in(0,1)$.

We focus now on the approximations of PDEs - seen as equations in the Hilbert space $H$ -  with equations in finite dimensional spaces $V_h$.



We consider the spatially semi discrete approximation of \eqref{eqY}: $(X_h(t))_{t\in\R^+}$, is a process taking values in $V_h$, such that
\begin{equation}\label{edpsef}
dX^h(t)=A_hX^h(t) dt+F^h(X^h(t))dt+P_hdW(t), \quad X^h(0)=P_hx=P_hX_0,
\end{equation}
where the non-linear coefficient $F^h:V_h\rightarrow V_h$ satisfies $F^h(x)=P_h(F(x))$ for any $x\in V_h$.

We remark that the regularity properties of Assumption \ref{hypG} and the dissipativity inequality \eqref{hypdiss} are satisfied if we replace $A$ (resp. $F$) with $A_h$ (resp. $F^h$).

This equation admits a unique mild solution, such that for any $0\leq t\leq T$
\begin{equation} \label{mildsolef}
X^h(t)=e^{tA_h}P_hx+\int_0^te^{(t-s)A_h}F^h(X^h(s))ds+\int_0^te^{(t-s)A_h}P_hdW(s).
\end{equation}

Notice that the stochastic integral is always well-defined, since for any $h\in(0,1)$ the linear operator $P_h$ has finite rank; on $V_h$, the noise process $P_hW$ is a standard $N_h$-dimensional Wiener process - as is easily seen by expanding $W$ in a complete orthonormal system $(q_{i})_{i\in\N}$ with $q_i=e_{i}^{h}$ for $0\leq i\leq N_h-1$.

To be able to state a convergence result of $X^h$ to $X$, and to give an order of convergence, we now express some important results - see \cite{AndersLars} for more details:
\begin{propo}\label{conv_FEM}
\textbf{(i)} We have an equivalence of norms: there exist two constants $c,C\in(0,+\infty)$, such that for any $h\in(0,1)$, any $\alpha\in[-1/2,1/2]$ and any $x^h\in V_h$,
\begin{equation}\label{estiA_h}
c|(-A_h)^\alpha x^h|\leq |(-A)^\alpha x^h|\leq C|(-A_h)^\alpha x^h|.
\end{equation}
Moreover, we have for any $h\in(0,1)$, $\alpha\in[-1/2,1/2]$, and $x\in H$,
\begin{equation}\label{estiA_hbis}
|(-A_h)^\alpha P_h x|\leq C|(-A)^\alpha x|.
\end{equation}
\textbf{(ii)} Let us denote by $R_h$ the so-called Ritz projector, defined as the orthogonal projection onto $V_h$ in $D((-A)^{1/2})$. We have the identity $R_h=(-A_{h})^{-1}P_h(-A)$ on $D(A)$, and
\begin{equation}\label{estiR_h}
\big|(-A)^{s/2}(I-R_h)(-A)^{-r/2}\big|_{\mathcal{L}(H)}\leq C_{r,s}h^{r-s}\quad \forall 0\leq s\leq 1\leq r\leq 2\\
\end{equation}
\textbf{(iii)} For $P_h$, we have the following error estimate:
\begin{equation}\label{estiP_h}
\big|(-A)^{s/2}(I-P_h)(-A)^{-r/2}\big|_{\mathcal{L}(H)}\leq C_{r,s}h^{r-s}\quad \forall 0\leq s\leq 1 \text{ and } 0\leq s\leq r\leq 2.
\end{equation}
\end{propo}


As a consequence, we get the following important result:
\begin{propo}\label{propoTrace}
For any $\kappa>0$, the linear operator on $H$ $P_h(-A_h)^{-1/2-\kappa}P_h$ is continuous, self-adjoint and semi-definite positive. Moreover,
$$\sup_{0<h<1}\text{Tr}\bigl(P_h(-A_h)^{-1/2-\kappa}P_h\bigr)<+\infty.$$
\end{propo}

The symmetry and the positivity are very important properties for our purpose: indeed, they allow to use inequalities like
$$|\text{Tr}(MN)|\leq |M|_{\mathcal{L}(H)}\text{Tr}(N),$$
for $M,N\in\mathcal{L}(H)$ such that $L$ is symmetric and semi-definite positive.

\underline{Proof}
The operator is well-defined on $H$, and self-adjointness is clear, since $(-A_h)^{-1/2-\kappa}\in\mathcal{L}(V_h)$ is symmetric.

Now from point (i) of Proposition \ref{conv_FEM}, the following linear operators are defined and continuous on $H$: $(-A)^\kappa(-A_h)^{-\kappa}P_h$, and $(-A)^{1/2}P_h(-A_h)^{-1/2}P_h$; their norm is uniformly bounded with respect to $h$.

By duality, the operator $P_h(-A_h)^{-1/2}P_h(-A)^{1/2}$ is well-defined on $H$ - by unique continuous extension from the dense subspace $D((-A)^{1/2})$ - and it has the same norm as $(-A)^{1/2}P_h(-A_h)^{-1/2}P_h$.

Finally, we write that for any $0<h<1$
\begin{align*}
\text{Tr}\bigl(P_h(-A_h)^{-1/2-\kappa}P_h\bigr)&=\text{Tr}\bigl((P_h(-A_h)^{-1/2}P_h(-A)^{1/2})(-A)^{-1/2-\kappa}((-A)^\kappa(-A_h)^{-\kappa}P_h)   \bigr)\\
&\leq \big|P_h(-A_h)^{-1/2}P_h(-A)^{1/2}\big|_{\mathcal{L}(H)}\text{Tr}\bigl((-A)^{-1/2-\kappa}\bigr)\big|(-A)^\kappa(-A_h)^{-\kappa}P_h\big|_{\mathcal{L}(H)}\\
&\leq C\text{Tr}\bigl((-A)^{-1/2-\kappa}\bigr).
\end{align*}
\qed

We now recall a few convergence results, valid on time-intervals of finite length $[0,T]$:
\begin{itemize}
\item in the deterministic case, the order of convergence is $2$;
\item in the stochastic case, the strong order of convergence is $1/2$ - see for instance \cite{Kruse}, \cite{Yan}:

$\forall 0<r<1/2$, $\exists C_{T,r}\in(0,+\infty)$, $\forall h\in(0,1)$, we have $\E|X^h(T)-X(T)|\leq C_{T,r} h^{1/2-r}$;
\item in the stochastic case, the weak order of convergence is $1$ - see \cite{AndersLars}: for any admissible test function $\phi$,

$\forall 0<r<1$, $\exists C_{T,r}\in(0,+\infty)$, $\forall h\in(0,1)$, we have $|\E\phi(X^h(T))-\E\phi(X(T))|\leq C_{T,r} h^{1-r}$.
\end{itemize}

To conclude this part, we introduce a notation which is useful to give some of the results in a compact way.
\begin{defi}\label{defX^0}
For $h=0$, we set $X^{0}=X$, as well as $V_{0}=H$, $A_0=A$, $P_{0}=Id_H$.
\end{defi}

\subsection{Another discretization in space: spectral Galerkin projection}

A tool in our proof will be an additional finite dimensional projection onto the subspaces $H_M$. This approximation allows to justify rigorously the required computations; even if the process $X^h$ takes values in a finite dimensional subspace of $H$, it is easier to prove some estimates with a process taking values in finite dimensional subspaces which are left invariant by the action of $A$ and of the noise. We define here the corresponding approximating processes, and give a few important convergence properties.

Let $M\in\N$. According to Definition \ref{def_HPM}, we can consider an approximate equation in the finite-dimensional subspace $H_M$:
\begin{equation}\label{edpsfini}
dX^{(M)}(t)=AX^{(M)}(t)dt+F_M(X^{(M)}(t))dt+P_MdW(t), \quad X^{(M)}(0)=P_Mx,
\end{equation}
where $F_M=P_M\circ F$ - see also Assumption \ref{hypG}. The process $W^{(M)}:=P_MW$ takes values in $H_M$, it is a standard Wiener process.

For any final time $T\in(0,+\infty)$, it admits a unique mild solution, taking values in $H_M\subset H$ - we recall that $H_M$ is stable by $A$:
$$X^{(M)}(t)=e^{tA}P_Mx+\int_0^te^{(t-s)A}F_M(X^{(M)}(s))ds+\int_0^te^{(t-s)A}P_MdW(s).$$

To study the convergence of $X^{(M)}$ to $X$, the following inequality is useful:
\begin{equation}\label{estiP_M}
|(I-P_M)A^{-r}|_{\mathcal{L}(H)}\leq  C_r\lambda_{M+1}^{-r},\quad 0\leq r\leq 1.
\end{equation}

We then have the following convergence results, for any $T\in(0,+\infty)$, in the stochastic case:
\begin{itemize}
\item the strong order of convergence is $1/4$:

$\forall 0<r<1/2$, $\exists C_{T,r}\in(0,+\infty)$, $\forall M\in\N$, we have $\E|X^{(M)}(T)-X(T)|\leq \frac{C_{T,r}}{\lambda_{M+1}^{1/4-r}}$;
\item the weak order of convergence is $1/2$: for any admissible test function $\phi$,

$\forall 0<r<1$, $\exists C_{T,r}\in(0,+\infty)$, $\forall M\in\N$, we have $|\E\phi(X^{(M)}(T))-\E\phi(X(T))|\leq \frac{C_{T,r}}{\lambda_{M+1}^{1/2-r}}$.
\end{itemize}
Those estimates can be proved with direct computations and the appropriate techniques from \cite{AndersLars} and \cite{deb}. Another possibility is to check that the projectors $P_N$ satisfy the estimates of Proposition \ref{conv_FEM} with $h=\lambda_{N+1}^{-1/2}$, see Example $3.4$ in \cite{Kruse}.
Once again, we define a value for $M=\infty$:
\begin{defi}\label{defXinfty}
For $M=\infty$, we set $X^{(\infty)}=X$, as well as $H_{\infty}=H$ and $P_{\infty}=Id_H$.
\end{defi}


%
%
%

\subsection{Discretization in time}
For each fixed mesh size $h\in(0,1)$, and for $h=0$, we now define a time approximation of the process $X_h$: denoting by $\tau>0$ a time step, we use a semi-implicit Euler scheme to define, for $k\in\mathbb{N}$,
\begin{gather*}
X^h_{k+1}(\tau,x)=X^h_k(\tau,x)+\tau A_hX^h_{k+1}(\tau,x)+\tau P_hF(X^h_{k}(\tau,x))+\sqrt{\tau}P_h\chi_{k+1}\\
X_0^h(\tau,x)=x,
\end{gather*}
where 
$\chi_{k+1}=\frac{1}{\sqrt{\tau}}(W((k+1)\tau)-W(k\tau)).$

To simplify the equations, most of the time we omit the dependence of $X_{k}^{h}$ on the time-step $\tau$ and on the initial condition $x$.

The previous equation can be replaced by
\begin{equation}\label{defYk}
X_{k+1}^{h}=S_{\tau,h} X_{k}^{h}+\tau S_{\tau,h} P_hF(X_{k}^{h})+\sqrt{\tau}S_{\tau,h}P_h\chi_{k+1},
\end{equation}
where $S_{\tau,h}$ is defined by
\begin{equation}\label{defRtau}
S_{\tau,h}=(I-\tau A_h)^{-1}.
\end{equation}
When $h=0$, the process is well-defined in $H$, since it is easily checked that $S_{\tau,0}$ is a Hilbert-Schmidt operator on $H$. When $h>0$, it is well-defined in the finite-dimensional space $V_h$.

%
For the analysis of the convergence of the scheme, we need the following technical estimates on the discrete-time semi-group $(S_{\tau,h}^{j})_{j\in\N}$ for $\tau>0$ and $h\geq 0$:

\begin{lemme}\label{lem6}
For any $0\leq \kappa\leq 1$, $h\in[0,1)$ and $j\geq 1$
$$|(-A_h)^{1-\kappa}S_{\tau,h}^{j}P_h|_{\mathcal{L}(H)}\leq \frac{1}{(j\tau)^{1-\kappa}}\frac{1}{(1+\lambda_0\tau)^{j\kappa}}.$$
Moreover, for any $\beta\geq 1$ and $j\geq \beta$
\begin{equation*}
|(-A_h)^{\beta}S_{\tau,h}^jP_h|_{\mathcal{L}(H)}\leq \frac{\beta^{\beta}}{(j\tau)^{\beta}},
\end{equation*}
and for any $0\leq \beta\leq 1$
\begin{equation*}
|(-A_h)^{-\beta}(S_{\tau,h}-I)P_h|_{\mathcal{L}(H)}\leq 2\tau^\beta.
\end{equation*}
\end{lemme}


\underline{Proof}
Using the notations of Proposition \ref{propA_h}, we have, for any $z\in H$,
\begin{align*}
|(-A_h)^{1-\kappa}S_{\tau,h}^{j}P_hz|_{H}^{2}&=\sum_{i=0}^{N_h}(\lambda_{i}^h)^{2(1-\kappa)}\frac{1}{(1+\lambda_i^h\tau)^{2j}}\langle z,f_i^h\rangle^2\\
&=\frac{1}{(j\tau)^{2(1-\kappa)}}\sum_{i=0}^{N_h}\langle z,f_i^h\rangle^2(\lambda_{i}^h)^{2(1-\kappa)}(j\tau)^{2(1-\kappa)}\frac{1}{(1+\lambda_i^h\tau)^{2j(1-\kappa)}}\frac{1}{(1+\lambda_i^h\tau)^{2j\kappa}}\\
&\leq \frac{1}{(j\tau)^{2(1-\kappa)}}\sum_{i=0}^{N_h}\left(\frac{\lambda_i^h j\tau}{1+\lambda_i^h j\tau}\right)^{2(1-\kappa)}\frac{1}{(1+\lambda_0^h\tau)^{2j\kappa}}\langle z,f_i^h\rangle^2\\
&\leq c|P_hz|_{H}^{2}\frac{1}{(j\tau)^{2(1-\kappa)}}\frac{1}{(1+\lambda_0^h\tau)^{2j\kappa}}\\
&\leq c|z|_{H}^{2}\frac{1}{(j\tau)^{2(1-\kappa)}}\frac{1}{(1+\lambda_0^h\tau)^{2j\kappa}}.
\end{align*}
Above we have used the notation $N_0=+\infty$. To conclude, we use that for any $h\in[0,1)$, $\lambda_0\leq \lambda_0^h$.
The proofs of the two other inequalities are similar - see \cite{B2} for the second one.\qed

\begin{rem}
Later, we often use the following expression for $X_k^{h}$:
\begin{equation}\label{exprYk}
X_{k}^{h}=S_{\tau,h}^{k}P_hx+\tau\sum_{l=0}^{k-1}S_{\tau,h}^{k-l}P_h F(X_{l}^{h})+\sqrt{\tau}\sum_{l=0}^{k-1}S_{\tau,h}^{k-l}P_h\chi_{l+1}.
\end{equation}
The following expression is also useful:
\begin{equation}\label{exprsto}
\sqrt{\tau}\sum_{l=0}^{k-1}S_{\tau,h}^{k-l}P_h\chi_{l+1}=\int_{0}^{t_k}S_{\tau,h}^{k-l_s}P_hdW(s),
\end{equation}
where $l_s=\lfloor \frac{s}{\tau} \rfloor$ - with the notation $\lfloor . \rfloor$ for the integer part.
\end{rem}

For $h\in(0,1)$, we finally introduce the following processes: for $0\leq k\leq m-1$ and $t_{k}\leq t\leq t_{k+1}$
\begin{equation}\label{tild}
\tilde{X}^{h}(t)=X_{k}^{h}+\int_{t_k}^{t}[A_hS_{\tau,h} X_{k}^{h}+S_{\tau,h} P_h F(X_{k}^{h})]ds+\int_{t_k}^{t}S_{\tau,h} P_hdW(s).
\end{equation}

The process $(\tilde{X}^{h}(t))_{t\in \mathbb{R}^+}$ is a natural interpolation in time of the numerical solution $(X_{k}^{h})_{k\in \mathbb{N}}$ defined by \eqref{defYk}: $\tilde{X}^{h}(t_k)=X_{k}^{h}$.


\subsection{A priori bounds on moments}\label{sectusefulestim}
We give a few results on the processes $(X(t))_{t\geq 0}$, $(X^{h}(t))_{t\in \mathbb{R}^+}$ and $(X_{k}^{h})_{k\in \mathbb{N}}$.

All the appearing constants are uniform with respect to $h\in(0,1)$.

\begin{lemme}\label{lem1}
For any $p\geq 1$, there exists a constant $C_p>0$ such that for every $h\in(0,1)$, $t\geq 0$ and $x\in H$
$$\E|X^{h}(t,x)|^p\leq C_p(1+|x|^p).$$
\end{lemme}

\begin{lemme}\label{lem2}
For any $p\geq 1$, $\tau_0>0$, there exists a constant $C>0$ such that for every $h\in(0,1)$, $0<\tau\leq \tau_0$, $k\in \N$, $t\geq 0$ and $x\in H$
$$\E|X_{k}^{h}|^{p}\leq C(1+|x|^p) \quad \text{and for } h\in(0,1) \text{ we have }\quad \E|\tilde{X}^{h}(t)|^p\leq C(1+|x|^p).$$
\end{lemme}

\underline{Proof of Lemmas \ref{lem1} and \ref{lem2}}
The case $h=0$ is treated in \cite{B2}. We thus only treat the case $h\in(0,1)$, with similar methods.

Using the mild formulation, in the continuous-time situation we need to control three terms:
$$|e^{tA_h}P_hx|\leq e^{-\lambda_0^h t}|P_hx|\leq e^{-\lambda_0 t}|x|,$$
since $\lambda_0^h\geq \lambda_0$ thanks to Proposition \ref{propA_h};
$$|\int_{0}^{t}e^{(t-s)A_h}P_hF(X_h(s))ds|\leq C\int_{0}^{t}e^{-\lambda_0^h(t-s)}ds\leq \frac{C}{\lambda_0^h}\leq \frac{C}{\lambda_0},$$
thanks to the boundedness of $F$;

\begin{align*}
\E[|\int_{0}^{t}e^{(t-s)A_h}P_hdW(s)|^p]&\leq C_p\E[|\int_{0}^{t}e^{(t-s)A_h}P_hdW(s)|^2]^{p/2}\\
&\leq C_p\bigl(\int_{0}^{t}\text{Tr}(e^{(t-s)A_h}P_he^{(t-s)A_h})ds\bigr)^{p/2}\\
&\leq C_p\bigl(\int_{0}^{t}\text{Tr}(P_hA_{h}^{-1/2-\kappa}P_h)|A_{h}^{1/2+\kappa}e^{2(t-s)A_h}P_h|_{\mathcal{L}(H)}ds\bigr)^{p/2}\\
&\leq C_p\bigl(\int_{0}^{t}\frac{1}{(t-s)^{1/2+\kappa}}e^{-\lambda_{0}^{h}(t-s)}ds\bigr)^{p/2}\leq C_p,
\end{align*}
where we have reduced the estimate to the case $p=2$ since the stochastic integral is a Gaussian random variable, and then we have used It\"o's isometry formula, Proposition \ref{proporegul} and Proposition \ref{propoTrace}, with any $\kappa>0$.

The proof of the first estimate of Lemma \ref{lem2} goes along the same way, using the discrete-time mild formulation \eqref{exprYk}, with also three quantities to control:
$$|S_{\tau,h}^{k}P_hx|\leq \frac{1}{(1+\lambda_{0}\tau)^k}|x|,$$
thanks to the first estimate of Lemma \ref{lem6} and the boundedness of $F$;
\begin{align*}
|\tau\sum_{l=0}^{k-1}S_{\tau,h}^{k-l}P_h F(X_{l}^{h})|&\leq C\tau\sum_{l=0}^{k-1}\frac{1}{(1+\lambda_0\tau)^{k-l}}\\
&\leq C\frac{\tau}{(1/(1+\lambda_0\tau))-1}\leq C;
\end{align*}
finally as before we only need to study the case $p=2$, and we use \eqref{exprsto} to get
\begin{align*}
\E|\sqrt{\tau}\sum_{l=0}^{k-1}S_{\tau,h}^{k-l}P_h\chi_{l+1}|^2&=\sum_{l=0}^{k-1}\tau\text{Tr}(S_{\tau,h}^{k-l}P_hS_{\tau,h}^{k-l})\\
&\leq \text{Tr}(P_hA_{h}^{-1/2-\kappa}P_h)\tau\sum_{l=0}^{k-1}|S_{\tau,h}^{2(k-l)}A_{h}^{1/2+\kappa}|_{\mathcal{L}(H)}\\
&\leq C\tau\sum_{l=0}^{k-1}\frac{1}{((k-l)\tau)^{1/2+\kappa}}\exp\big(-(1/2-\kappa)\frac{\log(1+\lambda_0\tau)}{\tau}(k-l)\tau\big)\\
&\leq C\int_{0}^{+\infty}\frac{1}{t^{1/2+\kappa}}\exp(-ct)dt,
\end{align*}
for some $c>0$.

The proof of the second estimate of Lemma \ref{lem2} using \eqref{tild} is then straightforward.
\qed

\section{Asymptotic behavior of the processes and invariant laws}\label{sectasymptu}

First, we focus on the existence of invariant measures for the continuous and discrete time processes. We use the well-known Krylov-Bogoliubov criterion - see \cite{DaP-Z2}. Tightness comes from two facts: $D(-A)^\gamma$ is compactly embedded in $H$ when $\gamma>0$, and when $\gamma<1/4$ we can control moments with the same techniques as for proving the Lemmas \ref{lem1} and \ref{lem2}:
\begin{lemme}\label{estimtight}
For any $0<\gamma<1/4$, $\tau>0$ and any $x\in H$, there exists $C(\gamma,\tau,x),C(\gamma,x)>0$ such that for every $h\in(0,1)$, $m\geq 1$ and $t\geq 1$
$$\E|X_{m}^{h}(\tau,x)|_{\gamma}^{2}\leq C(\gamma,\tau,x) \quad \text{ and }\quad \E|X^{h}(t,x)|_{\gamma}^{2}\leq C(\gamma,x).$$
\end{lemme}

For $h\in[0,1)$, uniqueness of the invariant probability measure for the continuous time process $(X^h(t))_{t\in\R^+}$ can be deduced from the well-known Doob Theorem - see \cite{DaP-Z2}. Indeed, since in equation \eqref{eqY} noise is additive and non-degenerate, the Strong Feller property and irreducibility can be easily proved. In the proof of the main Theorem \ref{th}, we also need speed of convergence, and thanks to a coupling argument we get the following exponential convergence result :
\begin{propo}\label{propoexpy1y2}
There exist $c>0$, $C>0$ such that for any bounded test function $\phi: H\rightarrow \R$, any $t\geq 0$, any $h\in[0,1)$ and any $x_1,x_2\in V_h$
\begin{equation}\label{cvexpy1y2}
|\E\phi(X^{h}(t,x_1))-\E\phi(X^{h}(t,x_2))|\leq C\|\phi\|_{\infty}(1+|x_1|^2+|x_2|^2)e^{-ct}.
\end{equation}
\end{propo}

\begin{rem}
A proof of this result can be found in Section $6.1$ in \cite{deb-hu-tess}. In this proof it is obvious that $c$ and $C$ are independent of $h$.
\end{rem}

The idea of coupling relies on the following formula: if $\nu_1$ and $\nu_2$ are two probability measures on a state space $S$, their total variation distance satisfies
$$d_{TV}(\nu_1,\nu_2)=\inf\left\{\PP(X_1\neq X_2)\right\},$$
which is an infimum over random variables $(X_1,X_2)$ defined on a same probability space, and such that $X_1\sim\nu_1$ and $X_2\sim \nu_2$.

Roughly speaking, the principle is to define a coupling $(Z_1(t,x_1,x_2),Z_2(t,x_1,x_2))_{t\geq 0}$ for the processes $(X(t,x_1))_{t\geq 0}$ and $X((t,x_2))_{t\geq 0}$ such that the coupling time $\mathcal{T}$ of $Z_1$ and $Z_2$ - i.e. the first time the processes are equal - has an exponentially decreasing tail.

This technique was first used in the study of the asymptotic behavior of Markov chains - see \cite{Brem}, \cite{Doe}, \cite{Lin}, \cite{Me-Twee} - and was later adapted for SDEs and more recently for SPDEs - see for instance \cite{ku-shi}, \cite{matt}, \cite{Mue}.

\begin{cor}\label{CorErgo}
For any $h\in[0,1)$, the process $X^{h}$ admits a unique invariant probability measure $\mu^{h}$, such that for any bounded test function $\phi: H\rightarrow \R$, $t\geq 0$ and $x\in V_h$ we have
\begin{equation}\label{cvexpy_int}
|\E\phi(X^{h}(t,x))-\int_{V_h}\phi d\mu^{h}|\leq C\|\phi\|_{\infty}(1+|x|^2)e^{-ct}.
\end{equation}
We use the notation $\mu=\mu^{0}$, when $h=0$: it is the invariant law of the non-discretized process, for which we show an approximation result.
\end{cor}

However, the situation is more complex for the discrete time approximations: given a time-step $\tau>0$, we do not know whether uniqueness also holds for the numerical approximation $(X_k^{h})_{k\in \N}$. In the following Remark \ref{remstrictergo} below, we describe a strict dissipativity assumption on the non-linear coefficient which ensures ergodicity by a straightforward argument. Without this assumption, it is not clear whether ergodicity holds for small time-steps $0<\tau\leq \tau_{\text{ergo}}$, where $\tau_{\text{ergo}}$ can be chosen indepently of $h\in[0,1)$; the answer to this question will be the subject of future works.



\begin{rem}\label{remstrictergo}
Let $h\in[0,1)$ be fixed. A sufficient condition for the uniqueness of the invariant probability measure of the discrete time process $(X_k^{h})_{k\in \N}$ is the strict dissipativity assumption
$$L_F<\lambda_0,$$
where we recall that $L_F$ denotes the Lipschitz constant of $F$.

Then trajectories of the processes $(X_t^{h})_{t\in\R^+}$ and $(X_k^{h})_{k\in\N}$ issued from different initial conditions $x_1$ and $x_2$ and driven by the same noise process are exponentially close when time increases: for any $\tau_0>0$, there exists $c>0$ such that for any $0<\tau\leq \tau_0$, $h\in[0,1)$, $k\geq 0$ and $t\geq 0$ we have almost surely
\begin{gather*}
|X^{h}(t,x_1)-X^{h}(t,x_2)|\leq e^{-(\lambda_0-L_F)t}|x_1-x_2|\\
|X_k^{h}(\tau,x_1)-X_k^{h}(\tau,x_2)|\leq e^{-ck\tau}|x_1-x_2|.
\end{gather*}
The proof of the uniqueness of the invariant law is now easy, and in particular we do not use Proposition \ref{propoexpy1y2}.
\end{rem}



Results are the same when we consider the spectral Galerkin discretization:
\begin{propo}\label{propoErgoM}
There exist $c>0$, $C>0$ such that for any bounded test function $\phi: H\rightarrow \R$, any $t\geq 0$, any $M\in\N\cup\left\{\infty\right\}$ and any $x_1,x_2\in H_M$
\begin{equation}\label{cvexpy1y2M}
|\E\phi(X^{(M)}(t,x_1))-\E\phi(X^{(M)}(t,x_2))|\leq C\|\phi\|_{\infty}(1+|x_1|^2+|x_2|^2)e^{-ct}.
\end{equation}
Moreover, for any $M\in\N\cup\left\{\infty\right\}$, the process $X^{(M)}$ admits a unique invariant probability measure $\mu^{(M)}$, such that for any bounded test function $\phi: H\rightarrow \R$, $t\geq 0$ and $x\in H_M$ we have
\begin{equation}\label{cvexpy_intM}
|\E\phi(X^{(M)}(t,x))-\int_{H_M}\phi d\mu^{(M)}|\leq C\|\phi\|_{\infty}(1+|x|^2)e^{-ct}.
\end{equation}
\end{propo}

With the notations of Definitions \ref{defX^0} and \ref{defXinfty}, we have $\mu^{(\infty)}=\mu=\mu^0$.

As consequence of Proposition \ref{propoErgoM}, we obtain the following Lemma:
\begin{lemme}\label{LemmeMinfini}
For any bounded test function $\phi\in\mathcal{C}(H)$, we have
\begin{equation*}
\overline{\phi}_M:=\int_{H_M}\phi(z)d\mu^{(M)}(z)\underset{M\to \infty}{\longrightarrow} \int_{H}\phi d\mu=:\overline{\phi}.
\end{equation*}
\end{lemme}

\underline{Proof of Lemma \ref{LemmeMinfini}:}
For any $t\geq 0$ and any fixed initial condition $x\in H$, we have
\begin{align*}
\int_{H_{M}}\phi(z)d\mu^{(M)}(z)-\int_{H}\phi(z)d\mu(z)&=\int_{H_{M}}\phi(z)d\mu^{(M)}(z)-\E\phi(X^{(M)}(t)))\\
&+\E\phi(X^{(M)}(t))-\E\phi(X(t))\\
&+\E\phi(X(t))-\int_{H}\phi(z)d\mu(z).
\end{align*}
We thus get that for any $t>0$
$$\limsup_{M\rightarrow +\infty}|\int_{H_{M}}\phi(z)d\mu^{(M)}(z)-\int_{H}\phi(z)d\mu(z)|\leq C\exp(-ct),$$
and it remains to take $t\rightarrow +\infty$. Notice that the constant $c$ does not depend on dimension $M$.
\qed

The speed of convergence in Lemma \ref{LemmeMinfini} will be given below in Remark \ref{RemSpeedLemmeMinfini}.

\section{The convergence results}\label{sectres}
We now state our main result, as well as a few important consequences.

For an admissible test function $\phi$ we define $\parallel \phi \parallel_{2,\infty}=\sup_{0\leq j\leq 2}(\parallel D^j\phi\parallel_{\infty})$.

\begin{theo}\label{th}
For any $0<\kappa<1/2$, $\tau_0$, there exists a constant $C>0$ such that for any $\mathcal{C}^2_b(H)$ function $\phi$,  $h\in(0,1)$, $N\geq 1$, $x\in H$ and $0<\tau\leq\tau_0$
\begin{equation*}
\Big|\frac{1}{N}\sum_{m=0}^{N-1}\Bigl(\E\phi(X_{m}^{h})-\overline{\phi}\Bigr)\Big|\leq C\parallel \phi \parallel_{2,\infty}(1+|x|^3)\Big(1+(N\tau)^{-1+\kappa}+(N\tau)^{-1}\Big)\Big(\tau^{1/2-\kappa}+h^{1-\kappa}+\frac{1}{N\tau}\Big),
\end{equation*}
where $\overline{\phi}=\int_{H}\phi(z)\mu(dz)$.
\end{theo}

This result can be interpreted with a \textit{statistical} point of view: $\frac{1}{N}\sum_{m=0}^{N-1}\E\phi(X_{m}^{h})$ is an \textit{estimator} of the average $\overline{\phi}=\int_{H}\phi(z)\mu(dz)$ of the admissible test function $\phi$ with respect to the invariant law $\mu$ of the SPDE. The Theorem \ref{th} gives an error bound on its \textit{bias}.

Of the two factors in parenthesis in the Theorem, only the second one is important - the presence of the first one is for technical estimates which degenerate at time $0$ whereas we are interested at the asymptotic behavior of the quantity. The main observation is that the orders of convergence with respect to $\tau$ and $h$ are given by the corresponding weak orders $1/2$ and $1$ in the approximation of $X(T)$ for a fixed value of the final time $T<+\infty$ - given in \cite{deb} and \cite{AndersLars}. The aim of this paper is to show how the corresponding error bounds are preserved asymptotically - under appropriate conditions.



An interesting supplementary result would concern the study of the \textit{statistical error}. In \cite{MattStuTre}, two more error bounds are proved: first in the mean-square sense, and then in an almost sure statement - thanks to an argument of Borel-Cantelli type. We have not been able to treat our problem in a similar way. We claim that it is for the following reason. The right order of convergence with respect to $\tau$ in Theorem \ref{th} is obtained thanks to an appropriate integration by parts formula - as explained in the Introduction; the study of the mean-square error - now in a stronger sense - implies that the use of such a technique seems impossible. To generalize the results of \cite{MattStuTre} in the infinite dimensional setting, new arguments should be found.


The additional term $\frac{1}{N\tau}$ corresponds to the bias introduced between the average in time and its limit when time increases.

From Theorem \ref{th}, it is easy to obtain approximation results when only one type of discretization is applied.

Results on distance between invariant laws of the various processes can now be given. As explained in Section \ref{sectasymptu}, without the time-discretization ergodicity holds for the spatially discretized, time-continuous process $X^h$ for any $h\in(0,1)$, while as soon as discretization in time is applied it is not clear whether it holds for small enough time-steps, uniformly with respect to $h$.

However, as a consequence of Theorem \ref{th} we obtain error bounds controlling the distance between the average of admissible test functions with respect to the possibly non unique ergodic invariant laws of the discretized process and the invariant law of the SPDE.

\begin{propo}\label{PrInv}
For any $0<\kappa<1/2$, $\tau_0>0$, there exists a constant $C>0$ such that the following holds:

for any $0<\tau<\tau_0$ and $h\in(0,1)$, assume that $\mu^{\tau,h}$ is an ergodic invariant law of $(Y_k^{h})_{k\in\N}$; then for any admissible test function $\phi$, we have
\begin{equation*}
|\int_{H}\phi(z)d\mu(z)-\int_{V_h}\phi(z)d\mu^{\tau,h}(z)|\leq C\parallel \phi \parallel_{2,\infty}\Big(\tau^{1/2-\kappa}+h^{1-\kappa}\Big).
\end{equation*}
\end{propo}
The proof of this result is easy - we let go $N$ to $\infty$ in the estimate of Theorem \ref{th}, and the convergence of the time-average for $\mu^{\tau,h}$ a.e. initial condition, see also \cite{B2}.

The above result also holds for invariant laws having a finite third order moment.

\begin{rem}\label{RemSpeedLemmeMinfini}
It is also possible to derive the error estimates as in Theorem \ref{th} and Proposition \ref{PrInv} when discretization in space is done with the spectral approximation in dimension $M$ instead of using a finite element method with mesh-size $h$. For instance, we precise the speed of convergence in  : for any $0<\kappa<1/2$ and $\phi\in\mathcal{C}^2_b(H)$, there exists a constant $C_{\kappa}$ such that
\begin{equation*}
|\int_{H}\phi d\mu-\int_{H_M}\phi_Md\mu^{(M)}|\leq C_{\kappa}\frac{1}{\lambda_{M+1}^{1/2-\kappa}}.
\end{equation*}
The strategy - developed in Section \ref{sectdescproof} - remains the same, where we use an auxiliary dimension variable, say $L$, going to $\infty$ for the ambient space where the Poisson equation is used, while $M$ is fixed.
\end{rem}

\section{Description of the proof}\label{sectdescproof}


We fix the time step $\tau$, as well as $N\in \N$; we then introduce the notation $T=N\tau$. We also define, for $k\in\mathbb{N}$, $t_k=k\tau$. $\kappa>0$ is a parameter, which is be supposed to be small enough. We also control $\tau$: for some $\tau_0>0$, $\tau\leq \tau_0$.

\subsection{Strategy}

The three key ingredients to prove Theorem \ref{th} are the use of an additional finite dimensional projection onto the subspaces $H_M$,  the use of the solution of the Poisson equation as in \cite{MattStuTre} (see Sub-section \ref{subsPoi}) and an integration by parts formula issued from Malliavin calculus as in \cite{B2,deb} (see Sub-Section \ref{sectmallia}).

We will use the Poisson equation in finite dimension, then we will use the following decomposition:
\begin{align*}
\frac{1}{N}\sum_{m=0}^{N-1}\E\phi(X_m^{h})-\overline{\phi}=&
\frac{1}{N}\sum_{m=0}^{N-1}\E\phi(P_MX_m^{h})-\overline{\phi}_M\\
&+\overline{\phi}_M-\overline{\phi}+\frac{1}{N}\sum_{m=0}^{N-1}\Big(\E\phi(X^h_m)-\E\phi(P_MX_m^{h})\Big),
\end{align*}
where we recall that $P_M$ is the orthogonal projection of $H$ into $H_M$  and $\overline{\phi}_M$ is define at the end of Section \ref{sectasymptu}.

It is obvious the the last term converge to $0$ when $M\rightarrow +\infty$. The convergence to $0$  of the second term is not difficult, it has been proved in Section \ref{sectasymptu}.
The proof of the estimate of the first term is very technical, so for pedagogy, in Sub-Section \ref{subsstrat} we introduce the decomposition of the error and identify the three terms which we control later in Sub-Section \ref{ssectesti0}, Sub-Section \ref{ssectesti2} and Sub-Section \ref{ssectesti1}.

\subsection{Some results on the Poisson equation in finite dimension}\label{subsPoi}
Let $M\in\left\{1,2,\ldots\right\}$. Let $\phi\in\mathcal{C}^2_b(H)$. We define $\Psi^{(M)}:H_M\rightarrow\R$ by  the unique solution of the Poisson equation
\begin{equation}\label{eqPoiss}
\mathcal{L}^{(M)}\Psi^{(M)}=\phi P_M-\overline{\phi}_M\text{ and } \int_{H_M}\Psi^{(M)} d\mu^{(M)}=0,
\end{equation} 
where $\mathcal{L}^{(M)}$ is the infinitesimal generator of the SPDE  \eqref{edpsfini} defined for functions of class $\mathcal{C}^2$ $\psi:H\rightarrow\R$ and for any $x\in H$ by
\begin{equation*}
\mathcal{L}^{(M)}\psi(x)=\langle AP_Mx+P_MF(x),D\psi(x)\rangle +\frac{1}{2}\text{Tr}(P_MD^2\psi(x)).
\end{equation*}

In the following, we will need to control the first and the second derivatives of $\Psi$.
The Proposition below is the essential result that we need. It is the same kind of estimation used in \cite{B2,deb} to obtain weak order of convergence $1/2$.

\begin{propo}\label{Prop220}
Let $M\in\left\{1,2,\ldots\right\}$. Let $\phi\in\mathcal{C}^2_b(H)$. The function $\Psi^{(M)}$ defined for any $x\in H_M$ by
\begin{equation*}
\Psi^{(M)}(x)=\int_0^{+\infty}\E\Big(\phi(X^{(M)}(t,x))-\overline{\phi}_{M}\Big)dt
\end{equation*}
is of class $\mathcal{C}^2$ and the unique solution of \eqref{eqPoiss}. Moreover, we have the estimates below: for $0\leq\beta,\gamma<1/2$ and $x\in H_M$
\begin{equation*}\label{eq221}
|\Psi^{(M)}(x)|\leq C(1+|x|^2)\parallel \phi\parallel_{\infty},
\end{equation*}
\begin{equation}\label{eq222}
|D\Psi^{(M)}(x)|_{\beta}\leq C_{\beta}(1+|x|^2)\parallel \phi\parallel_{1,\infty}
\end{equation}
and
\begin{equation}\label{eq223}
|(-A)^{\beta}D^2\Psi^{(M)}(x)(-A)^{\gamma}|_{\mathcal{L}(H_M)}\leq C_{\beta,\gamma}(1+|x|^2)\parallel\phi\parallel_{2,\infty},
\end{equation}
where $\parallel \phi \parallel_{i,\infty}=\sup_{0\leq j\leq i}(\parallel D^j\phi\parallel_{\infty})$.
\end{propo}

\begin{rem}\label{rem220}
In fact, the result on $D\Psi$ is also true for $\beta<1$ and the result on $D^2\Psi$ is also true for $\beta<1$, $\gamma<1$ and $\beta+\gamma<1$. Moreover, all the constants are uniform with respect to $M\in\left\{1,2,\ldots\right\}$.
\end{rem}
A proof of this result can be found in the Appendix.

\subsection{A Malliavin integration by parts formula}\label{sectmallia}
Let $h\in(0,1)$ be fixed - the case $h=0$ is not required in the calculations.

As explained in the Introduction, one of the key tools to obtain the right weak order is a transformation of some spatially irregular terms involving the stochastic integral with respect to the cylindrical Wiener process, into more suitable, deterministic ones, thanks to an integration by parts formula, issued from Malliavin calculus - see \cite{Nua}, \cite{Sanz}.

The notations here are the same as in \cite{deb}, where the following useful integration by parts formula is given - see Lemma $2.1$ therein:
\begin{lemme}\label{lemintbyparts}
For any $F\in \mathbb{D}^{1,2}(V_h)$, $u\in \mathcal{C}_{b}^{2}(V_h)$ and $\Psi\in L^2(\Omega\times[0,T],\mathcal{L}_2(V_h))$ an adapted process,
\begin{equation}\label{formulaintbyparts}
\E[Du(F).\int_{0}^{T}\Psi(s)dW^{(M)}(s)]=\E[\int_{0}^{T}\text{Tr}(\Psi(s)^*D^2u(F)\DM_sF)ds],
\end{equation}
where $\DM_sF:\ell\in H\mapsto \DM_{s}^{\ell}F\in V_h$ stands for the Malliavin derivative of $F$, and $\mathbb{D}^{1,2}(V_h)$ is the set of $H$-valued random variables $F=\sum_{i\in\N,i\leq N_h} F_if_i$, with $F_i\in \mathbb{D}^{1,2}$ the domain of the Malliavin derivative for $\R$-valued random variables for any $i$.
\end{lemme}

\begin{rem}
This Lemma remains valid if $u$ is not assumed to be bounded but only $u\in\mathcal{C}^2(V_h)$ provided the expectations and the integral above are well defined. This is easily seen by approximation of $u$ by bounded functions.
\end{rem}

Some care must be taken when controlling the Malliavin derivative of $\tilde{X}^{h}$: in general it is not possible to obtain uniform estimates with respect to time - unless for instance a strict dissipativity condition is satisfied, as in Remark \ref{remstrictergo}.


In the proof of technical estimates below, we circumvent this problem by using these derivatives only at times $t_k=k\tau$ and $s$ such that $t_{k-l_s}\leq 1$. The Lemma \ref{lem5} gives though the most general estimate.

\begin{lemme}\label{lem5}
For any $0\leq \beta<1$ and $\tau_0>0$, there exists a constant $C>0$ such that for every $h\in(0,1)$, $k\geq 1$, $0<\tau\leq \tau_0$ and $s\in \left[0,t_k\right]$
$$|(-A_h)^\beta\DM_{s}^{.}X_{k}^{h}|_{\mathcal{L}(V_h)}\leq C(1+L_F\tau)^{k-l_s}(1+\frac{1}{(1+\lambda_0\tau)^{(1-\beta)(k-l_s)}t_{k-l_s}^{\beta}}).$$
Moreover, if $t_k\leq t< t_{k+1}$, we have
$$|(-A_h)^\beta\DM_{s}^{\ell}\tilde{X}^{h}(t)|_{\mathcal{L}(V_h)}\leq C|(-A_h)^\beta\DM_{s}^{\ell}X_{k}^{h}|_{\mathcal{L}(V_h)}.$$
\end{lemme}

We want to emphasize that the constant in Lemma \ref{lem5} is uniform with respect to $h\in(0,1)$.

\underline{Proof}
According to the definition of $\DM_{s}F$ as a linear operator in $V_h$, we need to control $|(-A_h)^{\beta}\DM_{s}^{\ell}\tilde{X}^{h}(t)|$ and $|(-A_h)^{\beta}\DM_{s}^{\ell}X_{k}^{h}|$, uniformly with respect to $x\in V_h$ with $|\ell|\leq 1$.

The second inequality is a consequence of the following equality for $s\leq t_k\leq t<t_{k+1}$, thanks to \eqref{tild}:
$$\DM_{s}^{\ell}\tilde{X}^{h}(t)=
\DM_{s}^{\ell}X_{k}^{h}+(t-t_m)(\tau A_h S_{\tau,h} \DM_{s}^{\ell}X_{k}^{h}+R_\tau P_h F(X_{k}^{h})\DM_{s}^{\ell}X_{k}^{h}),$$
and the conclusion follows since $$\sup_{h\in(0,1)}|\tau A_h S_{\tau,h}|_{\mathcal{L}(V_h)}<+\infty.$$

Now we prove the first estimate, and we fix $h\in(0,1)$. For any $k\geq 1$, $\ell\in V_h$ and $s\in\left[0,t_k\right]$, using the chain rule for Malliavin calculus and expressions \eqref{exprYk} and \eqref{exprsto}, we have
$$\DM_{s}^{\ell}X_{k}^{h}=S_{\tau,h}^{k-l_s}\ell+\tau\sum_{i=l_s+1}^{k-1}S_{\tau,h}^{k-i}D(P_hF)(X_{i}^{h}).\DM_{s}^{\ell}X_{i}^{h}.$$
We recall that $l_s$ denotes the integer part of $\frac{s}{\tau}$, so that when $i\leq l_s$ we have $\DM_{s}^{\ell}X_{i}^{h}=0$.

As a consequence, the discrete Gronwall Lemma ensures that for $k\geq l_s+1$
$$|\DM_{s}^{\ell}X_{k}^{h}|\leq (1+L_F\tau)^{k-l_s}|\ell|.$$
Now using Lemma \ref{lem6}, we have
$$|(-A_h)^\beta \DM_{s}^{\ell}X_{k}^{h}|\leq \frac{1}{(1+\lambda_0\tau)^{(1-\beta)(k-l_s)}t_{k-l_s}^{\beta}}|\ell|+L_F\tau\sum_{i=l_s+1}^{k-1}\frac{(1+L_F\tau)^{l-l_s}}{(1+\lambda_0\tau)^{(1-\beta)(k-i)}t_{k-i}^{\beta}}|\ell|.$$
To conclude, we see that when $0<\tau\leq \tau_0$
$$\tau\sum_{i=l_s+1}^{k-1}\frac{1}{(1+\lambda_0\tau)^{(1-\beta)(k-i)}t_{k-i}^{\beta}}\leq C\int_{0}^{+\infty}t^{-\beta}\frac{1}{(1+\lambda_0\tau)^{(1-\beta)\frac{t}{\tau}}}dt\leq C<+\infty.$$
\qed

\subsection{Strategy for the estimate of $\lim_{M\rightarrow \infty}\frac{1}{N}\sum_{m=0}^{N-1}\E\phi(P_MX_m^{h})-\overline{\phi}_M$}\label{subsstrat}
 Let $M\in\left\{1,2,\ldots\right\}$, we define the function $\PsiM$  for $x\in H$ by
$$\tilde{\Psi}^{(M)}(x)=\Psi^{(M)}(P_Mx),$$ 
where $\Psi^{(M)}$ is the solution of the Poisson equation \eqref{eqPoiss}. Using the identifications introduced in Remark \ref{remident} we have for any $x\in H$
\begin{gather*}
D\tilde{\Psi}^{(M)}(x)=P_MD\Psi^{(M)}(P_Mx),\\
D^2\tilde{\Psi}^{(M)}(x)=P_MD^2\Psi^{(M)}(P_Mx)P_M,
\end{gather*}
then it is easy to show that $\PsiM$ and $\Psi^{(M)}$ verify the same estimates (see Proposition \ref{Prop220}).\\

The main idea is to give an expansion of $\E\PsiM(X_{m+1}^h)-\E\PsiM(X_m^h)$ to sum it for $m=1,...,N-1$ and divide by $N\tau$. The function $\PsiM$ is defined  for $x\in H$ by
\begin{equation*}
\tilde{\Psi}^{(M)}(x)=\Psi^{(M)}(P_Mx),
\end{equation*}
where $\Psi^{(M)}$ is the solution of the Poisson equation \eqref{eqPoiss}. It is easy to show that $\PsiM$ and $\Psi^{(M)}$ verify the same estimates (see Proposition \ref{Prop220}).

We need the continuous time interpolation of the numerical process $\tilde{X}^{h}$ defined by \eqref{tild}. For all $m\in\mathbb{N}$, we can associate to $\tilde{X}^{h}$ on $[t_m,t_{m+1}]$ the generator $\mathcal{L}^{\tau,m,h}$ defined for $x\in V_h$ and $\phi\in\mathcal{L}(H)$ by
\begin{equation}
\mathcal{L}^{\tau,m,h}\phi(x)=\langle S_{\tau,h}A_hX_m^{h}+S_{\tau,h}P_hF(X_m^{h}),D\phi(x)>+\frac{1}{2}\text{Tr}(S_{\tau,h}S_{\tau,h}^*P_hD^2\phi(x)).
\end{equation}

Thanks to the It\^o formula and Proposition \ref{Prop220}, we have for any integer $m$
$$
\E\PsiM(X_{m+1}^{h})-\E\PsiM(X_{m}^{h})=\int_{t_m}^{t_{m+1}}\E\mathcal{L}^{\tau,m,h}\PsiM(\tilde{X}^{h}(s))ds.
$$
We also need the markov generator $\mathcal{L}^h$ of the finite element solution $X^h$ to decompose this term. The generator $\mathcal{L}^h$ is given  for $x\in V_h$ by
\begin{equation*}
\mathcal{L}^h\phi(x)=\langle A_hx+P_hF(x),D_x\phi(x)\rangle+\frac{1}{2}Tr(P_hD^2_{xx}\phi(x)).
\end{equation*} 
We have the following decomposition:
\begin{align*}
\E\PsiM(X_{m+1}^{h})-\E\PsiM(X_{m}^{h})=&\int_{t_m}^{t_{m+1}}\E\Bigl(\mathcal{L}^{\tau,m,h}-\mathcal{L}^{h}\Bigr)\PsiM(\tilde{X}^{h}(s))ds\\
&+\int_{t_m}^{t_{m+1}}\E\Bigl(\mathcal{L}^{h}-\mathcal{L}^{(M)}\Bigr)\PsiM(\tilde{X}^{h}(s))ds\\
&+\int_{t_m}^{t_{m+1}}\E\mathcal{L}^{(M)}\PsiM(\tilde{X}^{h}(s))ds.
\end{align*}

Using the following equality for $x\in H$
\begin{equation*}
\mathcal{L}^{(M)}\PsiM(x)=\mathcal{L}^{(M)}\Psi^{(M)}(P_Mx)+\langle P_MF(x)-P_MF(P_Mx),D\Psi^{(M)}(P_Mx)\rangle
\end{equation*}
and the definition of $\Psi^{(M)}$, we get
\begin{align*}
\E\PsiM(X_{m+1}^{h})-\E\PsiM(X_{m}^{h})=&\int_{t_m}^{t_{m+1}}\E\Bigl(\mathcal{L}^{\tau,m,h}-\mathcal{L}^{h}\Bigr)\PsiM(\tilde{X}^{h}(s))ds\\
&+\int_{t_m}^{t_{m+1}}\E\Bigl(\mathcal{L}^{h}-\mathcal{L}^{(M)}\Bigr)\PsiM(\tilde{X}^{h}(s))ds\\
&+\int_{t_m}^{t_{m+1}}\E(\phi(P_M\tilde{X}^h(s))-\overline{\phi}_M)ds\\
&+\int_{t_m}^{t_{m+1}}\E\langle P_M\Big(F(\tilde{X}^h(s))-F(P_M\tilde{X}^h(s))\Big),D\Psi^{(M)}(P_M\tilde{X}^h(s))\rangle ds\\
=&\int_{t_m}^{t_{m+1}}\E\Bigl(\mathcal{L}^{\tau,m,h}-\mathcal{L}^{h}\Bigr)\PsiM(\tilde{X}^{h}(s))ds\\
&+\int_{t_m}^{t_{m+1}}\E\Bigl(\mathcal{L}^{h}-\mathcal{L}^{(M)}\Bigr)\PsiM(\tilde{X}^{h}(s))ds\\
&+\tau \big(\E\phi(P_MX_m^{h})-\overline{\phi}_M\big)\\
&+\int_{t_m}^{t_{m+1}}\E\Bigl(\phi(P_M\tilde{X}^{h}(s))\Bigr)-\E\Bigl(\phi(P_MX_m^{h})\Bigr)ds\\
&+\int_{t_m}^{t_{m+1}}\E\langle P_M\Big(F(\tilde{X}^h(s))-F(P_M\tilde{X}^h(s))\Big),D\Psi^{(M)}(P_M\tilde{X}^h(s))\rangle ds
\end{align*}

Then if we sum for $m=1,\ldots,N-1$ and divide by $N\tau$, we obtain

\begin{equation}
\begin{aligned}
\frac{1}{N}\sum_{m=0}^{N-1}\Bigl(\E\phi(P_MX_{m}^{h})-&\overline{\phi}_M\Bigr)=\frac{1}{N\tau}\bigl(\E\Psi^{(M)}(P_MX_{N}^{h})-\E\Psi^{(M)}(P_MX_{1}^{h})\bigr)\\
&+\frac{1}{N}\bigl(\phi(P_Mx)-\overline{\phi}_M\bigr)\\
&+\frac{1}{N\tau}\sum_{m=1}^{N-1}\int_{t_m}^{t_{m+1}}\E\Bigl(\mathcal{L}^{(M)}-\mathcal{L}^{h}\Bigr)\PsiM(\tilde{X}^{h}(s))ds\\
&+\frac{1}{N\tau}\sum_{m=1}^{N-1}\int_{t_m}^{t_{m+1}}\E\Bigl(\mathcal{L}^{h}-\mathcal{L}^{\tau,m,h}\Bigr)\PsiM(\tilde{X}^{h}(s))ds\\
&-\frac{1}{N\tau}\sum_{m=1}^{N-1}\int_{t_m}^{t_{m+1}}\Big(\E\phi(P_M\tilde{X}^{h}(s))-\E\phi(P_MX_m^{h})\Big)ds\\
&-\frac{1}{N\tau}\sum_{m=1}^{N-1}\int_{t_m}^{t_{m+1}}\E\langle P_M\Big(F(\tilde{X}^h(s))-F(P_M\tilde{X}^h(s))\Big),D\Psi^{(M)}(P_M\tilde{X}^h(s))\rangle ds\\ &:=I_1+I_2+I_3+I_4+I_5+I_6.
\end{aligned}
\end{equation}
The two first terms and the last term are easy to control. Indeed, using Proposition \ref{Prop220} and Lemma \ref{lem2}, we get for $0<\tau <\tau_0$,
\begin{equation*}
|I_1+I_2|\leq C(1+|x|^2)\frac{1}{N\tau},
\end{equation*}
where $\tau_0$ is any fixed positive real number. Using the fact that $F$ is lipschitz, Proposition \ref{Prop220} and Lemma \ref{lem2}, we get 
\begin{equation*}
\lim_{M\rightarrow\infty}I_6\rightarrow 0.
\end{equation*}

The control of of the three other terms will be the subject of three following Subsections. First, in Subsection \ref{ssectesti0}, the following estimate of $I_3$ is shown:
\begin{lemme}[Space-discretization error]\label{estiI3}
For any $0<\kappa<1/2$ and $\tau_0$, there exists a constant $C>0$ such that for any $\phi\in \mathcal{C}^2_b(H)$, $x\in H$ and $0<\tau\leq\tau_0$
\begin{equation*}
\limsup_{M\rightarrow\infty} \frac{1}{N\tau}\sum_{m=1}^{N-1}\int_{t_m}^{t_{m+1}}\E\Bigl(\mathcal{L}^{(M)}-\mathcal{L}^{h}\Bigr)\PsiM(\tilde{X}^{h}(s))ds\leq C(1+|x|^3)\parallel \phi \parallel_{2,\infty} h^{1-\kappa}(1+(N\tau)^{-1}).
\end{equation*}
\end{lemme}

In Subsection \ref{ssectesti1}, we will show the following estimate of $I_4$:
\begin{lemme}[Time-discretization error]\label{estiI4}
For any $0<\kappa<1/2$ and $\tau_0$, there exists a constant $C>0$ such that for any $\phi\in\mathcal{C}^2_b(H)$, $M\in\left\{1,2,\ldots\right\}$, $y\in H$ and $0<\tau\leq\tau_0$
\begin{equation*}
|\frac{1}{N\tau}\sum_{m=1}^{N-1}\int_{t_m}^{t_{m+1}}\E\Big(\mathcal{L}^{h}-\mathcal{L}^{\tau,m,h}\Big)\PsiM(\tilde{X}^{h}(t))dt|\leq C\parallel \phi\parallel_{2,\infty}(1+|x|^3)\tau^{1/2-\kappa}(1+(N\tau)^{-1+\kappa}+(N\tau)^{-1}).
\end{equation*}
\end{lemme}

Finally, the proof of the following estimate of $I_5$ is detailed in Subsection \ref{ssectesti2}:
\begin{lemme}[Additional time-discretization error]\label{estiI5}
For any $0<\kappa<1/4$ and $\tau_0$, there exists a constant $C>0$ such that for any $\phi\in\mathcal{C}^2_b(H)$, $M\in\left\{1,2,\ldots\right\}$, $y\in H$ and $0<\tau\leq\tau_0$
\begin{equation*}
|\frac{1}{N\tau}\sum_{m=1}^{N-1}\int_{t_m}^{t_{m+1}}\Big(\E\phi(P_M\tilde{X}^{h}(t))-\E\phi(P_MX_m^{h})\Big)dt| \leq C\parallel \phi\parallel_{2,\infty} \tau^{1/2-2\kappa}\big(1+\frac{|x|}{(N\tau)^{1-\kappa}}\big).
\end{equation*}
\end{lemme}

The order of the respective proofs may not seem natural. We have made the choice to put the proof of the Lemma \ref{estiI4} on the time-discretization error at the end since it essentially uses the same arguments as in \cite{B2}, while the others require new estimates, appearing for the first time in this context of approximations of invariant laws. We thus begin with the proof of Lemma \ref{estiI3} on the space discretization error, and go on with the proof of Lemma \ref{estiI5} concerning the additional time-discretization error induced by the use of the Poisson equation - instead of the the generators, we have the process alone.

\section{Detailed proof of the estimates}\label{sectdetproof}

We warn the reader that constants may vary from line to line during the proofs, and that in order to use lighter notations we usually forget to mention dependence on the parameters. We use the generic notation $C$ for such constants. All constants will depend on a parameter $\kappa>0$, which can be chosen as small as necessary.

To simplify the expressions, we have not mentioned the dependence of the error with respect to the test function $\phi$. However, thanks to Proposition \ref{Prop220} it is straightforward to give this precision.


\subsection{Control of the space dicretization error $\frac{1}{N\tau}\sum_{m=1}^{N-1}\int_{t_m}^{t_{m+1}}\E\Big(\mathcal{L}^{(M)}-\mathcal{L}^{h}\Big)\PsiM(\tilde{X}^{h}(t))dt$.}\label{ssectesti0}
\hspace{5cm}

\vspace{0.5cm}
In this Subsection, we prove the estimate of the Lemma \ref{estiI3}.

\subsubsection{Strategy}
To control this term, we will mix ideas describes in \cite{AndersLars} and in \cite{B2}. In \cite{AndersLars}, to prove weak convergence, the authors use estimation on $u$ the solution of the Kolmogorov equation associated to the SPDE. We can use the same ideas because, for each $M\in\{1,2,\ldots \}$, $\Psi^{(M)}$, the solution of the Poisson equation, and its derivatives verify the same kind of estimations than $u$.

Let $M\in\{1,2,\ldots\}$ be fixed. First, we will decompose $\mathcal{L}^{(M)}-\mathcal{L}^h$ in three terms. We have for any $x\in H$
\begin{align*}
\Big(\mathcal{L}^{(M)}-\mathcal{L}^h\Big)\tilde{\Psi}^{(M)}(x)=&\langle \big(A_M-A_h\big)x,D\tilde{\Psi}^{(M)}(x)\rangle\\
	&+\langle \big(P_M-P_h\big)F(x),D\PsiM(x)\rangle\\
	&+\frac{1}{2}Tr\Big(\big(P_M-P_h\big)D^2\PsiM(x)\Big)
\end{align*}
and we obtain
\begin{equation*}
\frac{1}{N\tau}\sum_{m=1}^{N-1}\int_{t_m}^{t_{m+1}}\E\Big(\mathcal{L}^{(M)}-\mathcal{L}^{h}\Big)\PsiM(\tilde{X}^{h}(t))dt=\frac{1}{N\tau}\sum_{m=1}^{N-1}(a^m+b^m+c^m),
\end{equation*}
where for $1\leq m\leq N-1$
\begin{gather*}
a^m=\E\int_{t_m}^{t_{m+1}}\langle (A_M-A_h)\tilde{X}^h(t),D\PsiM(\tilde{X}^h(t))\rangle dt\\
b^m=\E\int_{t_m}^{t_{m+1}}\langle (P_M-P_h)F(\tilde{X}^h(t)),D\PsiM(\tilde{X}^h(t))\rangle dt\\
c^m=\frac{1}{2}\E\int_{t_m}^{t_{m+1}}Tr\Big( (P_M-P_h)D^2\PsiM(\tilde{X}^h(t))\Big) dt.
\end{gather*}

\subsubsection{Estimate of $a^m$} The Ritz projection $R_h$ can be expressed in the form $R_h=A_h^{-1}P_hA$. Using this we can write
\begin{align*}
\langle(A_M-A_h)&\tilde{X}^h(t),D\PsiM(\tilde{X}^h(t))\rangle=\langle (A_MP_h-A_hP_h)\tilde{X}^h(t),D\PsiM(\tilde{X}^h(t))\rangle\\
	=&\langle \tilde{X}^h(t),(P_hA_M-A_hP_h)D\PsiM(\tilde{X}^h(t))\rangle=\langle \tilde{X}^h(t),A_hP_h(R_hP_M-I)D\PsiM(\tilde{X}^h(t))\rangle\\
	=&\langle \tilde{X}^h(t),A_hP_h(R_h-I)P_M D\PsiM(\tilde{X}^h(t))\rangle+\langle \tilde{X}^h(t),A_hP_h(P_M-I)D\PsiM(\tilde{X}^h(t))\rangle.
\end{align*}
The idea of this decomposition is to apply the error estimates \eqref{estiR_h} and \eqref{estiP_M} for $R_h$ and $P_M$ respectively. We now use formula \eqref{tild}  on $\tilde{X}^h(t)$. We then need to estimate the following five terms.
\begin{align*}
a^m= &\E \int_{t_m}^{t_{m+1}}\langle X^h_m,A_hP_h(R_h-I)P_M D\PsiM(\tilde{X}^h(t))\rangle dt\\
 &+\E \int_{t_m}^{t_{m+1}}(t-t_m)\langle A_hS_{\tau,h} X^h_m,A_hP_h(R_h-I)P_M D\PsiM(\tilde{X}^h(t))\rangle dt\\
&+\E \int_{t_m}^{t_{m+1}}(t-t_m)\langle S_{\tau,h}P_h F( X^h_m),A_hP_h(R_h-I)P_M D\PsiM(\tilde{X}^h(t))\rangle dt\\
&+\E \int_{t_m}^{t_{m+1}}\langle\int_{t_m}^t S_{\tau,h}P_h dW(s),A_hP_h(R_h-I)P_M D\PsiM(\tilde{X}^h(t))\rangle dt\\
&+\E \int_{t_m}^{t_{m+1}}\langle A_hX^h_m,(P_M-I) D\PsiM(\tilde{X}^h(t))\rangle dt\\
=:&a_1^{m,h}+a_2^{m,h}+a_3^{m,h}+a_4^{m,h}+a^{m,M}
 \end{align*}

\begin{enumerate}
\item \textbf{Estimate of $a_1^{m,h}$:}
We use expressions \eqref{exprYk} of $X_m^h$ and \eqref{exprsto} to decompose $a_1^{m,h}$:
\begin{align*}
a_1^{m,h}=&\E\int_{t_m}^{t_{m+1}}\langle S_{\tau,h}^mP_hx,A_hP_h(R_h-I)P_MD\PsiM(\tilde{X}^h(t))\rangle dt\\
	&+\E\int_{t_m}^{t_{m+1}}\tau\sum_{\ell=0}^{k-1}\langle S_{\tau,h}^{m-\ell}P_hF(X^h_{\ell}),A_hP_h(R_h-I)P_MD\PsiM(\tilde{X}^h(t))\rangle dt\\
	&+\E\int_{t_m}^{t_{m+1}}\langle\int_0^{t_m} S_{\tau,h}^{m-l_s}P_hdW(s),A_hP_h(R_h-I)P_MD\PsiM(\tilde{X}^h(t))\rangle dt\\
	=:&a_{1,1}^{m,h}+a_{1,2}^{m,h}+a_{1,3}^{m,h}.
\end{align*}

\begin{itemize}
\item \textit{Estimate of $a_{1,1}^{m,h}$:}
The ideas are to "share" $(-A_h)$ between different factors and to use regularization properties of the semi-group $(S^k_{\tau,h})_{k\in\mathbb{N}}$. Thanks to Proposition \ref{conv_FEM}, Proposition \ref{Prop220} for $\beta=\frac{1}{2}$, Lemma \ref{lem2} and Lemma \ref{lem6}, we get, for any small enough parameter $0<\kappa<1/2$,
\begin{align*}
|a_{1,1}^{m,h}|=&|\E\int_{t_m}^{t_{m+1}}\langle (-A_h)^{1-\kappa}S_{\tau,h}^mP_hx,(-A_h)^{\kappa}P_h(R_h-I)(-A)^{-1/2}P_M(-A)^{1/2}D\PsiM(\tilde{X}^h(t))\rangle dt|\\
	\leq &\E\int_{t_m}^{t_{m+1}}|(-A_h)^{1-\kappa}S_{\tau,h}^mP_h|_{\mathcal{L}(H)}|x||(-A_h)^{\kappa}P_h(R_h-I)(-A)^{-1/2}|_{\mathcal{L}(H)}\\
	&|P_M|_{\mathcal{L}(H)}|(-A)^{1/2}D\PsiM(\tilde{X}^h(t))|dt\\
	\leq&C\frac{1}{(m\tau)^{1-\kappa}}\frac{1}{(1+\lambda_0\tau)^{m\kappa}}|x||(-A)^{\kappa}(R_h-I)(-A)^{-1/2}|_{\mathcal{L}(H)}\int_{t_m}^{t_{m+1}}\E(1+|\tilde{X}^h(t)|^2)dt\\
	\leq&	C\tau\frac{1}{(m\tau)^{1-\kappa}}\frac{1}{(1+\lambda_0\tau)^{m\kappa}}(1+|x|^3)h^{1-2\kappa}.	
\end{align*}

We will now use the following useful inequality:  for $\tau\leq \tau_0$ and any $N\geq 1$
\begin{equation}\label{sumesti}
\tau\sum_{l=1}^{N}\frac{1}{(l\tau)^{1-\kappa}}\frac{1}{(1+\lambda_0\tau)^{l\kappa}}\leq C_\kappa.
\end{equation}

Indeed,
\begin{align*}
\tau\sum_{l=1}^{N}\frac{1}{(l\tau)^{1-\kappa}}\frac{1}{(1+\lambda_0\tau)^{l\kappa}}&\leq C\int_{0}^{t_N}\frac{1}{t^{1-\kappa}}\frac{1}{(1+\lambda_0\tau)^{\kappa\frac{t}{\tau}}}dt\\
&\leq \int_{0}^{\infty}\frac{1}{t^{1-\kappa}}e^{-t\frac{\kappa}{\tau}\log(1+\lambda_0\tau)}dt\\
&\leq \int_{0}^{\infty}\frac{1}{s^{1-\kappa}}e^{-s}ds \left(\frac{\tau}{\kappa\log(1+\lambda_0\tau)}\right)^\kappa\\
&\leq C_\kappa.
\end{align*}
Then, using \eqref{sumesti}, we get
\begin{equation}\label{ea11}
\frac{1}{N\tau}\sum_{m=1}^{N-1}|a_{1,1}^{m,h}|\leq C\frac{1}{T}h^{1-2\kappa}(1+|x|^{3}).
\end{equation}

\item \textit{Estimate of $a_{1,2}^{m,h}$:}
Using the same ideas than to estimate $a_{1,1}^{m,h}$, we have
\begin{multline*}
|a_{1,2}^{m,h}|\leq C\tau\E\int_{t_{m}}^{t_{m+1}}\sum_{l=0}^{m-1}|(-A_h)^{1-\kappa}S_{\tau,h}^{m-l}P_hF(X^h_l)||(-A_h)^{\kappa}P_h(R_h-I)(-A)^{-1/2}|_{\mathcal{L}(H)}\\
\times|(-A)^{1/2}D\PsiM(\tilde{X}^h(t))|dt.
\end{multline*}

Since $F$ is supposed to be bounded, the estimate \eqref{sumesti} yields
$$
|\tau(-A_h)^{1-\kappa}\sum_{l=0}^{m-1}S_{\tau,h}^{m-l}F(X^h_l)|\leq C\|F\|_{\infty}\tau\sum_{l=1}^{m}\frac{1}{(l\tau)^{1-\kappa}}\frac{1}{(1+\lambda_0\tau)^{l\kappa}}\leq C_\kappa.
$$

With Lemma \ref{conv_FEM}, Lemma \ref{lem2} and Proposition \ref{Prop220} for $\beta=\frac{1}{2}$, we can now write
$$|a_{1,2}^{m,h}|\leq C(1+|x|^2)h^{1-2\kappa}\tau,$$
and we get
\begin{equation}\label{ea12}
\frac{1}{N\tau}\sum_{m=0}^{N-1}|a_{1,2}^{m,h}|\leq C(1+|x|^2)h^{1-2\kappa}.
\end{equation}

\item \textit{Estimate of $a_{1,3}^{m,h}$:}
The analysis of this term is more complicated. We refer the reader to \cite{B2} for a discussion on the problem, and for detailed explications on the strategy of the proof - following the original idea of \cite{deb}.

We recall that the problem lies in the regularity in space of the process, due to whiteness in space of the driving noise. The strategy used to control $a^{m,h}_{1,1}$ and $a^{m,h}_{2,1}$ above would only give an order of convergence $1/4$, while we expect $1/2$ - our constants need to be uniform with respect to the mesh size $h$!

We decompose $a^{m,h}_{1,3}$ into two parts, corresponding to different intervals for the stochastic integration. On one of these parts, we can work directly. On the other, a Malliavin integration by parts is performed: it allows to use appropriate regularization properties, and to obtain the correct order of convergence $1/2$. We emphasize on the length of the interval where this integration by parts is applied: its maximal size is independent of $\tau$ and $h$, so that the possible exponential divergence when time increases of the Malliavin derivatives implies no trouble.

By using \eqref{exprsto}, we make the decomposition
\begin{align*}
a_{1,3}^{m,h}&=\E\int_{t_m}^{t_{m+1}}\langle \int_0^{t_m} S_{\tau,h}^{m-l_s}P_hdW(s),(-A_h)P_h(R_h-I)P_MD\PsiM(\tilde{X}^h(t))\rangle dt\\
&=\E\int_{t_m}^{t_{m+1}}\langle\int_0^{(t_m-3\tau_0)\vee 0} (-A_h)^{1-\kappa}S_{\tau,h}^{m-l_s}P_hdW(s),(-A_h)^{\kappa}P_h(R_h-I)P_MD\PsiM(\tilde{X}^h(t))\rangle dt\\
&+\E\int_{t_m}^{t_{m+1}}\langle\int_{(t_m-3\tau_0)\vee 0}^{t_m} P_M(R_h-I)P_h(-A_h)S_{\tau,h}^{m-l_s}P_hdW(s),D\PsiM(\tilde{X}^h(t))\rangle dt.
\end{align*}

For the first term - which is equal to $0$ when $t_m<3\tau_0$ - we use the Cauchy-Schwarz inequality and we directly get
\begin{multline*}
|\E\langle\int_{0}^{(t_m-3\tau_0)\vee 0}(-A_h)^{1-\kappa}S_{\tau,h}^{m-l_s}P_HdW(s),(-A_h)^{\kappa}P_h(R_h-I)P_MD\PsiM(\tilde{X}^h(t))\rangle|\\
\leq (\E|(-A_h)^{\kappa}P_h(R_h-I)(-A)^{-1/2}P_M(-A)^{1/2}D\PsiM(\tilde{X}^h(t))|^2)^{1/2}\\
\times (\E|\int_{0}^{(t_m-3\tau_0)\vee 0}(-A_h)^{1-\kappa}S_{\tau,h}^{m-l_s}P_hdW(s)|^2)^{1/2}.
\end{multline*}

We have the following inequality - we remark that in the integral below $t_{m-l_s}\geq 1$:
\begin{align*}
\E|\int_{0}^{(t_m-3\tau_0)\vee 0}(-A_h)^{1-\kappa}S_{\tau,h}^{m-l_s}P_hdW(s)|^2=&\int_{0}^{(t_m-3\tau_0)\vee 0}|(-(A_h)^{1-\kappa}S_{\tau,h}^{m-l_s}P_h|_{\mathcal{L}_2(H)}^{2}ds\\
=&\int_{0}^{(t_m-3\tau_0)\vee 0}\text{Tr}((-A_h)^{2-2\kappa}S_{\tau,h}^{2(m-l_s)}P_h)ds\\
\leq& \int_{0}^{(t_m-3\tau_0)\vee 0}|S_{\tau,h}^{(m-l_s)}P_h|_{\mathcal{L}(H)} |(-A_h)^{2+1/2+\kappa}S_{\tau,h}^{(m-l_s)}P_h|_{\mathcal{L}(H)}ds\\&\times\text{Tr}(P_h(-A_h)^{-1/2-\kappa}P_h)\\
\leq &C\int_{0}^{(t_m-3\tau_0)\vee 0}\frac{1}{(1+\lambda_0\tau)^{m-l_s}t_{m-l_s}^{2+1/2-\kappa}}ds\\
\leq &C\int_{0}^{(t_m-3\tau_0)\vee 0}\frac{1}{(1+\lambda_0\tau)^{m-l_s}}ds\\
\leq &C\int_{0}^{+\infty}\frac{1}{(1+\lambda_0\tau)^{s/\tau}}ds\\
\leq& C,
\end{align*}
when $\tau\leq \tau_0$ and thanks to Proposition \ref{propoTrace} and Lemma \ref{lem6} .
Then, thanks to Proposition \ref{conv_FEM}, Proposition \ref{Prop220} for $\beta=\frac{1}{2}$ and Lemma \ref{lem2}, we get
\begin{multline*}
|\E\int_{t_m}^{t_{m+1}}\langle\int_{0}^{(t_m-3\tau_0)\vee 0}(-A_h)^{1-\kappa}S_{\tau,h}^{m-l_s}P_hdW(s),(-A_h)^{\kappa}P_h(R_h-I)P_MD\PsiM(\tilde{X}^h(t))\rangle dt|\\
\leq C(1+|x|^2)\tau h^{1-2\kappa}.
\end{multline*}

For the second term, we use the Malliavin integration by parts formula (Lemma \ref{lemintbyparts}) to get
\begin{align*}
\E\int_{t_m}^{t_{m+1}}\langle&\int_{(t_m-3\tau_0)\vee 0}^{t_m} P_M(R_h-I)P_h(-A_h)S_{\tau,h}^{m-l_s}P_hdW(s),D\PsiM(\tilde{X}^h(t))\rangle dt\\
&=\E\int_{t_m}^{t_{m+1}}\int_{(t_m-3\tau_0)\vee 0}^{t_m}\text{Tr}\left(S_{\tau,h}^{m-l_s}(-A_h)P_h(R_h-I)P_M D^2\PsiM(\tilde{X}^h(t))\DM_s\tilde{X}^h(t)\right)dsdt.
\end{align*}

Thanks to both estimates of Lemma \ref{lem5}, we have for $(t_m-3\tau_0)\vee 0\leq s\leq t_m\leq t< t_{m+1}$
$$|(-A)^\alpha \DM_{s}^{\ell}\tilde{X}^h(t)|\leq C(1+L_F\tau)^{m-l_s}(1+\frac{1}{(1+\lambda_0\tau)^{(1-\alpha)(m-l_s)}t_{m-l_s}^{\alpha}}),$$
and we see that $(1+L_F\tau)^{m-l_s}$ is bounded by a constant.

We can then control the second term of $a_{3,1}^{m,h}$ with
\begin{multline*}
\E\int_{t_{m}}^{t_{m+1}}\int_{(t_m-3\tau_0)\vee 0}^{t_m}|(-A_h)^{1-3\frac{\kappa}{2}}S_{\tau,h}^{m-l_s}|_{\mathcal{L}(H)}|(-A_h)^{3\frac{\kappa}{2}}P_h(R_h-I)(-A)^{-1/2}|_{\mathcal{L}(H)}\\
\times|(-A)^{1/2}D^2\PsiM(\tilde{X}^h(t))(-A)^{1/2-\kappa/2}|_{\mathcal{L}(H)}\text{Tr}((-A)^{-1/2-\kappa/2})|(-A)^\kappa \DM_s\tilde{X}^h(t)|_{\mathcal{L}(H)}dsdt\\
\leq C
\int_{(t_m-3\tau_0)\vee 0}^{t_m}t_{m-l_s}^{-1+3\frac{\kappa}{2}}\frac{1}{(1+\lambda_0\tau)^{(m-l_s)3\frac{\kappa}{2}}}\Big(1+t_{m-l_s}^{-\kappa}\frac{1}{(1+\lambda_0\tau)^{(m-l_s)(1-\kappa)}}\Big)ds\\
\times\tau h^{1-3\kappa}(1+|x|^2),
\end{multline*}
using Proposition \ref{Prop220}, Lemmas \ref{lem5} and \ref{lem6}.

We have

$$\int_{(t_m-3\tau_0)\vee 0}^{t_m}t_{m-l_s}^{-1+3\frac{\kappa}{2}}\frac{1}{(1+\lambda_0\tau)^{(m-l_s)3\frac{\kappa}{2}}}ds\leq\int_{0}^{t_m}\frac{1}{s^{1-3\frac{\kappa}{2}}}\frac{1}{(1+\lambda_0\tau)^{3\frac{\kappa}{2} s/\tau}}ds\leq C<+\infty,$$
for $\tau\leq\tau_0$, thanks to \eqref{sumesti}.

Therefore
\begin{equation}\label{ea13}
\frac{1}{N\tau}\sum_{m=1}^{N-1}|a^{m,h}_{1,3}|\leq C(1+|x|^2)h^{1-3\kappa}.
\end{equation}

\end{itemize}

Using \eqref{ea11}, \eqref{ea12} and \eqref{ea13}, we have
\begin{equation}\label{estia1mh}
\frac{1}{N\tau}\sum_{m=1}^{N-1}|a_{1}^{m,h}|\leq C(1+\frac{1}{T})h^{1-3\kappa}(1+|x|^{3}).
\end{equation}

\item \textbf{Estimate of $a_2^{m,h}$:}
Since $(t-t_m)|(-A_h)S_{\tau,h}|_{\mathcal{L}(H)}\leq C$, $a_2^{m,h}$ is bounded by the same expression as $a_1^{m,h}$: by \eqref{estia1mh}, we have
\begin{equation}\label{estia2mh}
\frac{1}{N\tau}\sum_{m=1}^{N-1}|a_{2}^{m,h}|\leq C(1+\frac{1}{T})h^{1-3\kappa}(1+|x|^{3}).
\end{equation}

\item \textbf{Estimate of $a_3^{m,h}$:}
We have
\begin{multline*}
|a_3^{m,h}|\leq \E\int_{t_m}^{t_{m+1}}(t-t_m)|(-A_h)^{1-\kappa}S_{\tau,h}P_h|_{\mathcal{L}(H)}|F(X_{m}^h)|\\
|(-A_h)^{\kappa}P_h(R_h-I)(-A)^{-1/2}P_M(-A)^{1/2}D\PsiM(\tilde{X}^h(t))|dt
\end{multline*}
Since $F$ and $(t-t_m)|(-A_h)^{1-\kappa}S_{\tau,h}P_h|_{\mathcal{L}(H)}$ are bounded, using the same idea than to estimate $a_{1,1}^{m,h}$, we get
\begin{equation}\label{estia3mh}
\frac{1}{N\tau}\sum_{m=1}^{N-1}|a_{3}^{m,h}|\leq \frac{C}{T}h^{1-2\kappa}(1+|x|^{3}).
\end{equation}

\item \textbf{Estimate of $a_4^{m,h}$:}
We again use the integration by parts formula to rewrite $a_{4}^{m,h}$:
\begin{align*}
a_{4}^{m,h}&=-\E\int_{t_m}^{t_{m+1}}\langle\int_{t_m}^{t}S_{\tau,h}P_h dW(s),(-A_h)P_h(R_h-I)D\PsiM(\tilde{X}^h(t))\rangle dt\\
&=-\E\int_{t_{m}}^{t_{m+1}}\int_{t_m}^{t}\text{Tr}(S_{\tau,h}P_h(-A_h)P_h(R_h-I)\PsiM(\tilde{X}^h(t))\DM_s\tilde{X}^h(t))dsdt.
\end{align*}
From \eqref{tild}, for $t_m\leq s\leq t\leq t_{m+1}$ we have $\DM_{s}^{\ell}\tilde{X}^h(t)=S_{\tau,h} P_h\ell$; as a consequence, the situation is much simpler and we do not need to use the same trick as in the control of $a_{1,3}^{m,h}$.

Then, as previously, we have
\begin{align*}
|a_{4}^{m,h}|\leq&\E\int_{t_m}^{t_{m+1}}(t-t_m)\text{Tr}((-A_h)^{1-\kappa}S_{\tau,h}P_h (-A_h)^{\kappa}P_h(R_h-I)(-A)^{-1/2}P_M\\
	&\times (-A)^{1/2}  D^2\PsiM( \tilde{X}^h(t))(-A)^{1/2-\kappa/2}(-A)^{-1/2-\kappa/2}(-A)^{\kappa}S_{\tau,h})dt\\
\leq& c|(-A_h)^{1-\kappa}S_{\tau,h}P_h|_{\mathcal{L}(H)}\text{Tr}((-A)^{-1/2-\kappa/2})|(-A_h)^{\kappa}P_h(R_h-I)(-A)^{-1/2}|_{\mathcal{L}(H)}\\
&\times\E\int_{t_{m}}^{t_{m+1}}|(-A)^{1/2}D^2\PsiM(\tilde{X}^h(t))(-A)^{1/2-\kappa/2}|_{\mathcal{L}(H)}|(-A)^{\kappa}S_{\tau,h}P_h|_{\mathcal{L}(H)}dt\\
\leq& c(1+|x|^2)\tau h^{1-2\kappa}.
\end{align*}
Therefore
\begin{equation}\label{e7_2}
\frac{1}{N\tau}\sum_{m=1}^{N-1}|a_{4}^{m,h}|\leq C(1+|x|^2) h^{1-2\kappa}.
\end{equation}

\item \textbf{Estimate of $a^{m,M}$:}
Using Proposition \ref{Prop220}, Lemma \ref{lem2} and estimate \eqref{estiP_M}, we have
\begin{align*}
|a^{m,M}|&\leq\int_{t_m}^{t_{m+1}}\E\Big(|(-A_h)P_h|_{\mathcal{L}(H)}|\tilde{X}^h(t)||(P_M-I)(-A)^{-1/2+\kappa}||(-A)^{1/2+\kappa}D\PsiM(\tilde{X}^h(t))|\Big)dt \\
 	&\leq C_h\parallel \phi\parallel_{1,\infty}\lambda_M^{-1/2+\kappa}\int_{t_m}^{t_{m+1}}\E\Big(|\tilde{X}^h(t)|(1+|\tilde{X}^h(t)|^2)\Big)dt\\
	&\leq C_h\parallel\phi\parallel_{1,\infty}\lambda_M^{-1/2+\kappa}\tau(1+|x|^3)
\end{align*}
Then, we get
\begin{equation*}
\lim_{M\rightarrow\infty}\frac{1}{N\tau}\sum_{m=1}^{N-1}|a^{m,M}|=0
\end{equation*}
\end{enumerate}

With the previous estimates, we get
\begin{equation}\label{eaa}
\limsup_{M\rightarrow +\infty}\frac{1}{N\tau}\sum_{m=1}^{N-1}|a^{m}|\leq C\parallel \phi \parallel_{1,\infty}(1+|x|^3)(1+T^{-1})h^{1-3\kappa}.
\end{equation}

\subsubsection{Estimate of $b^m$} 
Writing $P_M-P_h=(P_M-I)+(I-P_h)$, we get the natural decomposition
\begin{align*}
b^m=&\E\int_{t_m}^{t_{m+1}}\langle (P_M-I)F(P_M\tilde{X}^h(t)),D\PsiM(\tilde{X}^h(t))\rangle dt\\
	&+\E\int_{t_m}^{t_{m+1}}\langle (I-P_h)F(P_M\tilde{X}^h(t)),D\PsiM(\tilde{X}^h(t))\rangle dt\\
	=:&b^{m,M}+b^{m,h}.
\end{align*}
Using the fact that $F$ is bounded, Proposition \ref{Prop220} and Lemma \ref{lem1}, we have
\begin{align*}
|b^{m,i}|=&|\E\int_{t_{m}}^{t_{m+1}}\langle F(P_M\tilde{X}^h(t)),(P_i-I)(-A)^{-1/2+\kappa}(-A)^{1/2-\kappa}D\PsiM(\tilde{X}^h(t))\rangle dt|\\
	\leq&\int_{t_m}^{t_{m+1}}\parallel F\parallel_{\infty}|(P_i-I)(-A)^{-1/2+\kappa}|_{\mathcal{L}(H)}\E|(-A)^{1/2-\kappa}D\PsiM(\tilde{X}^h(t))|dt\\
	\leq& C\tau (1+|x|^2)|(P_i-I)(-A)^{-1/2+\kappa}|_{\mathcal{L}(H)}.	
\end{align*}
Using \eqref{estiP_h} and \eqref{estiP_M}, we get
\begin{equation*}
|b^{m,h}|\leq C\tau (1+|x|^2)h^{1-2\kappa}
\end{equation*}
and
\begin{equation*}
|b^{m,M}|\leq C\tau (1+|x|^2)\lambda_M^{-1/2+\kappa}.
\end{equation*}
Finally, we have
\begin{equation*}
\frac{1}{N\tau}\sum_{m=1}^{N-1}|b^m|\leq C(1+|x|^2)(h^{1-2\kappa}+\lambda_M^{-1/2+\kappa})
\end{equation*}	
and
\begin{equation*}
\limsup_{M\rightarrow +\infty}\frac{1}{N\tau}\sum_{m=1}^{N-1}|b^m|\leq C(1+|x|^2)h^{1-2\kappa}.
\end{equation*}

\subsubsection{Estimate of $c^m$}
We use the same natural decomposition than for $b^m$: $c^m=c^{m,h}+c^{m,M}$, where for $i\in\{h,M\}$
\begin{multline*}
2|c^{m,i}|=|\E\int_{t_m}^{t_{m+1}}Tr\Big((-A)^{2\kappa}(P_i-I)(-A)^{-1/2+\kappa}(-A)^{1/2-\kappa}D^2\PsiM(\tilde{X}^h(t))(-A)^{1/2-\kappa}(-A)^{-1/2-\kappa}\Big)dt|\\
	\leq Tr((-A)^{-1/2-\kappa})|(-A)^{2\kappa}(P_i-I)(-A)^{-1/2+\kappa}|_{\mathcal{L}(H)}\\
	\times\int_{t_m}^{t_{m+1}}\E|(-A)^{1/2-\kappa}D^2\PsiM(\tilde{X}^h(t))(-A)^{1/2-\kappa}|dt.\end{multline*}
Using Assumptions \ref{hypB}, Proposition \ref{Prop220}, Lemma \ref{lem1}, commutativity of $A$ and $P_M$ and estimates \eqref{estiP_M} and \eqref{estiP_h}, we get
\begin{equation*}
2|c^{m,h}|\leq C\tau(1+|x|^2)\lambda_{M}^{-1/2+3\kappa}
\end{equation*}
and
\begin{equation*}
2|c^{m,M}|\leq C\tau(1+|x|^2)h^{1-6\kappa}.
\end{equation*}
Then, we have
\begin{equation*}
\frac{1}{N\tau}\sum_{m=1}^{N-1}|c^m|\leq C(1+|x|^2)(h^{1-6\kappa}+\lambda_M^{-1/2+3\kappa})
\end{equation*}
and
\begin{equation*}
\limsup_{M\rightarrow+\infty}\frac{1}{N\tau}\sum_{m=1}^{N-1}|c^m|\leq C(1+|x|^2)h^{1-6\kappa}.
\end{equation*}

\subsubsection{Conclusion}
With the above estimation, we get
\begin{equation}\label{ef3}
\limsup_{M\rightarrow\infty} \frac{1}{N\tau}\sum_{m=1}^{N-1}\int_{t_m}^{t_{m+1}}\E\Bigl(\mathcal{L}^{(M)}-\mathcal{L}^{h}\Bigr)\PsiM(\tilde{X}^{h}(s))ds\leq C(1+|x|^3) h^{1-\kappa}(1+T^{-1}).
\end{equation}

\subsection{Control of the additional time-discretization error $\E\phi(P_M\tilde{X}(t))-\E\phi(P_MX_{m}^{h})$, if $t_{m}\leq t<t_{m+1}$.}\label{ssectesti2}

\hspace{5cm}

\vspace{0.5cm}

In this Subsection, we prove the estimate of Lemma \ref{estiI5}.


This part of the error is due to the replacement of the continuous-time process $\tilde{X}$ with the discrete-time process from which it is built by interpolation.

If we compare with the other parts of the error, we observe that instead of $\Psi$ the expression involves the test function $\phi$. Since $\phi$ is only assumed to be of class $\mathcal{C}^2_b$, its derivatives do not satisfy estimates with a regularization in space like for $\Psi$. However, we are still able to distribute appropriately the powers of the operator $-A_h$, thus obtaining the good rate of convergence.

We define an auxiliary function $\tilde{\phi}_M:H\rightarrow \R$ with $\tilde{\phi}_M=\phi\circ P_M$. It is of class $\mathcal{C}_{b}^{2}$ and using the identifications introduced in Remark \ref{remident} we have for any $x\in H$
\begin{gather*}
D\tilde{\phi}_M(x)=P_MD\phi(P_Mx),\\
D^2\tilde{\phi}_M(x)=P_MD^2\phi(P_Mx)P_M.
\end{gather*}

Thanks to the It\^o's formula, from \eqref{tild} we get for $t_{m}\leq t<t_{m+1}$
\begin{align*}
\E\phi(P_M\tilde{X}^{h}(t))-\E\phi(P_MX_{m}^{h})&=\E\tilde{\phi}_M(\tilde{X}^{h}(t))-\E\tilde{\phi}_M(\tilde{X}(t_m))\\
&=\E\int_{t_m}^{t}\langle S_{\tau,h}A_hX_{m}^{h},D\tilde{\phi}_M(\tilde{X}^h(s))\rangle ds\\
&+\E\int_{t_m}^{t}\langle S_{\tau,h}P_hF(X_{m}^{h}),D\tilde{\phi}_M(\tilde{X}^h(s))\rangle ds\\
&+\E\int_{t_m}^{t}\frac{1}{2}\text{Tr}((S_{\tau,h}P_h)(S_{\tau,h}P_h)^{*}D^2\tilde{\phi}_M(\tilde{X}^h(s))ds.\\
&=:E_{1}(t)+E_{2}(t)+E_{3}(t).
\end{align*}
The error is naturally divided into three terms. We first treat the easiest ones: $E_2$ and $E_3$.

Using boundedness of the linear operator $S_{\tau,h}$, of the nonlinear coefficient $F$, of the orthogonal projectors $P_M$ and $P_h$ and of the first-order derivative of $\phi$, we easily obtain that for $t_{m}\leq t<t_{m+1}$
$$|E_{2}(t)|=|\E\int_{t_m}^{t}\langle S_{\tau,h}P_hF(X_{m}^{h}),D\tilde{\phi}_M(\tilde{X}^h(s))\rangle ds|\leq C\tau.$$

We now control $E_3(t)$. Using the boundedness of the second-order derivative of $\tilde{\phi}_M$, uniformly with respect to $M$, we have
\begin{align*}
|E_3(t)|&\leq C(t-t_m)\text{Tr}\bigl((S_{\tau,h}P_h)(S_{\tau,h}P_h)^{*}\bigr)\\
&\leq C\tau \text{Tr}\bigl[((-A_h)^{1/2+\kappa}S_{\tau,h}^{2}P_h)P_h(-A_h)^{-1/2-\kappa}P_h\bigr]\\
&\leq C\tau|(-A_h)^{1/2+\kappa}S_{\tau,h}^{2}P_h|_{\mathcal{L}(H)}\text{Tr}\bigl(P_h(-A_h)^{-1/2-\kappa}P_h\bigr)\\
&\leq C\tau^{1/2-\kappa},
\end{align*}
where $\kappa\in(0,1/2)$ is a small parameter, thanks to the first inequality of Lemma \ref{lem6} and to Proposition \ref{propoTrace}.

The treatment of the $E_1$ is the most complicated amongst the three terms, due to the presence of the unbounded operator $A_h$. We recall that $X_{m}^{h}$ is controlled in the norm of $(-A_h)^\alpha$, uniformly in $h$, only for $\alpha<1/4$; to obtain the correct weak order of convergence $1/2$ with respect to $\tau$, we need a careful control. One of the ingredients is the Malliavin integration by parts.

Thanks to \eqref{exprYk} and \eqref{exprsto}, $E_1$ is divided into three parts: $E_1(t)=E_{1,1}(t)+E_{1,2}(t)+E_{1,3}(t)$, with for $t_{m}\leq t<t_{m+1}$
\begin{gather*}
E_{1,1}(t)=\E\int_{t_m}^{t}\langle S_{\tau,h}^{m+1}A_hP_hx,D\tilde{\phi}_M(\tilde{X}^h(s))\rangle ds\\
E_{1,2}(t)=\E\int_{t_m}^{t}\langle \tau A_hS_{\tau,h}\sum_{k=0}^{m-1}S_{\tau,h}^{m-k}P_hF(X_{m}^{h}),D\tilde{\phi}_M(\tilde{X}^h(s))\rangle ds\\
E_{1,3}(t)=\E\int_{t_m}^{t}\langle A_hS_{\tau,h}\int_{0}^{t_m}S_{\tau,h}^{m-l_r}P_hdW(r),D\tilde{\phi}_M(\tilde{X}^h(s))\rangle ds.
\end{gather*}
We have isolated the stochastic part in $X_{m}^{h}$; then only the treatment of $E_{1,3}(t)$ is difficult.

First, using Lemma \ref{lem6}, we have if $m\geq 1$
\begin{align*}
|A_hS_{\tau,h}^{m+1}P_hx|&\leq |(-A_h)^\kappa S_{\tau,h}P_h|_{\mathcal{L}(H)}|(-A_h)^{1-\kappa}S_{\tau,h}^{m}P_h|_{\mathcal{L}(H)}|x|_{H}\\
&\leq C|x|_{H}\tau^{-\kappa}t_{m}^{-1+\kappa}.
\end{align*}
As a consequence, for $t_{m}\leq t<t_{m+1}$
$$|E_{1,1}(t)|\leq C|x|\frac{\tau^{1-\kappa}}{t_{m}^{1-\kappa}}.$$
The treatment of $E_{1,2}$ is similar: we have when $m\geq 1$
\begin{align*}
|\tau A_hS_{\tau,h}\sum_{k=0}^{m-1}S_{\tau,h}^{m-k}P_hF(X_{k}^{h})|&\leq C\tau|(-A_h)^\kappa S_{\tau,h}P_h|_{\mathcal{L}(H)}\sum_{k=0}^{m-1}|(-A_h)^{1-\kappa}S_{\tau,h}^{m-k}P_h|_{\mathcal{L}(H)}|F(X_{k}^{h})|_{H}\\
&\leq C\tau^{-\kappa}\tau\sum_{k=0}^{m-1}|(-A_h)^{1-\kappa}S_{\tau,h}^{m-k}P_h|_{\mathcal{L}(H)},
\end{align*}
$F$ being bounded. Now using Lemma \ref{lem6} and inequality \eqref{sumesti}
we obtain for $m\geq 1$ and $t_{m}\leq t<t_{m+1}$
$$|E_{1,2}(t)|\leq C\tau^{-\kappa}(t-t_m)\leq C\tau^{1-\kappa}.$$

It remains to control $E_{1,3}(t)$, which contains the stochastic term, with low regularity properties. We need to use a Malliavin integration by parts formula; however due to the weak dissipativity condition the behavior of the Malliavin derivatives is bad with respect to time. The solution is to split the stochastic integral factor into two parts: for any $t_m\leq s\leq t<t_{m+1}$
\begin{align*}
\E\langle A_hS_{\tau,h}\int_{0}^{t_m}S_{\tau,h}^{m-l_r}dW(r),D\tilde{\phi}_M(\tilde{X}^h(s))\rangle &=\E\langle A_hS_{\tau,h}\int_{0}^{(t_m-3\tau_0)\vee 0}S_{\tau,h}^{m-l_r}P_hdW(r),D\tilde{\phi}_M(\tilde{X}^h(s))\rangle \\
&+\E\langle A_hS_{\tau,h}\int_{(t_m-3\tau_0)\vee 0}^{t_m}S_{\tau,h}^{m-l_r}P_hdW(r),D\tilde{\phi}_M(\tilde{X}^h(s))\rangle \\
&=:E_{1,3,1}(s,t)+E_{1,3,2}(s,t).
\end{align*}

For the first error term, we directly use the Cauchy-Schwarz inequality and we have (see term $a_{1,3}^{m,h}$ of Sub-Section \ref{ssectesti0} for more details)
\begin{align*}
|E_{1,3,1}(s,t)|^2&\leq C(E|\int_{0}^{(t_m-3\tau_0)\vee 0}S_{\tau,h} A_hS_{\tau,h}^{m-l_r}P_hdW(r)|^2)(\E|D\tilde{\phi}_M(\tilde{X}^h(s))|^2)\\
&\leq C\int_{0}^{(t_m-3\tau_0)\vee 0}\text{Tr}\bigl(P_hA_{h}S_{\tau,h}^{(m-l_r)+1}S_{\tau,h}^{(m-l_r)+1}A_{h}P_h\bigr)dr\\
&\leq C\int_{0}^{(t_m-3\tau_0)\vee 0}\text{Tr}(P_h(-A_h)^{-1/2-\kappa}P_h)|(-A_{h})^{5/2+\kappa}S_{\tau,h}^{2(m-l_r)+1}P_h|_{\mathcal{L}(H)}dr\\
&\leq C.
\end{align*}

For the second error term, using the Malliavin integration by parts formula, we get for any $t_m\leq s\leq t<t_{m+1}$
\begin{align*}
\E\langle A_hS_{\tau,h}\int_{(t_m-3\tau_0)\vee 0}^{t_m}&S_{\tau,h}^{m-l_r}P_hdW(r),D\tilde{\phi}_M(\tilde{X}^h(s))\rangle \\
&=\E\int_{(t_m-3\tau_0)\vee 0}^{t_m}\text{Tr}(S_{\tau,h}^{m-l_r}A_hS_{\tau,h}P_hD^2\tilde{\phi}_M(\tilde{X}^h(s))\DM_r\tilde{X}^h(s))dr.
\end{align*}
We then write that
\begin{align*}
|\E&\int_{(t_m-3\tau_0)\vee 0}^{t_m}\text{Tr}(S_{\tau,h}^{m-l_r}A_hS_{\tau,h}P_hD^2\tilde{\phi}_M(\tilde{X}^h(s)))\DM_r\tilde{X}^h(s)dr|\\
&\leq \int_{(t_m-3\tau_0)\vee 0}^{t_m}\text{Tr}((-A_h)^{\kappa}S_{\tau,h}^{m-l_r}S_{h,\tau}P_h)\E[|(-A_h)^{1-\kappa}\DM_r\tilde{X}^h(s)|_{\mathcal{L}(H)}|D^2\tilde{\phi}_M(\tilde{X}^h(s))|]dr.
\end{align*}

Since $$\text{Tr}((-A_h)^{\kappa}S_{\tau,h}^{m-l_r}S_{h,\tau}P_h)\leq\text{Tr}(P_h(-A_h)^{-1/2-2\kappa}P_h)|S_{\tau,h}^{m-l_r}\bigl((-A_h)^{1/2+2\kappa}S_{h,\tau}P_h\bigr)|_{\mathcal{L}(H)},$$
we have (see term $a_{1,3}^{m,h}$ of Sub-Section \ref{ssectesti0} for more details)
\begin{align*}
|\E\int_{(t_m-3\tau_0)\vee 0}^{t_m}\text{Tr}&(S_{\tau,h}^{m-l_r}A_hS_{\tau,h}P_h\DM_r\tilde{X}^h(s)D^2\tilde{\phi}_M(\tilde{X}^h(s)))dr|\\
&\leq C\tau^{-1/2-2\kappa}\int_{(t_m-3\tau_0)\vee 0}^{t_m}\frac{1}{(1+\lambda_0\tau)^{m-l_r}}(1+L_F\tau)^{m-l_r}(1+\frac{1}{(1+\lambda_0\tau)^{\kappa(m-l_r)}t_{m-l_r}^{1-\kappa}})dr.
\end{align*}
Using that $(1+L_F\tau)^{m-l_r}\leq C$ for the range of $r$ used to compute the integral, we see that
$$|E_{1,3,2}(s,t)|\leq C\tau^{-1/2-2\kappa}.$$
After integration with respect to $s$, we obtain
$$|E_{1,3}(t)|\leq \int_{t_m}^{t}(|E_{1,3,1}(s,t)|+|E_{1,3,2}(s,t)|)ds\leq C(\tau+\tau^{1/2-2\kappa}),$$
and
$$|E_{1}(t)|\leq C(\tau^{1/2-2\kappa}+|x|\frac{\tau^{1-\kappa}}{t_{m}^{1-\kappa}}+\tau^{1-\kappa}).$$

Using the bounds on $E_2$ and $E_3$, we therefore obtain that when $m\geq 1$ and $t_m\leq t\leq t_{m+1}$
$$|\E\phi(P_M\tilde{X}^{h}(t))-\E\phi(P_MX_{m}^{h})|\leq C\tau^{1/2-2\kappa}(1+\frac{|x|}{(m\tau)^{1-\kappa}}).$$
As a consequence, we obtain
\begin{align}\label{ef2}
|\frac{1}{N\tau}\sum_{m=1}^{N-1}&\int_{t_m}^{t_{m+1}}\Big(\E\phi(P_M\tilde{X}^h(t))-\E\phi(P_MX_{m}^{h})\Big)dt|\nonumber\\
&\leq C\tau^{1/2-2\kappa}(1+|x|\frac{1}{N\tau}\int_{0}^{N\tau}\frac{1}{t^{1-\kappa}}dt)\nonumber\\
&\leq C\tau^{1/2-2\kappa}(1+\frac{|x|}{(N\tau)^{1-\kappa}}).
\end{align}

\subsection{Control of the time-discretization error $\frac{1}{N\tau}\sum_{m=1}^{N-1}\int_{t_m}^{t_{m+1}}\E\Big(\mathcal{L}^{h}-\mathcal{L}^{\tau,m,h}\Big)\tilde{\Psi}^{(M)}(\tilde{X}^{h}(t))dt$.}\label{ssectesti1}

\hspace{5cm}

\vspace{0.5cm}

To control this term, we will use ideas described in \cite{B2,deb} and in Sub-Section \ref{ssectesti0}.%
We decompose the error into five terms:
\begin{equation*}
\frac{1}{N\tau}\sum_{m=1}^{N-1}\int_{t_m}^{t_{m+1}}\E\Big(\mathcal{L}^h-\mathcal{L}^{\tau,m,h}\Big)\tilde{\Psi}^{(M)}(\tilde{X}^h(s))ds=\frac{1}{N\tau}\sum_{m=1}^{N-1}(a_m+b_m+c_m+d_m+e_m),
\end{equation*}
where for $1\leq m\leq N-1$
\begin{equation}\label{defabc_k}
\begin{gathered}
a_m=\E\int_{t_m}^{t_{m+1}}\langle (I-S_{\tau,h})A_h X_{m}^{h},D\tilde{\Psi}^{(M)}(\tilde{X}^h(t))\rangle dt,\\
b_m=\E\int_{t_m}^{t_{m+1}}\langle A_h(\tilde{X}^h(t)- X_{m}^{h}),D\tilde{\Psi}^{(M)}(\tilde{X}^h(t))\rangle dt,\\
c_m=\E\int_{t_m}^{t_{m+1}}\langle (I-S_{\tau,h}) P_hF(X_m^h),D\tilde{\Psi}^{(M)}(\tilde{X}^h(t))\rangle dt,\\
d_m=\E\int_{t_m}^{t_{m+1}}\langle F^h(\tilde{X}^h(t))-F^h(X_{m}^{h}),D\tilde{\Psi}^{(M)}(\tilde{X}^h(t))\rangle dt,\\
e_m=\frac{1}{2}\E\int_{t_m}^{t_{m+1}}\text{Tr}\bigl((P_hP_{h}^{*}-(S_{\tau,h}P_h)(S_{\tau,h}P_h)^{*})D^2\tilde{\Psi}^{(M)}(\tilde{X}^h(t))\bigr)dt.
\end{gathered}
\end{equation}

The estimates are:
\begin{lemme}
$$\frac{1}{\tau N}\sum_{m=1}^{N-1}|a_m|\leq C(1+|x|^3)(1+(N\tau)^{-1})\tau^{1/2-2\kappa},$$
$$\frac{1}{N\tau}\sum_{m=1}^{N-1}|b_{m}|\leq C(1+|x|^3)(1+T^{-1})\tau^{1/2-2\kappa},$$
$$\frac{1}{N\tau}\sum_{m=1}^{N-1}|c_{m}|\leq C(1+|x|^2)\tau^{1/2-\kappa},$$
$$\frac{1}{N\tau}\sum_{m=1}^{N-1}|d_m|\leq C\tau^{1/2-2\kappa}(1+|x|^3)(1+T^{-1+\kappa}),$$
$$\frac{1}{\tau N}\sum_{m=1}^{N-1}|e_m|\leq C(1+|x|^2)\tau^{1/2-3\kappa}.$$
\end{lemme}

\subsubsection{Estimate of $a_{m}$}\label{pourint}

We have the equality of linear operators in $\mathcal{L}(V_h)$: $(I-S_{\tau,h})A_h=-\tau S_{\tau,h}A_{h}^2$. Then, using \eqref{exprYk}, we decompose the error into three terms: $a_{m}=a_{m}^{1}+a_{m}^{2}+a_{m}^{3}$, with
\begin{gather*}
a_{m}^{1}=-\tau\E\int_{t_m}^{t_{m+1}}\langle S_{\tau,h} A_h^2S_{\tau,h}^{m}P_hx,P_hD\tilde{\Psi}^{(M)}(\tilde{X}^h(t))\rangle dt\\
a_{m}^{2}=-\tau\E\int_{t_m}^{t_{m+1}}\langle S_{\tau,h} A_h^2\tau\sum_{l=0}^{m-1}S_{\tau,h}^{m-l}F^h(X_{l}^{h}),P_hD\tilde{\Psi}^{(M)}(\tilde{X}^h(t))\rangle dt\\
a_{m}^{3}=-\tau\E\int_{t_m}^{t_{m+1}}\langle S_{\tau,h} A_h^2\sqrt{\tau}\sum_{l=0}^{m-1}S_{\tau,h}^{m-l}P_h\chi_{l+1},P_hD\tilde{\Psi}^{(M)}(\tilde{X}^h(t))\rangle dt;
\end{gather*}

The replacement of $D\tilde{\Psi}^{(M)}(\tilde{X}^h(t))\in H$ with its orthogonal projection $P_hD\tilde{\Psi}^{(M)}(\tilde{X}^h(t))$ is valid since the other factor in the scalar product belongs to $V_h$.

The main task is to control the operator $A_h^2$, using the benefits of the regularization properties of the semi-group $(S_{\tau,h}^{k})_{k\in \N}$ and of the derivatives of $\tilde{\Psi}^{(M)}$. The main difficulties appear in the control of $a_m^3$, where a Malliavin integration by parts is required in order to obtain the correct weak order of convergence. The control of the other terms $a_m^1$ and $a_m^2$ is technical but much easier.

\begin{enumerate}
\item \textbf{Estimate of $a_{m}^{1}$}

We write, for any small enough parameter $0<\kappa<1/2$,
\begin{align*}
|a_{m}^{1}|&\leq \tau\E\int_{t_m}^{t_{m+1}}|S_{\tau,h}(-A_h)^{1/2+2\kappa}P_h|_{\mathcal{L}(H)}|(-A_h)^{1-\kappa}S_{\tau,h}^{m}P_h|_{\mathcal{L}(H)}|P_hx|_{V_h}|(-A_h)^{1/2-\kappa}P_hD\tilde{\Psi}^{(M)}(\tilde{X}^h(t))|_{V_h}dt\\
&\leq C|x|_H\tau \tau^{-1/2-2\kappa}t_{m}^{-1+\kappa}(1+\lambda_0\tau)^{-m\kappa}\int_{t_m}^{t_{m+1}}\E|(-A)^{1/2-\kappa}D\tilde{\Psi}^{(M)}(\tilde{X}^h(t))|_{H}dt\\
&\leq C(1+|x|^3)\tau \tau^{-1/2-2\kappa}t_{m}^{-1+\kappa}(1+\lambda_0\tau)^{-m\kappa},
\end{align*}
thanks to Lemma \ref{lem6}, Lemma \ref{lem2}, Proposition \ref{conv_FEM} and Proposition \ref{Prop220}.

Thanks to estimate \eqref{sumesti}, we get
\begin{equation}\label{e2}
\frac{1}{N\tau}\sum_{m=1}^{N-1}|a_{m}^{1}|\leq C(1+|x|^3)\tau^{1/2-2\kappa}(N\tau)^{-1}.
\end{equation}

\item \textbf{Estimate of $a_{m}^{2}$}

First we write
$$|a_{m}^{2}|\leq C\tau\E\int_{t_{m}}^{t_{m+1}}|S_{\tau,h}(-A_h)^{1/2+2\kappa}|_{\mathcal{L}(V_h)}|\tau(-A_h)^{1-\kappa}\sum_{l=0}^{m-1}S_{\tau,h}^{m-l}P_hF^h(X_l^h)||(-A_h)^{1/2-\kappa}P_hD\tilde{\Psi}^{(M)}(\tilde{X}^h(t))|_{V_h}dt.$$

Since $F$ is supposed to be bounded by $\|F\|_{\infty}$, the estimate \eqref{sumesti} yields
$$
|\tau(-A_h)^{1-\kappa}\sum_{l=0}^{m-1}S_{\tau,h}^{m-l}P_hF(X_l^h)|\leq C\|F\|_{\infty}\tau\sum_{l=1}^{m}\frac{1}{(l\tau)^{1-\kappa}}\frac{1}{(1+\lambda_0\tau)^{l\kappa}}\leq C_\kappa.
$$

With Lemma \ref{lem2} and Proposition \ref{Prop220}, we easily obtain
\begin{equation}\label{e3}
\frac{1}{N\tau}\sum_{m=0}^{N-1}|a_{m}^{2}|\leq C(1+|x|^2)\tau^{1/2-2\kappa}.
\end{equation}

\item \textbf{Estimate of $a_{m}^{3}$}

This is where things become harder. The problem is the same than for estimate $a_{1,3}^{m,h}$ in Sub-Section \ref{ssectesti0}. As for $a^{m,h}_{1,3}$, using \eqref{exprsto}, we decompose $a_{m}^{3}=a_{m}^{3,1}+a_{m}^{3,2}$, with
\begin{gather*}
a_{m}^{3,1}=-\tau\E\int_{t_m}^{t_{m+1}}\langle \int_{0}^{(t_m-5\tau_0)\vee 0}S_{\tau,h} A_h^2S_{\tau,h}^{m-l_s}P_hdW(s),P_hD\tilde{\Psi}^{(M)}(\tilde{X}^h(t))\rangle dt,\\
a_{m}^{3,2}=-\tau\E\int_{t_m}^{t_{m+1}}\langle \int_{(t_m-5\tau_0)\vee 0}^{t_m}S_{\tau,h} A_h^2S_{\tau,h}^{m-l_s}P_hdW(s),P_hD\tilde{\Psi}^{(M)}(\tilde{X}^h(t))\rangle dt.
\end{gather*}

We remark that $a_{m}^{3,1}=0$ if $t_m< 5\tau_0$. Thanks to the Cauchy-Schwarz inequality, to Proposition \ref{Prop220} (for $\beta=0$) and Lemma \ref{lem2}, we get
\begin{align*}
|\E\langle \int_{0}^{(t_m-5\tau_0)\vee 0}&S_{\tau,h} A_h^2S_{\tau,h}^{m-l_s}P_hdW(s),P_hD\tilde{\Psi}^{(M)}(\tilde{X}^h(t))\rangle |\\
&\leq (\E|\int_{0}^{(t_m-5\tau_0)\vee 0}S_{\tau,h} A_h^2S_{\tau,h}^{m-l_s}P_hdW(s)|_H^2)^{1/2}(\E|P_hD\tilde{\Psi}^{(M)}(\tilde{X}^h(t))|_H^2)^{1/2}\\
&\leq C(1+|x|^2);
\end{align*}
indeed we have the following inequality for $0<\tau\leq \tau_0$:
$$\E|\int_{0}^{(t_m-5\tau_0)\vee 0}R_\tau A_h^2S_{\tau,h}^{m-l_s+1}P_hdW(s)|^2\leq C.$$
The proof is easily adapted from the corresponding one in the estimate of $a_{1,3}^{m,h}$ in Sub-Section \ref{ssectesti0}.

Then
$$\frac{1}{N\tau}\sum_{m=0}^{N-1}|a_{m}^{3,1}|\leq  C(1+|x|^2)\tau.$$

For $a_{m}^{3,2}$, we use the integration by parts formula of Lemma \ref{lemintbyparts} to get
\begin{align*}
a_{m}^{3,2}&=-\tau\E\int_{t_m}^{t_{m+1}}\langle \int_{(t_m-5\tau_0)\vee 0}^{t_m}S_{\tau,h} A_h^2S_{\tau,h}^{m-l_s}P_hdW(s),P_hD\tilde{\Psi}^{(M)}(\tilde{X}^h(t))\rangle dt\\
&=-\tau\E\int_{t_m}^{t_{m+1}}\int_{(t_m-5\tau_0)\vee 0}^{t_m}\text{Tr}\left(S_{\tau,h}^{m-l_s}A_h^2S_{\tau,h}P_hD^2\tilde{\Psi}^{(M)}(\tilde{X}^h(t))\DM_s\tilde{X}^h(t)\right)dsdt.
\end{align*}

Thanks to both estimates of Lemma \ref{lem5}, we have for $(t_m-5\tau_0)\vee 0\leq s\leq t_m\leq t< t_{m+1}$
$$|(-A_h)^\beta \DM_{s}\tilde{X}^h(t)|\leq C(1+L_F\tau)^{m-l_s}(1+\frac{1}{(1+\lambda_0\tau)^{(1-\beta)(m-l_s)}t_{m-l_s}^{\beta}}),$$
and we see that $(1+L_F\tau)^{m-l_s}$ is bounded by a constant.

Seen in $\mathcal{L}(H)$, we have $A_{h}^{-1/2-\kappa}P_h=A_{h}^{-1/2}P_h(-A)^{1/2}(-A)^{-1/2-\kappa}(-A)^{\kappa}(-A_{h})^{-\kappa}P_h$; then we use the fact that $\text{Tr}((-A)^{-1/2-\kappa})<+\infty$, and the equivalence of norms of Proposition \ref{conv_FEM}, so that $\sup_{0<h<1}\text{Tr}((-A_h)^{-1/2-\kappa}P_h)<+\infty$; then
\begin{align*}
&|a_{m}^{3,2}|\leq\tau\E\int_{t_{m}}^{t_{m+1}}\int_{(t_m-5\tau_0)\vee 0}^{t_m}|S_{\tau,h} (-A_h)^{1/2+2\kappa}|_{\mathcal{L}(V_h)}|(-A_h)^{1-3\frac{\kappa}{2}}S_{\tau,h}^{m-l_s}|_{\mathcal{L}(V_h)}\text{Tr}((-A)^{-1/2-\frac{\kappa}{2}})\\
&\hspace{1cm}\times|(-A_h)^{1/2-\kappa/2}P_hD^2\tilde{\Psi}^{(M)}(\tilde{X}^h(t))(-A)^{1/2-\kappa/2}|_{\mathcal{L}(H)}|(-A)^\kappa \DM_s\tilde{X}^h(t)|_{\mathcal{L}(H)}dsdt\\
&\leq C\tau^{1/2-2\kappa}\int_{t_{m}}^{t_{m+1}}\int_{(t_m-5\tau_0)\vee 0}^{t_m}t_{m-l_s}^{-1+3\frac{\kappa}{2}}\frac{1}{(1+\lambda_0\tau)^{(m-l_s)3\frac{\kappa}{2}}}\Big(1+t_{m-l_s}^{-\kappa}\frac{1}{(1+\lambda_0\tau)^{(m-l_s)(1-\kappa)}}\Big)dsdt(1+|x|^2)\\
&\leq C\tau\tau^{1/2-2\kappa}(1+|x|^2),
\end{align*}
using estimate \eqref{sumesti}, Assumption \ref{hypB}, Proposition \ref{Prop220}, Lemmas \ref{lem5} and \ref{lem6}.

We obtain $$\frac{1}{N\tau}\sum_{m=1}^{N-1}|a_{m}^{3,2}|\leq C(1+|x|^2)\tau^{1/2-2\kappa}.$$

Therefore
\begin{equation}\label{e4}
\frac{1}{N\tau}\sum_{m=1}^{N-1}|a_{m}^{3}|\leq C(1+|x|^2)\tau^{1/2-2\kappa}.
\end{equation}

\end{enumerate}

With the previous estimates on $a^{1}$, $a^{2}$ and $a^{3}$, we get
\begin{equation}\label{ea}
\frac{1}{N\tau}\sum_{m=1}^{N-1}|a_{m}|\leq C(1+|x|^3)(1+T^{-1})\tau^{1/2-2\kappa}.
\end{equation}

\subsubsection{Estimate of $b_{m}$}

We decompose the corresponding term into three parts: $b_{m}=b_{m}^{1}+b_{m}^{2}+b_{m}^{3}$, with
\begin{gather*}
b_{m}^{1}=\E\int_{t_m}^{t_{m+1}}(t-t_m)\langle A_h(I-S_{\tau,h})X_{m}^h,P_hD\tilde{\Psi}^{(M)}(\tilde{X}^h(t))\rangle dt\\
b_{m}^{2}=\E\int_{t_m}^{t_{m+1}}(t-t_m)\langle A_hS_{\tau,h} F(X_m^h),P_hD\tilde{\Psi}^{(M)}(\tilde{X}^h(t))\rangle dt\\
b_{m}^{3}=\E\int_{t_m}^{t_{m+1}}\langle \int_{t_m}^{t}A_hS_{\tau,h}P_hdW(s),P_hD\tilde{\Psi}^{(M)}(\tilde{X}^h(t))\rangle dt;
\end{gather*}

\begin{enumerate}
\item \textbf{Estimate of $b_{m}^{1}$}

$b_{m}^{1}$ is bounded by the same expression as $a_{m}$: by \eqref{ea} we have
\begin{equation}\label{e5}
\frac{1}{N\tau}\sum_{m=1}^{N-1}|b_{m}^{1}|\leq C(1+|y|^3)(1+T^{-1})\tau^{1/2-2\kappa}.
\end{equation}

\item \textbf{Estimate of $b_{m}^{2}$}

We have 
\begin{align*}
|b_{m}^{2}|&\leq \tau \E\int_{t_{m}}^{t_{m+1}}|(-A_h)^{1/2+\kappa}S_{\tau,h}P_h|_{\mathcal{L}(H)}|P_hF(X_m^h)||(-A_h)^{1/2-\kappa}P_hD\tilde{\Psi}^{(M)}(\tilde{X}^h(t))|dt\\
&\leq \|F\|_\infty\tau^{1/2-\kappa}\tau(1+|x|^2).
\end{align*}
We then have
\begin{equation}\label{e6}
\frac{1}{N\tau}\sum_{m=1}^{N-1}|b_{m}^{2}|\leq C\tau^{1/2-\kappa}(1+|x|^2).
\end{equation}

\item \textbf{Estimate of $b_{m}^{3}$}

We again use the integration by parts formula to rewrite $b_{m}^{3}$:
\begin{align*}
b_{m}^{3}&=\E\int_{t_m}^{t_{m+1}}\langle \int_{t_m}^{t}A_hS_{\tau,h}P_hdW(s),P_hD\tilde{\Psi}^{(M)}(\tilde{X}^h(t))\rangle dt\\
&=\E\int_{t_{m}}^{t_{m+1}}\int_{t_m}^{t}\text{Tr}(S_{\tau,h}A_hP_hD^2\tilde{\Psi}^{(M)}(\tilde{X}^h(t))\DM_s\tilde{X}^h(t))dsdt.
\end{align*}
From \eqref{tild}, for $t_m\leq s\leq t\leq t_{m+1}$ we have $\DM_{s}^{\ell}\tilde{X}^h(t)=S_{\tau,h}P_h\ell$; as a consequence, the situation is much simpler and we do not need to use the same trick as in the control of $a_{m}^{3}$.

Then we have
\begin{align*}
|b_{m}^{3}|&\leq\E\int_{t_m}^{t_{m+1}}(t-t_m)\text{Tr}(S_{\tau,h} A_h P_hD^2\tilde{\Psi}^{(M)}(\tilde{X}^h(t))S_{\tau,h}P_h)dt\\
&\leq C\tau\int_{t_{m}}^{t_{m+1}}|S_{\tau,h}(-A_h)^{1/2+\kappa/2}P_h|_{\mathcal{L}(H)}\text{Tr}((-A)^{-1/2-\kappa/2})|(-A)^\kappa S_{\tau,h}P_h|_{\mathcal{L}(H)}\\
&|(-A_h)^{1/2-\kappa/2}P_hD^2\tilde{\Psi}^{(M)}(\tilde{X}^h(t))(-A)^{1/2-\kappa/2}|_{\mathcal{L}(H)}dt\\
&\leq C(1+|x|^2)\tau^{1/2-3\kappa/2}\tau.
\end{align*}
Therefore
\begin{equation}\label{e7}
\frac{1}{N\tau}\sum_{m=1}^{N-1}|b_{m}^{3}|\leq C(1+|x|^2)\tau^{1/2-3\kappa/2}.
\end{equation}

\end{enumerate}

With the previous estimates, we get
\begin{equation}\label{eb}
\frac{1}{N\tau}\sum_{m=1}^{N-1}|b_{m}|\leq C(1+|x|^3)(1+T^{-1})\tau^{1/2-2\kappa}.
\end{equation}

\subsubsection{Estimate of $c_{m}$}

This term is easy to treat: we have
\begin{align*}
|c_{m}|&\leq \E\int_{t_{m}}^{t_{m+1}}|(-A_h)^{-1/2+\kappa}(I-S_{\tau,h})|_{\mathcal{L}(V_h)}|P_hF(X_m^h)||(-A_h)^{1/2-\kappa}P_hD\tilde{\Psi}^{(M)}(\tilde{X}^h(t))|dt\\
&\leq C\tau^{1/2-\kappa}\tau(1+|x|^2),
\end{align*}
where we have used Proposition \ref{Prop220}, Assumption \ref{hypG} and Lemma \ref{lem6}. 
Then we see that
\begin{equation}\label{e8}
\frac{1}{N\tau}\sum_{m=1}^{N-1}|c_{m}|\leq C\tau^{1/2-\kappa}(1+|x|^2).
\end{equation}

\subsubsection{Estimate of $d_{m}$}

The term $d_m$ contains the error between $F^h(\tilde{X}^h(t))$ and $F^h(X_m^h)$; we recall that $F^h=P_h\circ F$. We perform an expansion with respect to an orthonormal basis $(e_i^h)_{i\in \N}$ of $H$, such that $(e_i^h)_{i=0}^{N_h-1}$ is the orthonormal basis of $V_h$ introduced in Proposition \ref{propA_h} - the vectors $e_i^h$ for $i\geq N_h$ do not matter.

In this orthonormal system, the cylindrical Wiener process is expanded as
\begin{equation}\label{expandW_h}
W(t)=\sum_{i\in\N}\beta_i^h(t)e_i^h,
\end{equation}
with a family $(\beta_i^h)_{i\in\N}$ of independent standard one-dimensional Wiener processes.


Let $F_i^h:H\rightarrow \R$ denote the function such that $F_i^h(x)=\langle F^h(x),e_i^h\rangle $. We also denote, for any $i\in\N$, by $\partial_i^h$ the operator such that $\partial_i^h\phi(x)=\langle D\phi(x),e_i^h\rangle \in\R$, for any $x\in H$, where $\phi:H\rightarrow\R$ is of class $\mathcal{C}^1$. Then
$$\langle F^h(\tilde{X}^h(t))-F^h(X_m^h),D\tilde{\Psi}^{(M)}(\tilde{X}^h(t))\rangle =\sum_{i=0}^{N_h-1}(F_i^h(\tilde{X}^h(t))-F_i^h(X_m))\partial_i^h\tilde{\Psi}^{(M)}(\tilde{X}^h(t)).$$



The It\^o formula gives for $t_m\leq t<t_{m+1}$ and $0\leq i\leq N_h-1$
\begin{align*}
F_i^h(\tilde{X}^h(t))-F_i^h(X_m^h)&=\frac{1}{2}\int_{t_m}^{t}\text{Tr}(S_{\tau,h}S_{\tau,h}^{*}P_hD^2F_i^h(\tilde{X}^h(s)))ds\\
&+\int_{t_m}^{t}\langle A_hS_{\tau,h} X_m^h,DF_i^h(\tilde{X}^h(s))\rangle ds\\
&+\int_{t_m}^{t}\langle S_{\tau,h} F^h(X_m^h),DF_i^h(\tilde{X}^h(s))\rangle ds\\
&+\int_{t_m}^{t}\langle DF_i^h(\tilde{X}^h(s)),S_{\tau,h}P_h dW(s)\rangle .
\end{align*}

Thanks to this, $d_{m}^{j}$ for $j\in\left\{1,2,3,4\right\}$ are naturally defined, and we now control each term.

\begin{enumerate}

\item \textbf{Estimate of $d_{m}^{1}$}

By definition, we have
$$d_{m}^{1}=\int_{t_m}^{t_{m+1}}\E\frac{1}{2}\int_{t_m}^{t}\sum_{i=0}^{N_h-1}\text{Tr}(S_{\tau,h} S_{\tau,h}^{*}P_hD^2F_i^h(\tilde{X}^h(s)))ds\partial_i^h\tilde{\Psi}^{(M)}(\tilde{X}^h(t))dt.$$
Expanding the trace thanks to the complete orthonormal system $(e_{i}^{h})_{i\in\N}$,
\begin{align*}
\sum_{i=0}^{N_h-1}\text{Tr}(P_hS_{\tau,h}^{*}&P_hD^2F_i^h(\tilde{X}^h(s))S_{\tau,h}P_h)\partial_i^h\tilde{\Psi}^{(M)}(\tilde{X}^h(t))\\
&=\sum_{i=0}^{N_h-1}\sum_{j=0}^{N_h-1}\langle D^2F_i^h(\tilde{X}^h(s))\frac{1}{(1+\lambda_j^h\tau)^2}e_j^h,e_j^h\rangle \partial_i^h\tilde{\Psi}^{(M)}(\tilde{X}^h(t))\\
&=\sum_{i=0}^{N_h-1}\sum_{j=0}^{N_h-1}\frac{1}{(1+\lambda_j^h\tau)^2}D^2F_i^h(\tilde{X}^h(s)).(e_j^h,e_j^h)\partial_i^h\tilde{\Psi}^{(M)}(\tilde{X}^h(t)).
\end{align*}

For any fixed $0\leq j\leq N_h-1$, the Cauchy-Schwarz inequality yields
\begin{multline*}
|\sum_{i=0}^{N_h-1}D^2F_i^h(\tilde{X}^h(s)).(e_j^h,e_j^h)\partial_i^h\tilde{\Psi}^{(M)}(\tilde{X}^h(t))|\\
\leq \left(\sum_{i=0}^{N_h-1}\frac{|D^2F_i^h(\tilde{X}^h(s)).(e_j^h,e_j^h)|^2}{(\lambda_{i}^{h})^{2\eta}}\right)^{1/2}\left(\sum_{i=0}^{N_h-1}(\lambda_{i}^{h})^{2\eta}|\partial_i^h\tilde{\Psi}^{(M)}(\tilde{X}^h(t))|^2\right)^{1/2},
\end{multline*}
where $\eta\leq 1/2$ is defined in Assumption \ref{hypG}.

The second factor is bounded from above by $|(-A_h)^{\eta}P_hD\tilde{\Psi}^{(M)}(\tilde{X}^h(t))|\leq C|(-A)^{\eta}D\tilde{\Psi}^{(M)}(\tilde{X}^h(t))|$, thanks to Proposition \ref{conv_FEM}; the right-hand side is then controlled thanks to Proposition \ref{Prop220}.

To control the first factor, thanks to Assumption \ref{hypG} we get
\begin{align*}
\left(\sum_{i=0}^{N_h-1}\frac{|D^2F_i^h(\tilde{X}^h(s)).(e_j^h,e_j^h)|^2}{(\lambda_{i}^{h})^{2\eta}}\right)^{1/2}&=|(-A_h)^{-\eta}P_hD^2F(\tilde{X}^{h}(s)).(e_j^h,e_j^h)|\\
&\leq C|(-A)^{-\eta}D^2F(\tilde{X}^{h}(s)).(e_j^h,e_j^h)|\leq C|e_j^h|_H|e_j^h|_H\leq C,
\end{align*}
since $(e_j^h)_j$ is an orthonormal system, and using Proposition \ref{conv_FEM}.


Therefore we obtain
\begin{align*}
|\sum_{i=0}^{N_h-1}\text{Tr}(S_{\tau,h} S_{\tau,h}^{*}&P_hD^2F_i^h(\tilde{X}^h(s)))\partial_i^hD\tilde{\Psi}^{(M)}(\tilde{X}^h(t))|\\
&\leq C(1+|x|^2)\sum_{j=0}^{N_h-1}\frac{1}{(1+\lambda_j^h\tau)^2}\\
&\leq C(1+|x|^2) \tau^{-1/2-\kappa}\sum_{j=0}^{N_h-1}\frac{(\lambda_j^h\tau)^{1/2+\kappa}}{(1+\lambda_j^h\tau)^2}\frac{1}{(\lambda_{j}^{h})^{1/2+\kappa}}\\
&\leq C(1+|x|^2)\tau^{-1/2-\kappa}\text{Tr}\bigl(P_h(-A_h)^{-1/2-\kappa}P_h).
\end{align*}
Then thanks to Proposition \ref{propoTrace}, we have
$$|d_{m}^{1}|\leq C(1+|x|^2)\tau^{1/2-\kappa}\tau,$$
and
\begin{equation}\label{e9}
\frac{1}{N\tau}\sum_{m=1}^{N-1}|d_{m}^{1}|\leq C(1+|x|^2)\tau^{1/2-\kappa}.
\end{equation}

\item \textbf{Estimate of $d_{m}^{2}$}

Thanks to \eqref{exprYk} and \eqref{exprsto}, we have
\begin{align*}
d_{m}^{2}&=\int_{t_m}^{t_{m+1}}\int_{t_m}^{t}\sum_{i=0}^{N_h-1}\langle A_hS_{\tau,h} X_m^h,DF_i^h(\tilde{X}^h(s))\rangle \partial_i^h\tilde{\Psi}^{(M)}(\tilde{X}^h(t))dsdt\\
&=\E\int_{t_{m}}^{t_{m+1}}\int_{t_m}^{t}\sum_{i}\langle A_h S_{\tau,h}^{m+1}P_hx+\tau\sum_{l=0}^{m-1}A_hS_{\tau,h}^{m-l+1}F^h(X_l^m),DF_i^h(\tilde{X}^h(s))\rangle \partial_i^h\tilde{\Psi}^{(M)}(\tilde{X}^h(t))dsdt\\
&+\E\int_{t_{m}}^{t_{m+1}}\int_{t_m}^{t}\sum_{i}\langle A_h\int_{0}^{t_m}S_{\tau,h}^{m-l_r+1}P_hdW(r),DF_i^h(\tilde{X}^h(s))\rangle \partial_i^h\tilde{\Psi}^{(M)}(\tilde{X}^h(t))dsdt\\
&:=d_{m}^{2,1}+d_{m}^{2,2}.
\end{align*}

{\em (i)} For the first term, since $F$ and $DF$ are bounded on $H$ and $\tau\leq \tau_0$, we have, using \eqref{sumesti},
\begin{align*}
|d_{m}^{2,1}|&=|\E\int_{t_{m}}^{t_{m+1}}\int_{t_m}^{t}\langle DF^h(\tilde{X}^h(s)).(A_h S_{\tau,h}^{m+1}P_hx+A_h\tau\sum_{l=0}^{m-1}S_{\tau,h}^{m-l+1}F^h(X_l^h)),D\Psi^{(M)}(\tilde{X}^h(t))\rangle dsdt|\\
&\leq\E\int_{t_{m}}^{t_{m+1}}\int_{t_m}^{t}|(-A_h)^\kappa S_{\tau,h}P_h|_{\mathcal{L}(H)}\Big(|(-A_h)^{1-\kappa}S_{\tau,h}^{m}P_hx|+\tau\sum_{l=0}^{m-1}|(-A_h)^{1-\kappa}S_{\tau,h}^{m-l}P_h|_{\mathcal{L}(H)}|P_hF(X_l^h)|\Big)\\
&\times|D\Psi^{(M)}(\tilde{X}^h(t))|dsdt\\
&\leq C\tau^{1-\kappa}\tau(1+|x|^2)\Big(t_{m}^{-1+\kappa}|x|+\tau\sum_{l=0}^{m-1}t_{m-l}^{-(1-\kappa)}\frac{1}{(1+\lambda_0\tau)^{(m-l)\kappa}}\Big)\\
&\leq C\tau^{1-\kappa}(1+|x|^3)\tau(\frac{1}{t_{m}^{1-\kappa}}+1).
\end{align*}
Therefore
\begin{align}\label{eq10bis}
\frac{1}{N\tau}\sum_{m=1}^{N-1}|d_{m}^{2,1}|&\leq C\tau^{1-\kappa}(|x|^3+1)\frac{\tau}{T}\sum_{m=1}^{N-1}(\frac{1}{t_{m}^{1-\kappa}}+1)\nonumber\\
&\leq C\tau^{1-\kappa}(1+|x|^3)\big(\frac{1}{T}\int_{0}^{T}\frac{1}{t^{1-\kappa}}dt+1\big)\nonumber\\
&\leq C\tau^{1-\kappa}(1+|x|^3)\big(1+T^{-1+\kappa}\int_{0}^{1}\frac{1}{s^{1-\kappa}}ds\big)\nonumber\\
&\leq C\tau^{1-\kappa}(1+|x|^3)\big(1+T^{-1+\kappa}\big).
\end{align}

{\em (ii)} For the second term, we again use an integration by parts, after a decomposition of the time interval - as in the estimates for $a_{m}^{3}$.
First,
\begin{align*}
&d_{m}^{2,2}=\E\int_{t_{m}}^{t_{m+1}}\int_{t_m}^{t}\sum_{i=0}^{N_h-1}\langle A_h\int_{0}^{t_m}S_{\tau,h}^{m-l_r+1}P_hdW(r),DF_i^h(\tilde{X}^h(s))\rangle \partial_i^h\tilde{\Psi}^{(M)}(\tilde{X}^h(t))dsdt\\
&=\E\int_{t_{m}}^{t_{m+1}}\int_{t_m}^{t}\sum_{i=0}^{N_h-1}\langle A_h\int_{0}^{(t_m-5\tau_0)\vee 0}S_{\tau,h}^{m-l_r+1}P_hdW(r),DF_i^h(\tilde{X}^h(s))\rangle \partial_i^h\tilde{\Psi}^{(M)}(\tilde{X}^h(t))dsdt\\
&+\E\int_{t_m}^{t_{m+1}}\int_{t_m}^{t}\sum_{i,j,m=0}^{N_h-1}\langle A_h\int_{(t_m-5\tau_0)\vee 0}^{t_m}S_{\tau,h}^{m-l_r+1}e_m^h,e_j^h\rangle d\beta_m^h(r)\partial_j^hF_i^h(\tilde{Y}(s))\partial_i^h\tilde{\Psi}^{(M)}(\tilde{X}^h(t))dsdt\\
&=:d_{m}^{2,2,1}+d_{m}^{2,2,2}.
\end{align*}
Recall that $(e_i^h)_{0\leq i\leq N_h-1}$ is the orthonormal basis of $V_h$ introduced in Proposition \ref{propA_h}. See also \eqref{expandW_h}.
For $d_{m}^{2,2,1}$, we can work directly as for the similar part in $a_m^3$ and we see that
\begin{align*}
|d_{m}^{2,2,1}|&\leq |\E\int_{t_{m}}^{t_{m+1}}\int_{t_m}^{t}\sum_{i=0}^{N_h-1}\langle A_h\int_{0}^{(t_m-5\tau_0)\vee 0}S_{\tau,h}^{m-l_r+1}P_hdW(r),DF_i^h(\tilde{X}^h(s))\rangle \partial_i^h\tilde{\Psi}^{(M)}(\tilde{X}^h(t))dsdt|\\
&\leq \int_{t_m}^{t_{m+1}}\int_{t_m}^{t}\E|\langle DF^h(\tilde{X}^h(s)).A_hS_{\tau,h}\int_{0}^{(t_m-5\tau_0)\vee 0}S_{\tau,h}^{m-l_r}P_hdW(r),D\tilde{\Psi}^{(M)}(\tilde{X}^h(t))\rangle |dsdt\\
&\leq \int_{t_m}^{t_{m+1}}\int_{t_m}^{t}(\E|A_hS_{\tau,h}\int_{0}^{(t_m-5\tau_0)\vee 0}S_{\tau,h}^{m-l_r}P_hdW(r)|^{2})^{1/2}(\E|D\tilde{\Psi}^{(M)}(\tilde{X}^h(t))|^2)^{1/2}dsdt\\
&\leq C\tau^2(1+|x|^2),
\end{align*}
thanks to Lemmas \ref{lem2}, Proposition \ref{Prop220} and to the estimate proved in the control of $a_m^3$:
$$\E|A_hS_{\tau,h}\int_{0}^{(t_m-5\tau_0)\vee 0}S_{\tau,h}^{m-l_r}P_hdW(r)|^{2}\leq\E|A_h^2S_{\tau,h}\int_{0}^{(t_m-5\tau_0)\vee 0}S_{\tau,h}^{m-l_r}P_hdW(r)|^{2}\leq C.$$

For $d_{m}^{2,2,2}$, we can write thanks to the Malliavin integration by parts \eqref{formulaintbyparts} and with the chain rule
\begin{align*}
&d_{m}^{2,2,2}=\E\int_{t_m}^{t_{m+1}}\int_{t_m}^{t}\sum_{i,j,m=0}^{N_h-1}\langle A_h\int_{(t_m-5\tau_0)\vee 0}^{t_m}S_{\tau,h}^{m-l_r+1}e_m^h,e_j^h\rangle d\beta_m^h(r)\partial_j^hF_i^h(\tilde{Y}(s))\partial_i^h\tilde{\Psi}^{(M)}(\tilde{X}^h(t))dsdt\\
&=\E\int_{t_{m}}^{t_{m+1}}\int_{t_{m}}^{t}\int_{(t_m-5\tau_0)\vee 0}^{t_m}\sum_{i,j,m,n=0}^{N_h-1}\langle A_h S_{\tau,h}^{m-l_r+1}e_m^h,e_j^h\rangle \partial_{j}^{h}\partial_{n}^{h}F_i^h(\tilde{X}^h(s))\langle \DM_{r}^{e_m^h}\tilde{X}^h(s),e_n^h\rangle \partial_i^h\tilde{\Psi}^{(M)}(\tilde{X}^h(t))drdsdt\\
&+\E\int_{t_{m}}^{t_{m+1}}\int_{t_{m}}^{t}\int_{(t_m-5\tau_0)\vee 0}^{t_m}\sum_{i,j,m,n=0}^{N_h-1}\langle A_h S_{\tau,h}^{m-l_r+1}e_m^h,e_j^h\rangle \partial_j^hF_i^h(\tilde{X}^h(s))\partial_{n}^{h}\partial_{i}^{h}\tilde{\Psi}^{(M)}(\tilde{X}^h(t))\langle \DM_{r}^{e_m^h}\tilde{X}^h(t),e_n^h\rangle drdsdt\\
&=\E\int_{t_{m}}^{t_{m+1}}\int_{t_{m}}^{t}\int_{(t_m-5\tau_0)\vee 0}^{t_m}\sum_{i,m=0}^{N_h-1}D^2F_i^h(\tilde{X}^h(s))(A_h S_{\tau,h}^{m-l_r+1}e_m^h,\DM_{r}^{e_m^h}\tilde{X}^h(s))\partial_i^h\tilde{\Psi}^{(M)}(\tilde{X}^h(t))drdsdt\\
&+\E\int_{t_{m}}^{t_{m+1}}\int_{t_{m}}^{t}\int_{(t_m-5\tau_0)\vee 0}^{t_m}\sum_{i,m=0}^{N_h-1}\langle \mathcal{B}_i(s,t)A_h S_{\tau,h}^{m-l_r+1}e_m^h,\DM_{r}^{e_m^h}\tilde{X}^h(t)\rangle drdsdt\\
&=\E\int_{t_{m}}^{t_{m+1}}\int_{t_{m}}^{t}\int_{(t_m-5\tau_0)\vee 0}^{t_m}\sum_{i=0}^{N_h-1}\text{Tr}\left(P_h(\DM_{r}\tilde{X}^h(s))^{*}D^2F_i^h(\tilde{X}^h(s))A_h S_{\tau,h}^{m-l_r+1}P_h\right)\partial_i^h\tilde{\Psi}^{(M)}(\tilde{X}^h(t))drdsdt\\
&+\E\int_{t_{m}}^{t_{m+1}}\int_{t_{m}}^{t}\int_{(t_m-5\tau_0)\vee 0}^{t_m}\sum_{i=0}^{N_h-1}\text{Tr}\left(P_h(\DM_r\tilde{X}^h(t)^{*}\mathcal{B}_i(s,t)A_h S_{\tau,h}^{m-l_r+1}P_h\right)drdsdt,
\end{align*}
where we define, for each $0\leq i\leq N_h-1$, a linear operator on $V_h$ by
\begin{align*}
\langle \mathcal{B}_i(s,t)x_h,y_h\rangle &=\langle DF_i^h(\tilde{X}^h(s)),x_h\rangle \sum_{n=0}^{N_h-1}\partial_{n}^{h}\partial_{i}^{h}\tilde{\Psi}^{(M)}(\tilde{X}^h(t))\langle y_h,e_n^h\rangle \\
&=\langle DF_i^h(\tilde{X}^h(s)),x_h\rangle \langle D^2\tilde{\Psi}^{(M)}(\tilde{X}^h(t)).e_i^h,y_h\rangle .
\end{align*}

We have, for any $x_h,y_h\in V_h$, $\sum_{i=0}^{N_h-1}\langle \mathcal{B}_i(s,t)x_h,y_h\rangle =D^2\tilde{\Psi}^{(M)}(\tilde{X}^h(t)).(DF^h(\tilde{X}^h(s)).x_h,y_h)$
and, using Proposition \ref{Prop220} and Assumption \ref{hypG},
$$|\sum_{i=0}^{N_h-1}\mathcal{B}_i(s,t)|_{\mathcal{L}(V_h)}\leq|DF(\tilde{X}^h(s))|_{\mathcal{L}(H)}|D^2\tilde{\Psi}^{(M)}(\tilde{X}^h(t))|_{\mathcal{L}(H)};$$
so we can write, for $(t_m-5\tau_0)\vee 0\leq r\leq t_m$
\begin{align*}
&|\sum_{i=0}^{N_h-1}\text{Tr}\left(P_h(\DM_r\tilde{X}^h(t)^{*}\mathcal{B}_i(s,t)A_h S_{\tau,h}^{m-l_r+1}P_h\right)|\\
&\leq|\DM_r\tilde{X}^h(t)|_{\mathcal{L}(V_h)}|\sum_{i=0}^{N_h-1}\mathcal{B}_i(s,t)|_{\mathcal{L}(V_h)}|(-A_h)^{1-3\kappa/2}S_{\tau,h}^{m-l_r}P_h|_{\mathcal{L}(H)}|S_{\tau,h}(-A_h)^{1/2+2\kappa}P_h|_{\mathcal{L}(H)}\text{Tr}(P_h(-A_h)^{-1/2-\kappa/2}P_h)\\
&\leq C\tau^{-1/2-2\kappa}t_{m-l_r}^{-1+3\kappa/2}\frac{1}{(1+\lambda_0\tau)^{(m-l_r)3\kappa/2}},
\end{align*}
using Lemma \ref{lem5} - since $(1+L_F\tau)^{m-l_r}\leq C$ - and Lemma \ref{lem6}.

For the other term, we have to deal with the poor regularity of $F$ at second order. We proceed as in the control of $d_m^1$, and expand the trace with respect to the orthonormal basis $(e_i^h)_{i\in\N}$ of $H$.

\begin{align*}
&|\sum_{i=0}^{N_h-1}\text{Tr}\left(P_h(\DM_{r}\tilde{X}^h(s))^{*}D^2F_i^h(\tilde{X}^h(s))A_h S_{\tau,h}^{m-l_r+1}P_h\right)\partial_i^h\tilde{\Psi}^{(M)}(\tilde{X}^h(t))|\\
&\leq |\DM_r\tilde{X}^h(s)|_{\mathcal{L}(V_h)}\sum_{i,j=0}^{N_h-1}\frac{|D^2F_i^h(\tilde{X}^h(s)).(e_j^h,e_j^h)|}{(\lambda_{i}^{h})^{\eta}}\frac{\lambda_j^h}{(1+\lambda_j^h\tau)^{1+m-l_r}}(\lambda_{i}^{h})^{\eta}|\partial_i^h\tilde{\Psi}^{(M)}(\tilde{X}^h(t))|\\
&\leq |\DM_r\tilde{X}^h(s)|_{\mathcal{L}(V_h)}|(-A_h)^\eta P_hD\tilde{\Psi}^{(M)}(\tilde{X}^h(t))|_{H}\sum_{j=0}^{N_h-1}|(-A_h)^{-\eta}D^2F^h(\tilde{X}^h(s)).(e_j^h,e_j^h)|\frac{\lambda_j^h}{(1+\lambda_j^h\tau)^{1+m-l_r}},
\end{align*}
thanks to the Cauchy-Schwarz inequality.

By using the same analysis as in the estimation of $d_{m}^{1}$, we see that the above expression is bounded by
$$C|\DM_r\tilde{X}^h(s)|_{\mathcal{L}(V_h)}|(-A_h)^\eta P_hD\tilde{\Psi}^{(M)}(\tilde{X}^h(t))|_{H}\sum_{j=0}^{N_h-1}\frac{\lambda_j^h}{(1+\lambda_j^h\tau)^{1+m-l_r}};$$
but the last sum is equal to $\text{Tr}(P_hA_h S_{\tau,h}^{m-l_r+1})$, so that we see that indeed the two expressions in $d_{m}^{2,2}$ are bounded by the same expression.

Therefore
\begin{align*}
&|d_{m}^{2,2,2}|\\
&\leq \E\int_{t_{m}}^{t_{m+1}}\int_{t_{m}}^{t}\int_{(t_m-5\tau_0)\vee 0}^{t_m}C\tau^{-1/2-2\kappa}t_{m-l_r}^{-1+3\kappa/2}\frac{1}{(1+\lambda_0\tau)^{(m-l_r)3\kappa/2}}(1+|x|^2)drdsdt\\
&\leq C(1+|x|^2)\tau^{1/2-2\kappa}\tau\int_{0}^{t_m}t_{m-l_r}^{-1+3\kappa/2}\frac{1}{(1+\lambda_0\tau)^{(m-l_r)3\kappa/2}}dr\\
&\leq C(1+|x|^2)\tau^{1/2-2\kappa}\tau,
\end{align*}
as already proved - see \eqref{sumesti}.

Now gathering estimates for $d_{m}^{2,2,1}$ and $d_{m}^{2,2,2}$, we obtain
\begin{equation}\label{e10}
\frac{1}{N\tau}\sum_{m=1}^{N-1}|d_{m}^{2,2}|\leq C(1+|x|^2)\tau^{1/2-2\kappa}.
\end{equation}

\item \textbf{Estimate of $d_{m}^{3}$}
We have
\begin{align*}
d_{m}^{3}&=\E\int_{t_{m}}^{t_{m+1}}\int_{t_m}^{t}\sum_{i=0}^{N_h-1}\langle S_{\tau,h} F^h(X_m^h),DF_i^h(\tilde{X}^h(s))\rangle \partial_i^h\tilde{\Psi}^{(M)}(\tilde{X}^h(t))dsdt\\
&=\E\int_{t_{m}}^{t_{m+1}}\int_{t_m}^{t}\langle P_hD\tilde{\Psi}^{(M)}(\tilde{X}^h(t)),DF^h(\tilde{X}^h(s)).(S_{\tau,h} F^h(X_m^h))\rangle dsdt.
\end{align*}
We have assumed $F$ and $DF$ to be bounded, so we easily get
$$|d_{m}^{3}|\leq C(1+|x|^2)\tau^2$$
and
\begin{equation}\label{e11}
\frac{1}{N\tau}\sum_{m=1}^{N-1}|d_{m}^{3}|\leq C(1+|x|^2)\tau.
\end{equation}

\item \textbf{Estimate of $d_{m}^{4}$}

Finally, thanks to the integration by parts formula of Proposition \ref{lemintbyparts}, we have
\begin{align*}
d_{m}^{4}&=\E\int_{t_{m}}^{t_{m+1}}\int_{t_m}^{t}\sum_{i=0}^{N_h-1}\langle DF_i^h(\tilde{X}^h(s)),S_{\tau,h}P_h dW(s)\rangle \partial_i^h\tilde{\Psi}^{(M)}(\tilde{X}^h(t))dt\\
&=\E\int_{t_m}^{t_{m+1}}\int_{t_m}^{t}\text{Tr}\left((\DM_s\tilde{X}^h(t))^{*}P_hD^2\tilde{\Psi}^{(M)}(\tilde{X}^h(t))DF^h(\tilde{X}^h(s))S_{\tau,h}P_h\right)dsdt\\
&=\E\int_{t_m}^{t_{m+1}}\int_{t_m}^{t}\text{Tr}\left(S_{\tau,h}P_hD^2\tilde{\Psi}^{(M)}(\tilde{X}^h(t))DF^h(\tilde{X}^h(s))S_{\tau,h}\right)dsdt;
\end{align*}
indeed, we have $\DM_{s}^{\ell}\tilde{Y}(t)=S_{\tau,h}\ell$ for all $\ell\in V_h$, when $t_m\leq s\leq t\leq t_{m+1}$ - see also the control of $b_m^3$. Now,
\begin{align*}
|d_{m}^{4}|&\leq\E\int_{t_m}^{t_{m+1}}\int_{t_m}^{t}|(S_{\tau,h}(-A_h)^{1/2+\kappa}|_{\mathcal{L}(V_h)}|DF(\tilde{X}^h(s))|_{\mathcal{L}(H)}|S_{\tau,h}|_{\mathcal{L}(V_h)}\\
&\times |D^2\tilde{\Psi}^{(M)}(\tilde{X}^h(t))|_{\mathcal{L}(H)}\text{Tr}(P_h(-A_h)^{-1/2-\kappa}P_h)dsdt\\
&\leq C(1+|x|^2)\tau^{1/2-\kappa}\tau,
\end{align*}
and
\begin{equation}\label{e12}
\frac{1}{N\tau}\sum_{m=1}^{N-1}|d_{m}^{4}|\leq  C(1+|x|^2)\tau^{1/2-\kappa}.
\end{equation}

\item \textbf{Estimate of $d_m$: conclusion}
With \eqref{e9}, \eqref{eq10bis}, \eqref{e10}, \eqref{e11} and \eqref{e12}, we get
\begin{equation}\label{eb_2}
\frac{1}{N\tau}\sum_{m=1}^{N-1}|d_m|\leq C\tau^{1/2-2\kappa}(1+|x|^3)(1+T^{-1+\kappa}).
\end{equation}

\end{enumerate}

\subsubsection{Estimate of $e_m$}

Using the symmetry of the operators $P_h\in\mathcal{L}(H)$ and $S_{\tau,h}\in\mathcal{L}(V_h)$, the commutativity of $P_h$ and $S_{\tau,h}$ and the fact that $P_h$ is a projector, we have in $\mathcal{L}(H)$
$$P_hP_{h}^{*}-(S_{\tau,h}P_h)(S_{\tau,h}P_h)^{*}=2P_hS_{\tau,h}P_h(I-S_{\tau,h})P_h+P_h(I-S_{\tau,h})P_h(I-S_{\tau,h})P_h;$$
we then decompose $e_m$ into two parts: $e_m=e_{m}^{1}+e_{m}^{2}$, with
\begin{gather*}
e_{m}^{1}:=\E\int_{t_m}^{t_{m+1}}\text{Tr}\bigl(P_hS_{\tau,h}P_h(I-S_{\tau,h})P_hD^2\tilde{\Psi}^{(M)}(\tilde{X}^h(t))\bigr);\\
e_{m}^{2}:=\frac{1}{2}\E\int_{t_m}^{t_{m+1}}\text{Tr}\bigl(P_h(I-S_{\tau,h})P_h(I-S_{\tau,h})P_hD^2\tilde{\Psi}^{(M)}(\tilde{X}^h(t))\bigr).
\end{gather*}

\begin{enumerate}

\item \textbf{Estimate of $e_{m}^{1}$}

\begin{align*}
|e_{m}^{1}|
&\leq \E\int_{t_m}^{t_{m+1}}\big|\text{Tr}\bigl((-A)^{-1/2-\kappa}(-A)^{2\kappa}P_hS_{\tau,h}P_h(I-S_{\tau,h})P_hD^2\tilde{\Psi}^{(M)}(\tilde{X}^h(t))(-A)^{1/2-\kappa}\bigr)\big|ds\\
&\leq \E\int_{t_m}^{t_{m+1}}\text{Tr}((-A)^{-1/2-\kappa})\big|(-A)^{2\kappa}P_hS_{\tau,h}P_h(I-S_{\tau,h})(-A_h)^{-1/2+\kappa}\big|_{\mathcal{L}(H)}\\
&\hspace{2cm}\times\big|(-A_h)^{1/2-\kappa}P_hD^2\tilde{\Psi}^{(M)}(\tilde{X}^h(t))(-A)^{1/2-\kappa}\big|_{\mathcal{L}(H)}ds\\
&\leq C(1+|x|^2)\tau \big|(-A_h)^{2\kappa}S_{\tau,h}P_h(I-S_{\tau,h})(-A_h)^{-1/2+\kappa}\big|_{\mathcal{L}(V_h)}\\
&\leq C(1+|x|^2)\tau \tau^{-2\kappa}\tau^{1/2-\kappa}=C\tau^{1+1/2-3\kappa}.
\end{align*}

\item \textbf{Estimate of $e_{m}^{2}$}

\begin{align*}
|e_{m}^{2}|&=\big|\frac{1}{2}\E\int_{t_m}^{t_{m+1}}\text{Tr}\bigl(P_h(I-S_{\tau,h})P_h(I-S_{\tau,h})P_hD^2\tilde{\Psi}^{(M)}(\tilde{X}^h(t))\bigr)\big|\\
&\leq \frac{1}{2}\E\int_{t_m}^{t_{m+1}}\big|\text{Tr}\bigl((-A)^{-1/2-\kappa}(-A)^{2\kappa}P_h(I-S_{\tau,h})P_h(I-S_{\tau,h})\\
&\hspace{2cm} \times P_h(-A_h)^{-1/2+\kappa}(-A_h)^{1/2-\kappa}D^2\tilde{\Psi}^{(M)}(\tilde{X}^h(t))(-A)^{1/2-\kappa}\bigr) \big|ds\\
&\leq C(1+|x|^2)\tau\text{Tr}((-A)^{-1/2-\kappa})\big|(-A)^{2\kappa}P_h(I-S_{\tau,h})P_h(I-S_{\tau,h})P_h(-A_h)^{-1/2+\kappa}\big|_{\mathcal{L}(H)}\\
&\leq C(1+|x|^2)\tau\big|(-A_h)^{2\kappa}(I-S_{\tau,h})P_h(I-S_{\tau,h})(-A_h)^{-1/2+\kappa}\big|_{\mathcal{L}(H)}\\
&\leq C(1+|x|^2)\tau \tau^{-2\kappa}\tau^{1/2-\kappa}=C\tau^{1+1/2-3\kappa}.
\end{align*}

\item \textbf{Estimate of $e_m$: conclusion}

We thus have
\begin{equation}\label{ec}
\frac{1}{\tau N}\sum_{m=1}^{N-1}|e_m|\leq C(1+|x|^2)\tau^{1/2-3\kappa}.
\end{equation}

\end{enumerate}

\subsubsection{Conclusion} With the above estimations, we get
\begin{equation}\label{ef1}
|\frac{1}{N\tau}\sum_{m=1}^{N-1}\int_{t_m}^{t_{m+1}}\E\Big(\mathcal{L}^h-\mathcal{L}^{\tau,m,h}\Big)\PsiM(\tilde{X}^h(t))dt|\leq C(1+|x|^3)\tau^{1/2-\kappa}(1+T^{-1+\kappa}+T^{-1}) .
\end{equation}

\subsection{Conclusion} With \eqref{ef1}, \eqref{ef3} and \eqref{ef2}, we get
\begin{equation*}
\frac{1}{N}\sum_{m=0}^{N-1}\Bigl(\phi(X_m^{h})-\overline{\phi}\Bigr)\leq C(1+|x|^3)\tau^{1/2-\kappa}(1+T^{-1+\kappa}+T^{-1})(1+h^{1-\kappa}),
\end{equation*}
where $C$ does not depend of $T$, $h$ and $M$.

\newpage
\begin{appendix}
\section{Study of the (finite-dimensional) Poisson equation}

\vspace{0.5cm}

This Section is devoted to the proof of Proposition \ref{Prop220}.

\vspace{0.5cm}

Let $\phi\in \mathcal{C}^2_b$. For lighter notation, we will assume in this appendix that $\overline{\phi}=0$. We define the function $u$ for any $t>0$ and $\in H$ by
\begin{equation}\label{defu}
u(t,x)=\E[\phi(X(t,x))],
\end{equation}
which is solution of a finite dimensional Kolmogorov equation associated with the Galerkin finite dimensional approximation of \eqref{eqY}:
$$\frac{\text{d}u}{\text{d}t}(t,x)=Lu(t,x)=\frac{1}{2}\text{Tr}\left(D^2u(t,x)\right)+\langle Ax+F(x),Du(t,x)\rangle .$$
Since $\phi$ is of class $\mathcal{C}^2$, bounded and with bounded derivatives, we are able to prove that with respect to $y$ the function $u$ is twice differentiable. Then, using the It\^o formula, we can show that $\Psi$ is solution of \eqref{eqPoiss}. 

To prove Proposition \ref{Prop220}, we only need to show that $u\in\mathcal{C}^2$ and that $u$ and its two first derivates have estimates which are integrable with respect to $t$. In fact, we will show the result below:

\begin{propo}\label{lem3}
Let $\phi\in\mathcal{C}^2_b$ such that $\overline{\phi}=0$ and $u$ defined by \eqref{defu}. There exist constants $C$, $c$ and $\tilde{\mu}>0$ such that for any $0\leq\beta,\gamma<1/2$ there exist constants $C_{\beta}$ and $C_{\beta,\gamma}$ such that for any $t>0$ and $y\in H$
\begin{equation}\label{decreaseU}
|u(t,x)|\leq C(1+|x|^2)\parallel \phi \parallel_{\infty}
\end{equation}
\begin{equation}\label{decreaseDU2}
|Du(t,x)|_{\beta}\leq C_\beta(1+\frac{1}{t^\beta})e^{-\tilde{\mu}t}(1+|x|^2)\parallel \phi \parallel_{1,\infty}.
\end{equation}
and
\begin{equation}\label{decreaseDU3}
|(-A)^\beta D^2 u(t,x)(-A)^\gamma|_{\mathcal{L}(H)}\leq C_{\beta,\gamma} (1+\frac{1}{t^\eta}+\frac{1}{t^{\beta+\gamma}})e^{-\tilde{\mu}t}(1+|x|^2)\parallel \phi \parallel_{2,\infty},
\end{equation}
where $\eta<1$ is defined in the Assumption \ref{hypG}
\end{propo}

\begin{rem}
In fact the estimation \eqref{decreaseDU2}, is true for $\beta<1$.
\end{rem}

The proof of this result is similar to the proof done in \cite{B2}.

\begin{rem}
 Since $\phi$ is of class $\mathcal{C}^2$, bounded and with bounded derivatives, we are able to prove that with respect to $y$ the function $u$ is twice differentiable, and that the derivatives can be calculated in the following way:
\begin{itemize}
\item For any $h\in H$, we have
\begin{equation}\label{formDu}
Du(t,x).h=\E[D\phi(X(t,x)).\eta^{h,x}(t)],
\end{equation}
where $\eta^{h,x}$ is the solution of
\begin{gather*}
\frac{d\eta^{h,x}(t)}{dt}=A\eta^{h,x}(t)+DF(X(t,x)).\eta^{h,x}(t),\\
\eta^{h,y}(0)=h.
\end{gather*}
\item For any $h,k\in H$, we have
\begin{equation}\label{formD^2u}
D^2u(t,x).(h,k)=\E[D^2\phi(X(t,x)).(\eta^{h,x}(t),\eta^{k,x}(t))+D\phi(X(t,x)).\zeta^{h,k,x}(t)],
\end{equation}
where $\zeta^{h,k,x}$ is the solution of
\begin{gather*}
\frac{d\zeta^{h,k,x}(t)}{dt}=A\zeta^{h,k,x}(t)+DF(X(t,x)).\zeta^{h,k,x}(t)+D^2F(X(t,x)).(\eta^{h,x}(t),\eta^{k,x}(t)),\\
\zeta^{h,k,x}(0)=0.
\end{gather*}
\end{itemize}

Morover, we already have the equation \eqref{decreaseU} (see \eqref{cvexpy1y2}).
\end{rem}

We will now show the equations \eqref{decreaseDU2} and \eqref{decreaseDU3}.
The singularity $t^{-\eta}$ in \eqref{decreaseDU3} is a consequence of the regularity properties satisfied by $F$. 
The proofs require several steps. First in Lemma \ref{lem000} below we prove estimates for a finite horizon and general $0\leq\beta,\gamma<1/2$; then in Lemma \ref{lem00} we study the long-time behaviour in the particular case $\beta=\gamma=0$; we finally conclude with the proofs of Proposition \ref{lem3}.

First, we prove estimates of these quantities for $0<t\leq 1$ - see Lemmas $4.4$ and $4.5$ in \cite{deb}, with a difference coming from the assumptions made on the nonlinear coefficient $F$:
\begin{lemme}\label{lem000}
For any $0\leq \beta<1/2$, $0\leq \gamma<1/2$, there exist constants $C_{\beta}$ and $C_{\beta,\gamma}$ such that for any $y\in H$ and any $0<t\leq 1$
\begin{gather*}
|Du(t,x)|_{\beta}\leq \frac{C_\beta}{t^\beta}\parallel D\phi \parallel_{\infty}\\
|(-A)^\beta D^2 u(t,x)(-A)^\gamma|_{\mathcal{L}(H)}\leq C_{\beta,\gamma}(\frac{1}{t^\eta}+\frac{1}{t^{\beta+\gamma}})(\parallel D\phi\parallel_{\infty}+\parallel D^2\phi \parallel_{\infty}).
\end{gather*}
\end{lemme}
\begin{rem}
If we take another time interval $]0,T_{max}]$ instead of $]0,1]$, the constants $C_{\beta}$ and $C_{\beta,\gamma}$ are \textit{a priori} exponentially increasing in $T_{max}$.
\end{rem}

\underline{Proof}
Owing to \eqref{formDu} and \eqref{formD^2u}, we only need to prove the following almost sure estimates, for some constants $C_\beta$ and $C_{\beta,\gamma}$ - which may vary from line to line below: for any $0<t\leq 1$
\begin{equation}\label{deriv_imp}
\begin{gathered}
|\eta^{h,x}(t)|\leq \frac{C_\beta}{t^\beta}|h|_{-\beta}\\
|\zeta^{h,k,x}(t)|\leq C_{\beta,\gamma}\min(\frac{1}{t^\eta},\frac{1}{t^{\beta+\gamma}})|h|_{-\beta}|k|_{-\gamma},
\end{gathered}
\end{equation}
where the parameter $\eta$ is defined in Assumption \ref{hypG}.

We use mild formulations, and the regularization properties of the semi-group given in Proposition \ref{proporegul}:
\begin{align*}
|\eta^{h,x}(t)|&=|e^{tA}h+\int_{0}^{t}e^{(t-s)A}DF(X(s,y)).\eta^{h,x}(s)ds|\\
&\leq \frac{C_\beta}{t^\beta}|h|_{-\beta}+C\int_{0}^{t}|\eta^{h,x}(s)|ds,
\end{align*}
and by the Gronwall Lemma we get the result.

For the second-order derivative, we moreover use the properties of $F$ in Assumption \ref{hypG} to get
\begin{align*}
|\zeta^{h,k,x}(t)|&=|\int_{0}^{t}e^{(t-s)A}DF(X(s,x)).\zeta^{h,k,x}(s)ds\\
&+\int_{0}^{t}e^{(t-s)A}D^2F(X(s,x)).(\eta^{h,x}(s),\eta^{k,x}(s))ds|\\
&\leq C\int_{0}^{t}|\zeta^{h,k,x}(s)|ds+\int_{0}^{t}\frac{C_{\beta,\gamma}}{(t-s)^\eta}|\eta^{h,y}(s)||\eta^{k,x}(s)|ds\\
&\leq C\int_{0}^{t}|\zeta^{h,k,x}(s)|ds+C_{\beta,\gamma}|h|_{-\beta}|k|_{-\gamma}t^{1-\eta-\beta-\gamma}\int_{0}^{1}\frac{1}{(1-s)^\eta s^{\beta+\gamma}}ds.
\end{align*}
To conclude, it remains to use the Gronwall Lemma, since for any $0<t\leq 1$ we get $t^{1-\eta-\beta-\gamma}\leq t^{-\eta}$, thanks to the assumption $\beta+\gamma<1$.
\qed

Thanks to the dissipativity property expressed in Proposition \ref{propodiss}, we can prove the result in the case $\beta=\gamma=0$. We notice that the proof would be straightforward under a strict dissipativity assumption - since then $\eta^{h,x}(t)$ and $\zeta^{h,k,x}(t)$ would decrease exponentially in $t$; in the general case $\eta^{h,x}(t)$ and $\zeta^{h,k,x}(t)$ are exponentially increasing in time so that we can not work directly. Here the result comes from the estimate \eqref{cvexpy1y2} of Proposition \ref{propoexpy1y2}.
\begin{lemme}\label{lem00}
There exist constants $C$ and $c>0$ such that for any $t\geq 0$ and any $y\in H$
\begin{equation}\label{decreaseDU1}
|Du(t,x)|\leq Ce^{-ct}(1+|x|^2)\parallel \phi \parallel_{\infty} \quad \text{ and } \quad |D^2u(t,x)|_{\mathcal{L}(H)}\leq Ce^{-ct}(1+\frac{1}{t^\eta})(1+|x|^2)\parallel \phi \parallel_{\infty}.
\end{equation}
\end{lemme}

\underline{Proof}
The Bismut-Elworthy-Li formula states that if $\Phi:H\rightarrow \R$ is a function of class $\mathcal{C}^2$ with bounded derivatives and with at most quadratic growth - i.e. there exists $M(\Phi)>0$ such that for any $y\in H$ we have $|\Phi(y)|\leq M(\Phi)(1+|x|^2)$ - then we can calculate the first and the second order derivatives of $(t,x)\mapsto v(t,x):=\E\Phi(X(t,x))$ with respect to $y$.
First, we have for any $y\in H$ and $h\in H$
\begin{equation}\label{BEL1}
\begin{aligned}
Dv(t,x).h&=\frac{1}{t}\E[\int_{0}^{t}\langle \eta^{h,x}(s),dW(s)\rangle \Phi(X(t,x))]\\
&=\frac{2}{t}\E[\int_{0}^{t/2}\langle \eta^{h,x}(s),dW(s)\rangle v(t/2,X(t/2,x))];
\end{aligned}
\end{equation}
the second equality is a consequence of the identity $v(t,x)=\E v(t/2,X(t/2,y))$ obtained with the Markov property, and of the first equality applied with the function $v(t/2,.)$.

Using the second formula of \eqref{BEL1}, we obtain a formula for the second order derivative: for any $y\in H$ and $h,k\in H$,
\begin{equation}\label{BEL2}
\begin{aligned}
D^2v(t,x).(h,k)&=\frac{2}{t}\E[\int_{0}^{t/2}\langle \zeta^{h,k,x}(s),dW(s)\rangle v(t/2,X(T/2,x))]\\
&+\frac{2}{t}\E[\int_{0}^{t/2}\langle \eta^{h,x}(s),dW(s)\rangle Dv(t/2,X(t/2)).\eta^{k,x}(t/2)].
\end{aligned}
\end{equation}
We then see, using Lemmas \ref{lem1} and \ref{lem000} - with $\beta=\gamma=0$ - that there exists $C>0$ such that for any $0\langle t\leq 1$, $x\in H$, $h,k\in H$
\begin{equation}\label{consBEL}
\begin{gathered}
|Dv(t,x).h|\leq \frac{C}{\sqrt{t}}M(\Phi)(1+|x|^2)|h|,\\
|D^2v(t,x).(h,k)|\leq \frac{C}{t}M(\Phi)(1+|x|^2)|h||k|.
\end{gathered}
\end{equation}
Now when $t\geq 1$ the Markov property implies that $u(t,x)=\E u(t-1,X(1,x))$, and by \eqref{cvexpy_int} we have
$$|u(t-1,x)-\int_H\phi d\overline{\mu}|\leq Ce^{-c(t-1)}(1+|x|^2)\parallel \phi\parallel_{\infty}.$$
If we choose $\Phi_t(x)=u(t-1,x)-\int_H\phi d\overline{\mu}$, we have $u(t,x)=\E\Phi_t(X(1,x))+\int_H\phi d\overline{\mu}$, with $M(\Phi_t)\leq Ce^{-c(t-1)}\parallel \phi \parallel_{\infty}$. With \eqref{consBEL} at time $1$, we obtain for $t\geq 1$
\begin{gather*}
|Du(t,x).h|\leq C\parallel \phi \parallel_{\infty}e^{-c(t-1)}(1+|x|^2)|h|\\
|D^2u(t,x).(h,k)|\leq C\parallel \phi \parallel_{\infty}e^{-c(t-1)}(1+|x|^2)|h||k|.
\end{gather*}
Moreover by Lemma \ref{lem000} we have a control when $0\leq t\leq 1$, so that with a change of constants we get the result.
\qed

We can finally prove the Proposition \ref{lem3}. The key tool is the Markov property of the process $X$ which yields the following formula: for any $t\geq 1$
\begin{equation}\label{umarkov}
u(t,x)=\E[u(t-1,X_1(x))].
\end{equation}

To get the exponential decreasing, we use Lemma \ref{lem00} at time $t-1$ when $t\geq 1$, while $|h|_{-\beta}$ appears from $\eta_{h,y}(1)$, and with estimates coming from Lemma \ref{lem000}.

\vspace{0.5cm}

\underline{Proof of Propositions \ref{lem3}}

Using equation \eqref{umarkov} and Lemma \ref{lem00}, for any $t\geq 1$ we have
$$|Du(t,x).h|\leq C\parallel \phi \parallel_{\infty}e^{-c(t-1)}\E[(1+|X(1,x)|^2)|\eta^{h,x}(1)|]\leq C\parallel \phi \parallel_{\infty}e^{-c(t-1)}(1+|x|^2)|h|_{-\beta},$$
where the last estimate comes from Lemmas \ref{lem1} and \ref{lem000}.

Combining this estimate with the result of Lemma \ref{lem000}, which gives an estimate for $t\leq 1$, we obtain \eqref{decreaseDU2}.

For the second order derivatives, Lemma \ref{lem000} gives an estimate for $t\leq 1$, and for $t\geq 1$ we use \eqref{umarkov} to see that
\begin{align*}
D^2u(t,x).(h,k)&=\E[D^2[u(t-1,X(1,x))].(h,k)]\\
&=\E D^2u(t-1,X(1,x)).(\eta^{h,x}(1),\eta^{k,x}(1))+\E Du(t-1,X(1,x)).\zeta^{h,k,x}(1).
\end{align*}
Using Lemma \ref{lem00}, we get an exponential decreasing; thanks to Lemma \ref{lem1} and to the estimates in the proof of Lemma \ref{lem000} at time $1$, we obtain
$$|D^2u(t,x).(h,k)|\leq \parallel \phi \parallel_{\infty}e^{-c(t-1)}(1+|x|^2)|h|_{-\beta}|k|_{-\gamma}.$$
Then \eqref{decreaseDU3} easily follows.

\qed

\section{Proof of some estimates}\label{SectProofsUsefulestim}

We give the detailed proofs of some estimates on the processes $(X^{h}(t))_{t\in \mathbb{R}^+}$ and $(X_{k}^{h})_{k\in \mathbb{N}}$, given in Section \ref{sectusefulestim}.

We omit the reference to the parameter of the spatial discretization $h\in(0,1)$, but it is clear from the proofs that the constants are uniform with respect to $h$.

\underline{Proof of Lemma \ref{lem1}}
If we define $Z(t)=X(t)-W^A(t)$, we have $Z(0)=X(0)=x$, and
$$\frac{dZ(t)}{dt}=AZ(t)+F(X(t)),$$
and by Proposition \ref{propodiss}
\begin{align*}
\frac{1}{2}\frac{d|Z(t)|^2}{dt}&=\langle AZ(t)+F(X(t)),Z(t)\rangle \\
&=\langle AZ(t)+F(Z(t)),Z(t)\rangle +\langle F(X(t))-F(Z(t)),Z(t)\rangle \\
&\leq -c|Z(t)|^2+C+\|F\|_{\infty}|Z(t)|\\
&\leq -c'|Z(t)|^2+C',
\end{align*}
for some new constants $c',C'$.

Then almost surely we have for any $t\geq 0$
$$|Z(t)|\leq C(1+|x|).$$
Thanks to \eqref{momWBp}, the conclusion easily follows.
\qed

\underline{Proof of Lemma \ref{lem2}}
As in the proof of Lemma \ref{lem1} above, we introduce $Z_m=X_m-w_m$, where the process $(w_m)$ is the numerical approximation of $W^A$ with the numerical scheme \eqref{defYk} - with $F=0$; it is defined by
$$w_{m+1}=S_{\tau,h} w_m+\sqrt{\tau}S_{\tau,h}\chi_{m+1}.$$
Using Theorem $3.2$ of \cite{prin2}, giving the strong order $1/4$ for the numerical scheme - when the initial condition is $0$, with no nonlinear coefficient, with a constant diffusion term and under the assumptions made here - we obtain the following estimate: for any $p\geq 1$, $\tau_0>0$ and $0<r<1/2$ there exists $C>0$ such that for any $0<\tau\leq \tau_0$ and $m\geq 0$
\begin{equation}\label{estimnoisenum}
\E|w_{m}-W^A(m\tau)|^{2p}\leq C\tau^{(1/2-r)p}.
\end{equation}
Thanks to \eqref{momWBp} and \eqref{estimnoisenum}, we get that for any $\tau_0>0$, there exists $C>0$ such that for $0<\tau\leq \tau_0$ and $m\geq 0$
\begin{equation}\label{estimunifwm}
\E|w_{m}|^{2}\leq C.
\end{equation}
Now $Z_m$ defined above satisfies $Z_0=X_0=x$ and
$$Z_{m+1}=S_{\tau,h} Z_m+\tau S_{\tau,h} F(X_m);$$
since, for $h\in(0,1)$, $|S_{\tau,h}|_{\mathcal{L}(H)}\leq \frac{1}{1+\lambda_0\tau}$, we obtain the almost sure estimates
$$|Z_{m+1}|\leq \frac{1}{1+\lambda_0\tau}|Z_m|+C\tau$$
and
$$|Z_m|\leq C(1+|x|).$$
Thanks to \eqref{estimunifwm}, we therefore obtain the result.
\qed

\end{appendix}

\end{document}